\documentclass[preprint,12pt]{elsarticle}
\usepackage{bbding}
\usepackage{dcolumn}
\usepackage{graphicx}
\usepackage{amsmath}
\usepackage{amsfonts}
\usepackage{amssymb}
\usepackage{psfrag}
\usepackage{epsfig}
\usepackage{wrapfig}
\usepackage{subfigure}
\usepackage{makeidx}
\usepackage{bm}
\usepackage{epsf}
\usepackage{epsfig}
\usepackage{leftidx}
\usepackage{float}
\usepackage{extarrows}
\usepackage{color}
\usepackage{cases}
\newdefinition{rmk}{Remark}
\begin{document}

\title{A Unified Gas-kinetic Scheme for Continuum and Rarefied Flows IV: full Boltzmann and Model Equations}


\author[HKUST]{Chang Liu}
\ead{cliuaa@ust.hk}

\author[HKUST]{Kun Xu\corref{cor1}}
\ead{makxu@ust.hk}

\author[imech]{Quanhua Sun}
\ead{qsun@imech.ac.cn}

\author[pku]{Qingdong Cai}
\ead{caiqd@mech.pku.edu.cn}

\address[HKUST]{Department of Mathematics and Department of Mechanical and Aerospace Engineering, Hong Kong University of Science and Technology, Clear Water Bay, Kowloon, Hong Kong}
\address[imech]{State Key Laboratory of High-temperature Gas Dynamics, Institute of Mechanics, Chinese Academy of Sciences,
No. 15 Beisihuan Xi Rd, Beijing 100190, China}
\address[pku]{Department of Mechanics and Aerospace Engineering, College of Engineering, Peking University, Beijing 100871, China}
\cortext[cor1]{Corresponding author}

\begin{abstract}

Fluid dynamic equations are valid in their respective modeling scales, such as the particle mean free path scale of the Boltzmann equation
and the hydrodynamic scale of the Navier-Stokes (NS) equations.
With a variation of the modeling scales, theoretically there should have a
continuous spectrum of fluid dynamic equations.
Even though the Boltzmann equation is claimed to be valid in all scales,
many Boltzmann solvers, including direct simulation Monte Carlo method, require the cell resolution to the particle mean free path scale.
Therefore, they are still single scale methods.
In order to study multiscale flow evolution efficiently,
the dynamics in the computational fluid has to be changed with the scales.
A direct modeling of flow physics with a changeable scale may become an appropriate approach.
The unified gas-kinetic scheme (UGKS) is a direct modeling method in the mesh size scale,
and its underlying flow physics depends on the resolution of the
cell size relative to the particle mean free path.
The cell size of UGKS is not limited by the particle mean free path.
With the variation of the ratio between the numerical cell size and local particle mean free path, the UGKS recovers the flow dynamics from the
particle transport and collision in the kinetic scale to the wave propagation in the hydrodynamic scale.
 The previous UGKS is mostly constructed from the evolution solution of kinetic model equations.
Even though the UGKS is very accurate and effective in the low transition and continuum flow regimes with
the time step being much larger than the particle mean free time, it is still necessary to develop more accurate flow solver in the rarefied regime,
where the time step is comparable  with the local particle mean free time.
In such a scale, there is dynamic difference from the full Boltzmann collision term and the model equations.
This work is about the further development of the UGKS with the implementation of the full Boltzmann collision term in the region where it is needed.
The central ingredient of the UGKS is the coupled treatment of particle transport and collision in the flux evaluation across a cell interface, where
a continuous flow dynamics from  kinetic to hydrodynamic scales is modeled.
The newly developed UGKS has the asymptotic preserving (AP) property of recovering the NS solutions in the continuum flow regime, and the full Boltzmann solution in the rarefied regime.
In the mostly unexplored transition regime, the UGKS itself provides a valuable tool for the flow study in this regime.
The mathematical properties of the scheme, such as stability, accuracy, and the asymptotic preserving, will be analyzed in this paper as well.

\end{abstract}
\maketitle

\section{Introduction}

The flow regime is categorized according to the Knudsen number $\mbox{Kn}$, which is defined as the ratio of the molecular mean free path
to a characteristic length scale. The value of the Knudsen number determines the validity of different approaches in the description of gas flow.
The whole flow regime is qualitatively divided into continuum ($\mbox{Kn}<0.001$), transitional ($0.001<\mbox{Kn}<10$), and free molecular regimes ($\mbox{Kn}>10$).
Numerically, all solutions obtained are in the mesh size scale.
As a variation of mesh size, different dynamics, such as particle free transport and wave propagation, should appear automatically in an idealized numerical method.
An appropriate Knudsen number for a numerical scheme may be the cell Knudsen number, which can be defined as the particle mean free path over the numerical cell size.
Due to the relative change of the cell's Knudsen number, different flow dynamics should be captured.
The aim of the unified gas-kinetic
scheme (UGKS) is to capture such a flow evolution in terms of the cell's Knudsen number.

A UGKS based on the kinetic BGK and Shakhov models has been developed
in the past \cite{xu2010,Huangxuyu2,Huangxuyu3,Songze, Shaliu,liu-d}.
The unified scheme is a multi-scale method with coupled particle transport and collision in its numerical flux modeling.
A time evolution solution of the kinetic model equation has been used to construct the flux transport across a cell interface.
This time evolution solution covers the flow physics from the kinetic scale particle free transport to the hydrodynamic scale wave propagation,
and the weight between these two limiting solutions depends on the ratio of time step to the local particle mean free time $\tau$.
As a result, both kinetic and hydrodynamic solutions can be automatically obtained in a unified way.
In the continuum flow regime, due to the un-splitting treatment of particle transport and collision,
and automatic recovering of the Chapman-Enskog gas distribution function for the flux evaluation,
the viscous effect can be captured by UGKS without the constraints on the
cell size and time step, i.e. cell size $\Delta x \leq l_{mfp}$ and time step $\Delta t \leq \tau$,
which are required by many Boltzmann solvers for NS solutions.
Therefore, the UGKS is more efficient than the single scale-based method, such as the direct simulation Monte Carlo (DSMC) method,
in the low transition and continuum flow regime with the mesh size on the order of tens or hundreds of particle mean free path.

If a uniform time step is used for the continuum flow computation, and at the same time the highly non-equilibrium shock structure needs to be well-resolved, the accuracy requirement enforces the UGKS to use a small cell
size at the shock region, which consequently limits the use of the overall large time step. Therefore, under this kind of situation
the UGKS will not be efficient, even in the continuum flow regime.
However, in the continuum and near continuum flow regime, we may not need to resolve the highly non-equilibrium shock layer.
In most cases, when a large cell size and time step are used, the UGKS
will become a shock capturing scheme, where the shock structure with a few mesh points is an artificial one.
If the shock capturing property can be accepted in the continuum flow regime,
the efficiency of the UGKS can be kept.
For the NS viscous boundary layer solution at high Reynolds number, the boundary layer thickness is much larger than the particle mean free path.
The boundary layer can be captured accurately by UGKS with a few mesh points without resolving the mean free path scale.
For a multiscale unsteady flow problem with both continuum and rarefied regimes,
such as the gas expansion from a nozzle into vacuum,
the advantage of the UGKS is obvious.
With a unform time step $\Delta t$,
which can be much larger than the local particle mean free time inside the nozzle for the NS solution,
and can be much smaller than the particle mean free time outside the nozzle for the Boltzmann solution \cite{Songze},
a smooth dynamic transition in the gas expansion process can be accurately captured by UGKS.
Moreover, for steady flow simulation, a non-uniform local time step and implicit discretization can be obtained in UGKS to improve its efficiency without loss of its accuracy.
Overall, the UGKS provides a general framework to construct multi-scale method for transport process, such as
the recent extension to the radiative transfer \cite{luc,sun-wj,sun-wj2}.

The previous development of UGKS is based on the kinetic model equations which approximate the full Boltzmann equation.
The Boltzmann equation is a modeling equation in the kinetic scale, i.e., the scale to identify the distinguishable process of  particle transport and collision.
In such a kinetic scale, there is difference in dynamics between the kinetic collision model and the full Boltzmann collision term.
One of the objectives of this paper is to quantitatively evaluate such a difference and to use the full Boltzmann collision term in the region where it is needed.
For the unified scheme, the ratio of the time step $\Delta t$ over the local particle mean free time $\tau$ can be varied significantly from the kinetic scale regime $\Delta t \leq \tau$
to the hydrodynamic scale regime $\Delta t \gg \tau$.
In the regime with $\Delta t \gtrsim \tau$,  the solution difference from the full Boltzmann collision term and the
kinetic model equation diminishes. Based on this observation, a UGKS  with the implementation of a hybrid particle collision terms can be constructed.
The use of both full Boltzmann collision term and the kinetic model equation is close to the idea of penalty method \cite{filbet2010class}, but with a distinguishable consideration in the
design of UGKS.
For the UGKS, the full Boltzmann collision term is only used in the local kinetic regime, where the time step is less than a critical time interval $t_c$ which is related to the deviation of a distribution
function from equilibrium.
When a particle encounters multiple collisions within a time step,
the evolution of a gas distribution function will not be sensitive to individual particle collision at all. Therefore, the kinetic model equations can be faithfully used in the
regime beyond the kinetic one, such as those regions with $\Delta t \geq t_c$.
The current modeling method can give accurate Boltzmann solution in the rarefied regime and exact NS solution in the continuum flow regime.
We believe that the UGKS also presents accurate solution in the whole transition regime, where a continuous spectrum of
gas dynamics from the kinetic to the hydrodynamic scale is recovered \cite{xu-book}.
For example, based on UGKS, the cavity flow solution in the whole transition regime from $\mbox{Kn} =10$ to $\mbox{Kn}=10^{-4}$ is obtained and compared with the DSMC and NS solutions.
At the same time, the solution differences between UGKS and NS in the near continuum regime, such as from Reynolds numbers $5$ to $50$, are explicitly presented.

The UGKS in this paper is equipped with the full Boltzmann collision operator in the kinetic regime with small ratio of $\Delta t/\tau$.
The computational cost of solving the five-fold quadratic integral operator is so large that fast algorithms are required.
In the literature, a lot of  effective methods have been developed for solving the Boltzmann equation.
A pioneering work has been done by Bobylev in 1988 using Fourier transform techniques in the analysis of the Boltzmann equation for Maxwell molecules \cite{bobylevanalysis}.
Later,  Bobylev and Rjasanow  developed a numerical method to solve the collision operator for Maxwell molecules with a computational cost of the order $O(N^4)$
in 1996 \cite{bobylevmaxwell}, where $N$ denotes the number of discrete velocity points in each dimension, and extended the method to hard sphere molecules in 1999 \cite{bobylev1999hs}.
In 2000, Pareschi and Russo  developed an algorithm to solve the collision operator for the variable hard sphere (VHS) model with a computational cost of $O(N^6)$ \cite{pareschiVHS}.
In 2002, Ibragimov and Rjasanow solved the collision operator on a uniform grid with computational cost $O(N^6\log(N))$
for general model and achieved an accuracy of the order $O(N^{-2})$ \cite{ibragimov2002numerical}.
Later, in 2006 by means of the Carleman-like representation and the fast Fourier transform, Mouhot and Pareschi  developed a fast spectral algorithm with computational cost $O(M^2N^3\log N)$ \cite{mouhot2006fast}, where $M$ denotes the number of grid points in the discretization of the unit sphere.
Gamba and Tharkabhushanam extended the fast spectral method to the non-elastic collision by spectral-Lagrangian method \cite{gamba2009spectral}.
The fast spectral method has been applied to space non-homogeneous problems in two-dimensional velocity space
\cite{filbet2006solving, filbet2012deterministic, filbet2003high} as well as to the quantum collision operators \cite{filbet2012numerical, hu2012fast}.
Algorithms other than the spectral method have also been developed for solving the Boltzmann equation, such as
the finite element methods developed by the kinetic group in Kyoto \cite{ohwada2004kinetic, tcheremissine2005direct, tcheremissine2006solution},
and the discontinuous Galerkin method \cite{alekseenko2014deterministic, armandonumerical,majorana2011}.
For the test cases in this paper, we adopt the fast spectral method of  a recent paper by Wu et al. \cite{wu}.
The UGKS targets to obtain accurate solutions in all flow regimes.
In the rarefied regime, the Boltzmann collision operator plays an important role in capturing the peculiar highly non-equilibrium gas distribution function.
In the hydrodynamic regime $\Delta t \gg \tau$, due to the intensive particle collisions the full Boltzmann collision term and many kinetic model equations
can present identical Chapman-Enskog gas distribution function.
Here the use of the full Boltzmann collision term is not necessary at all, and the scheme based on the model equations can become more efficient than that based on the full Boltzmann collision term
without sacrificing the accuracy.

This paper is organized in the following. The full Boltzmann equation and the kinetic model equations will be introduced in Section 2.
Section 3 is about the numerical experiments on time evolution of gas distribution functions from full Boltzmann collision term and kinetic model equation.
Based on this observation, a unified scheme is proposed in Section 4. Section 5 is about stability, accuracy, and asymptotic preserving analysis of UGKS.
Numerical experiments are presented in Section 6.
The last section is conclusion.

\section{Boltzmann equation and kinetic model equations}

The Boltzmann equation is a fundamental equation which statistically models the gas dynamics in the kinetic scale, i.e. the scales of particle mean free time and mean free path.
In this work, we focus on monatomic gas with binary elastic collisions.
For space variable $\mathbf{x} \in \mathcal{R}^3$, particle velocity $\mathbf{u}=(u,v,w)^t \in \mathcal{R}^3$, the corresponding Boltzmann equation reads:
\begin{equation}\label{Bzm}
\frac{\partial f}{\partial t}+ \mathbf{u}\cdot\nabla_\mathbf{x} f = Q(f,f),
\end{equation}
where $f: =f(\mathbf{x},t,\mathbf{u})$ is the time-dependent particles distribution function in the phase space.
The collision term $Q(f,f)$ is a quadratic collision operator,
\begin{equation}\label{Col}
Q(f,f) =\int_{\mathcal{R}^3}\int_{\mathcal{S}^2} (f'_*f'-f_*f) |\mathbf{u_r}|  \sigma d\Omega d\mathbf{u}_*.
\end{equation}
Here the short hand notation $f'_*=f(\mathbf{x},t,\mathbf{u}'_*)$ is used, similarly for $f'$ and $f_*$.
Based on conservation of momentum and energy, the pre-collision particle velocities  $\mathbf{u}$, $\mathbf{u}_*$ and the corresponding post-collision velocities $\mathbf{u}'$, $\mathbf{u}'_*$ satisfy the follow relations
\begin{equation}
\begin{aligned}
  \mathbf{u}'&=\frac{\mathbf{u}+\mathbf{u}_*}{2}+\frac{|\mathbf{u}-\mathbf{u}_*|}{2}\Omega=\mathbf{u}+\frac{|\mathbf{u_r}|\Omega-\mathbf{u_r}}{2}, \\
  \mathbf{u}'_*&=\frac{\mathbf{u}+\mathbf{u}_*}{2}-\frac{|\mathbf{u}-\mathbf{u}_*|}{2}\Omega=\mathbf{u}_*-\frac{|\mathbf{u_r}|\Omega-\mathbf{u_r}}{2},
\end{aligned}
 \end{equation}
 where $\mathbf{u_r}=\mathbf{u}-\mathbf{u}_*$ is the relative pre-collision velocity and $\Omega$ is a unit vector in $\mathcal{S}^2$ along the relative post-collision velocity $\mathbf{u'}-\mathbf{u'}_*$.

The differential cross section $\sigma$ measures the probability of collision which depends on the strength of relative velocity and deflection angle between pre-collision and post-collision velocities.
In this paper, three collision models are used in calculating the cross section, namely the hard sphere (HS) model, variable hard sphere (VHS) model, and an anisotropic collision model proposed by Mouhot and Pareschi \cite{mouhot2006fast}.
The differential cross section for the hard sphere molecules can be written down as
\begin{equation}\label{hscross}
  \sigma=\frac{5}{64\sqrt{\pi}}\frac{\sqrt{mkT_{ref}}}{\mu_{ref}},
\end{equation}
and the dynamic viscosity coefficient satisfies
\begin{equation}\label{hsviscous}
  \mu=\frac{5}{16}\sqrt{\frac{2\pi kT}{m}}\rho\ell,
\end{equation}
where $k$ is the Boltzmann constant, T stands for temperature, and $\ell$ denotes the mean free path at equilibrium state.
The viscosity coefficient of hard sphere molecules is proportional to $T^{0.5}$.
When we simulate the argon gas with viscosity $\mu\varpropto T^{0.81}$, we need to employ the VHS model whose cross section is
\begin{equation}\label{vhscross}
  \sigma=\frac{15m u_r}{64\sqrt{\pi}\Gamma\left(\frac92-\omega\right)\mu_{ref}}\left(\frac{4kT_{ref}}{u_r^2}\right)^\omega,
\end{equation}
and the corresponding dynamic viscosity coefficient follows,
\begin{equation}\label{vhsvisvous}
  \mu=\frac{15\rho\ell \sqrt{2\pi kT}}{2(7-2\omega)(5-2\omega)\sqrt{m}}.
\end{equation}
The VHS model or a simple inverse power law model is a phenomenological model. In reality the potential between monatomic gas is better described by the Lennard-Jones (L-J) potential. The L-J potential for argon gas is
\begin{equation}\label{LJpotential}
  \phi(r)=4\epsilon\left[\left(\frac{d_{LJ}}{r}\right)^{12}-\left(\frac{d_{LJ}}{r}\right)^{6}\right],
\end{equation}
with potential depth $\epsilon=119.18 k$, and $d_{LJ}=3.42\times10^{-10}m$.
A generalized anisotropic collision model suggested by Mouhot and Pareschi \cite{mouhot2006fast} is used in this paper to recover the L-J potential, and the differential cross section of which can be written as
\begin{equation}\label{anisotropic}
  \sigma=\sum_j C'_{\alpha_j} \sin^{\alpha_j-1}\left(\frac{\theta}{2}\right)|u_r|^{\alpha_j-1},
\end{equation}
where $C'_{\alpha_j}$ is a constant.
The dynamic viscosity of argon with L-J potential is numerically fitted by Wu \emph{et al.} \cite{wu} as
\begin{equation}\label{LJviscous}
  \mu=\frac{5\sqrt{\pi m kT}}{8d^2_{LJ}\sum_jb_i(kT/\epsilon)^{\frac{\alpha_j-1}{2}}},
\end{equation}
with the coefficient $\alpha_1=0.2$, $\alpha_2=0.1$, $\alpha_3=0$ and $b_1=407.4$, $b_2=-881.9$, $b_3=414.4$.
The fitted viscosity coefficient corresponds to the following differential cross section,
\begin{equation}\label{LJ}
\sigma_{LJ}=\frac{d_{LJ}^2}{32\pi}\sum_{j=1}^3\frac{(m/4\epsilon)^{(\alpha_j-1)/2}b_j}{\Gamma\left(\frac{3+\alpha_j}{2}\right)}\sin^{\alpha_j-1}\left(\frac{\theta}{2}\right)|u_r|^{\alpha_j-1}.
\end{equation}
For the HS, VHS and generalized anisotropic collision models described above,
the corresponding Boltzmann collision operator as well as the local collision frequency can be solved by the FFT based fast spectral method \cite{wu}.

Due to the stiffness nature of the Boltzmann collision operator,
the implicit treatment is preferred to stabilize the scheme, especially in the continuum regime with intensive particle collisions.
 The nonlinear Boltzmann collision operator is not convenient to be treated implicitly.
 Fortunately, in the regime with the scale of multiple collisions, the accumulating effect from the full Boltzmann equation and the kinetic model equations present the same result on the evolution of
 a gas distribution function, see section 3 for the numerical experiments.
 { Therefore, the collision process is modeled using an explicit Boltzmann collision operator and an implicit Shakhov model.
 The domains of the explicit and implicit discretization are determined from the comparison between the time step $\Delta t$ and a critical time interval $t_c$,
 which is related to the deviation of a gas distribution function from local equilibrium.
 The stability analysis in Section 5 shows that the time step will not be constrained by the stiffness of the Boltzmann collision term in the region beyond the kinetic scale.}

In the literature, many kinetic models are proposed, among which the ES-BGK \cite{holway} and Shakhov \cite{shakhov1968generalization} are two popular ones.
These two models can be combined as a generalized model \cite{chen-combined} with the introduction of one more degree of freedom.
In this paper, we will use the full Boltzmann and Shakhov model to construct UGKS.
In general, the kinetic model can be written as
\begin{equation}\label{mdl}
  f_t + \mathbf{u}\cdot\nabla_\mathbf{x} f = S(f),
\end{equation}
with a relaxation-type source term $S(f)$,
\begin{equation}
  S(f)=\frac{\tilde{M}(f)-f}{\tau_{r}}\nonumber,
\end{equation}
where $\tau_r$ denotes the relaxation time that distinguishes from the local mean free time obtained from the Boltzmann collision term.
For Shakhov model, $\tilde{M}(f)$ is defined as
\begin{equation}\label{shakhov-M}
 \begin{aligned}
   \tilde{M}(f)&=M(f)+g^1(f),\\
   M(f)&=\rho\left(\frac{\lambda}{\pi}\right)^{\frac{3}{2}}\mathbf{e}^{-\lambda(\mathbf{u}-\mathbf{U})^2},\\
   g^1(f)&=M(f)(1-\text{Pr})\mathbf{c}\cdot
   \mathbf{q}\left(\frac{\mathbf{c}^2}{RT}-5\right)/(5pRT),
 \end{aligned}
\end{equation}
 where $\rho$, $\mathbf{U}$, $T$ and $\mathbf{q}$ are the corresponding density, velocity,
 temperature, and heat flux of $f$. Here $\mathbf{c}$ is the random velocity, Pr is the Prandtl number,
 and $R$ is specific gas constant. Thus $M(f)$ is the corresponding equilibrium state of $f$ and $\tilde{M}(f)$
 is the modified equilibrium state with Prandtl number correction. In continuum regime, the heat flux is of the order $\tau_r$, and therefore $g^1(f)$ has the same order of $\tau_r$.

{In the highly non-equilibrium rarefied regime,
the kinetic model equation may not be able to fully describe
the time evolution of the distribution as accurate as the Boltzmann equation.
In the hydrodynamic regime
where the velocity distribution can be written as an asymptotic expansion with respect to Kn,
i.e. $f=f^0+ \mbox{Kn} f^1+O({\mbox{Kn}}^2)$,
both the Boltzmann equation and Shakhov model give the same solution up to the order $O(\text{Kn})$,
which recovers the NS equations \cite{xu-book}.
For a given viscosity coefficient, the cross section in Boltzmann equation is calculated from Eq.\eqref{hsviscous}, \eqref{vhsvisvous}, and \eqref{LJviscous},
and the relaxation parameter $\tau_r$ in Shakhov equation is defined as the ratio of dynamic viscosity coefficient to pressure.
Numerically, in the highly non-equilibrium regime with the scales of particle mean free path and mean collision time,
the Boltzmann collision operator is essential to capture the flow physics.
In the hydrodynamic scale where the local mean free time is much smaller than the characteristic time scale,
the numerical flux at cell interface plays an important role in capturing the NS solution other than Boltzmann collision term,
and both transport and collision process have to be taken into account in flux modeling \cite{chen-xu}.}

In this paper, the flux in the UGKS is calculated from the time evolving solution of the Shakhov model.
The collision term is treated with a hybridization of the full Boltzmann collision operator and the Shakhov model,
and the switching function depends on the ratio between local particle mean free time and numerical time step.
Section 3 presents the numerical experiments to validate such a switching function modeling.
Physically, with the change of the ratio between the time step and the particle mean free time, a continuous variation of flow physics in different regime will emerge automatically.

\section{Distribution function evolution from the full Boltzmann and Shakhov collision terms}

In order to compare the collision effect from the full Boltzmann collision term and the Shakhov model, we study the homogeneous flow relaxation process.
Theoretical, based on Wild's analysis \cite{wild1951boltzmann},
 starting from a non-equilibrium initial condition, the solutions of Boltzmann and Shakhov equations shall
 get close when time becomes larger than the local mean free time,
 which is consistent with a numerical experiment recently done by Sun \emph{et al.} \cite{sun}.
 Here we are going to quantitatively evaluate the differences between full Boltzmann solution and Shakhov solution in specific cases.
 These experiments and observations are the basis for the construction of UGKS.
As a remark, the development of numerical method can be based on the physical laws and the experimental observation, which is similar to the derivation of the fluid dynamic equations.

Three kinds of relaxation problems are considered.
The first one is the evolution of an anisotropic Maxwellian distribution.
Specifically, the distribution for each velocity component is Maxwellian,
but has different temperature in different directions.
The second one is double half-normal distribution,
where a full distribution is comprised of two half-normal distributions in one velocity space,
and is Maxwellian type in  other velocity directions.
This test is used to show the evolution of a discontinuous distribution function.
The third one is a tailored half-Maxwellian distribution,
which is similar to the second case except that the discontinuity is removed by adjusting the amplitudes of half distributions.
The third test is  rather general with a continuous distribution, but asymmetric.
In the previous study \cite{sun}, the solutions from the kinetic model equations are compared with DSMC solutions.
Here the comparison with the full Boltzmann solution will be presented.

The working gas is monatomic argon with viscosity coefficient $\mu \propto T^{0.81}$.
The collision model used is the VHS model, with viscosity of Eq.\eqref{vhsvisvous}
and differential cross section of Eq.\eqref{vhscross}. The relaxation parameter in Shakhov model are calculated by
\begin{equation}
\tau_{r}=\frac{\mu}{P}=\frac{30}{(7-2\omega)(5-2\omega)}\tau=1.65\tau\nonumber,
\end{equation}
where $\tau$ is the local mean free time at equilibrium state.

\subsection{Relaxation of anisotropic Maxwellian distribution}

The initial anisotropic Maxwellian distribution is specified as follows
$$  f_0(u,v,w)=\frac{\beta_1}{\sqrt{\pi}}\mathrm{e}^{-\beta_1^2u^2}\frac{\beta_2}{\sqrt{\pi}}\mathrm{e}^{-\beta_2^2v^2}\frac{\beta_3}{\sqrt{\pi}}\mathrm{e}^{-\beta_3^2w^2}, $$
where $\beta_i=\sqrt{m/2kT_i}$ for $i=1,2,3$.
Two cases are tested with the following initial conditions of $T_1=273K$, $T_2=373K$, $T_3 = 273K$,  and $T_1=273K$, $T_2=5460K$, $T_3=273K$.

\begin{figure}
\center
\includegraphics[width=6cm]{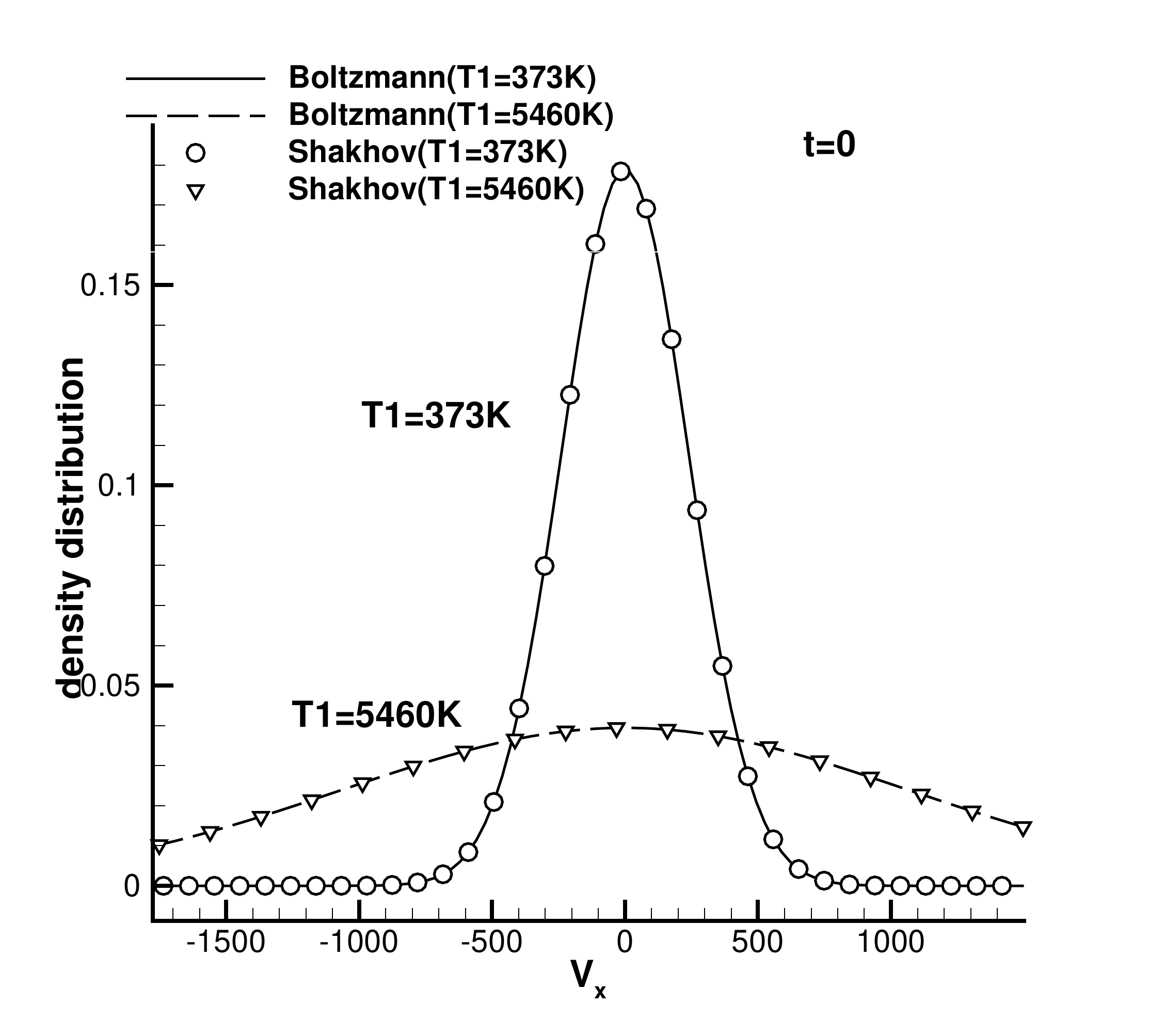}
\includegraphics[width=6cm]{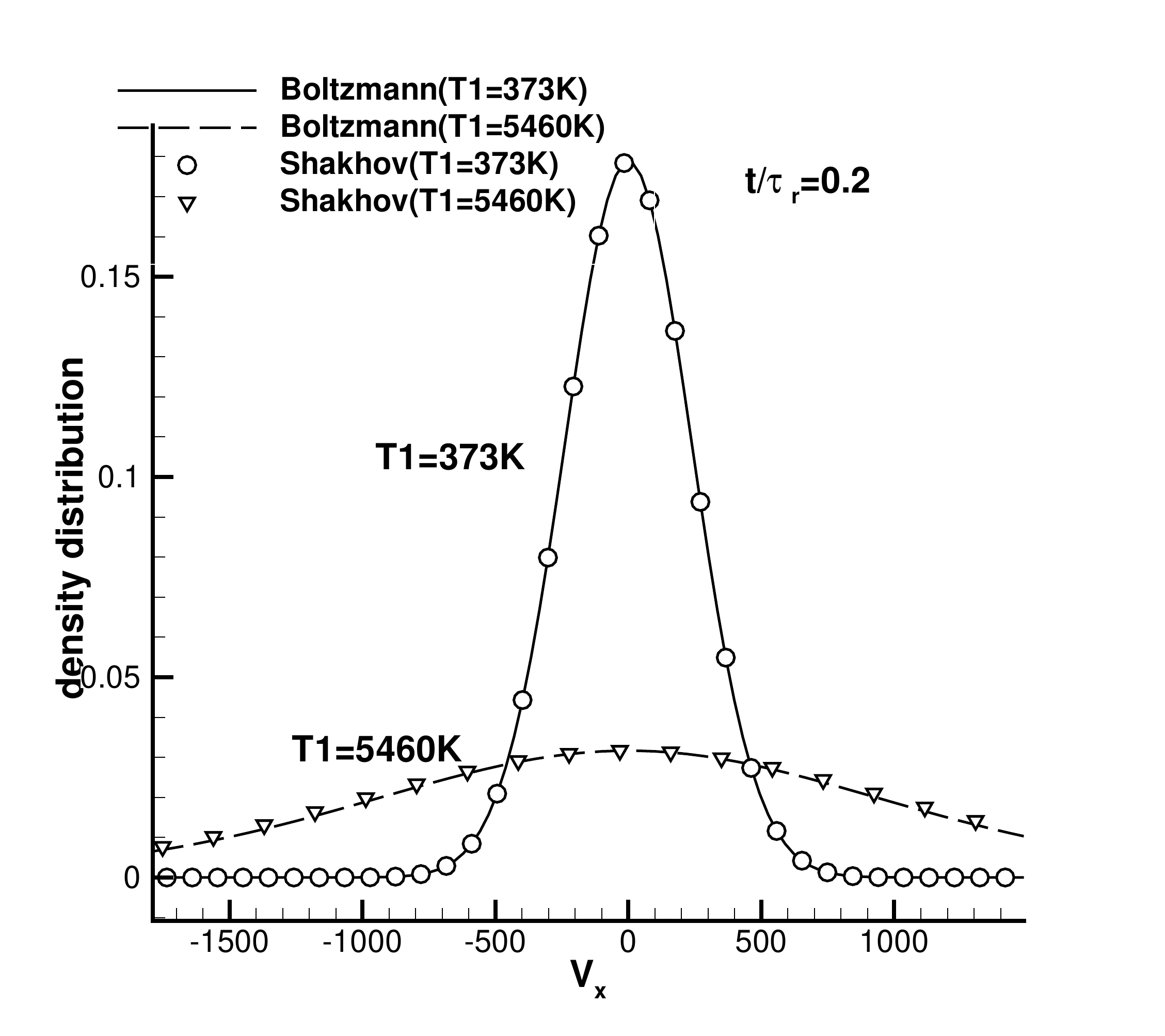}
\includegraphics[width=6cm]{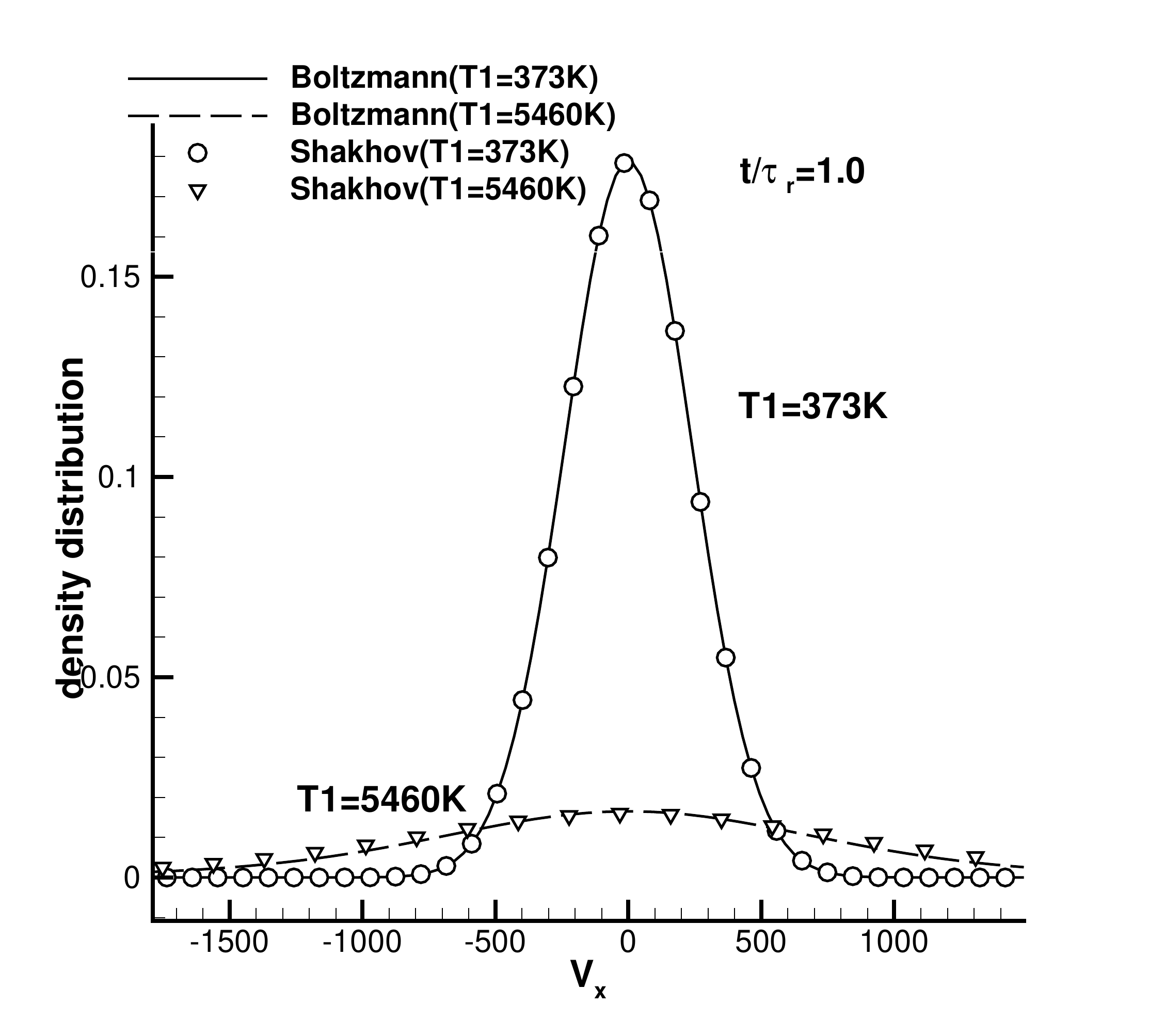}
\includegraphics[width=6cm]{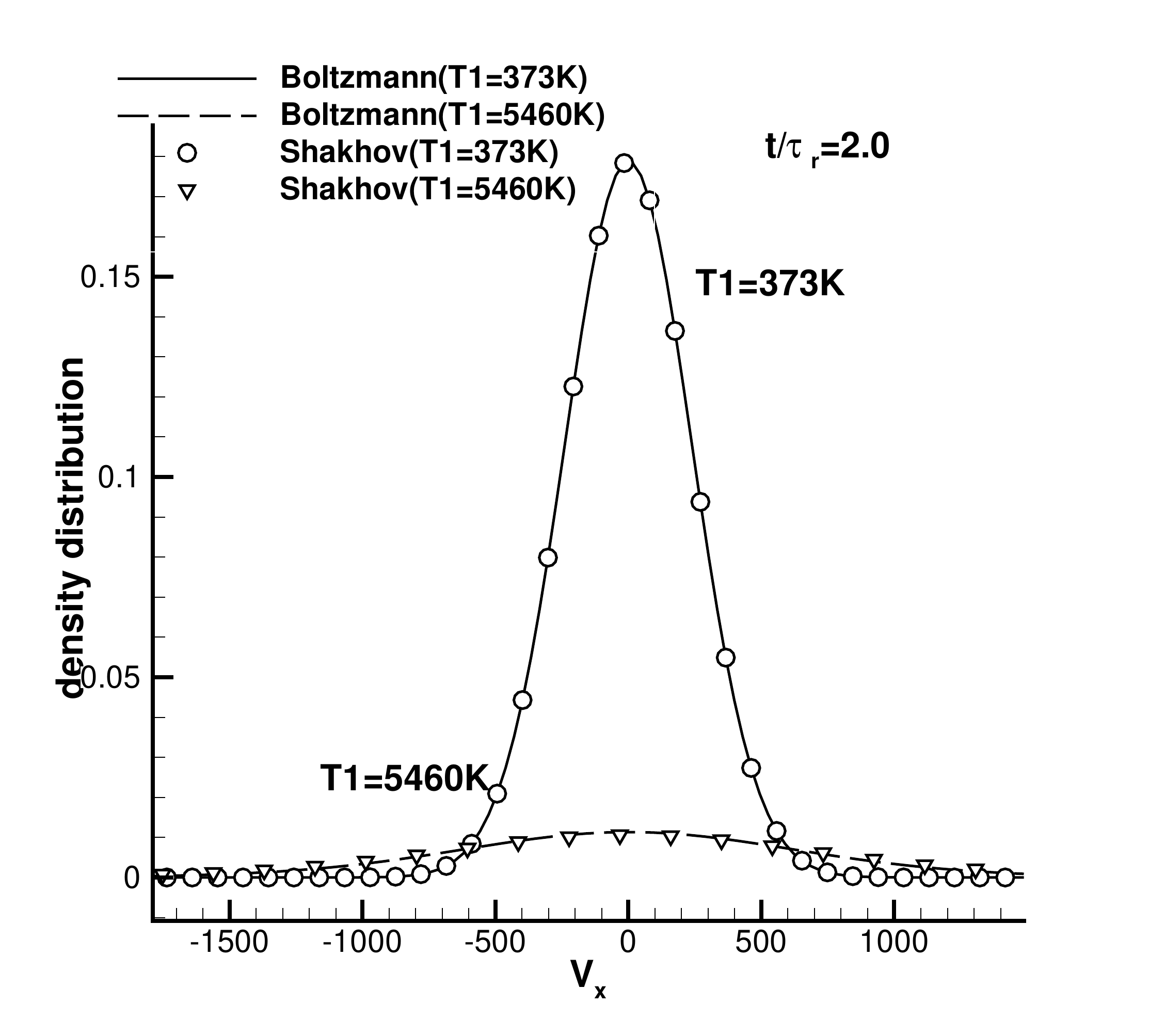}
\includegraphics[width=6cm]{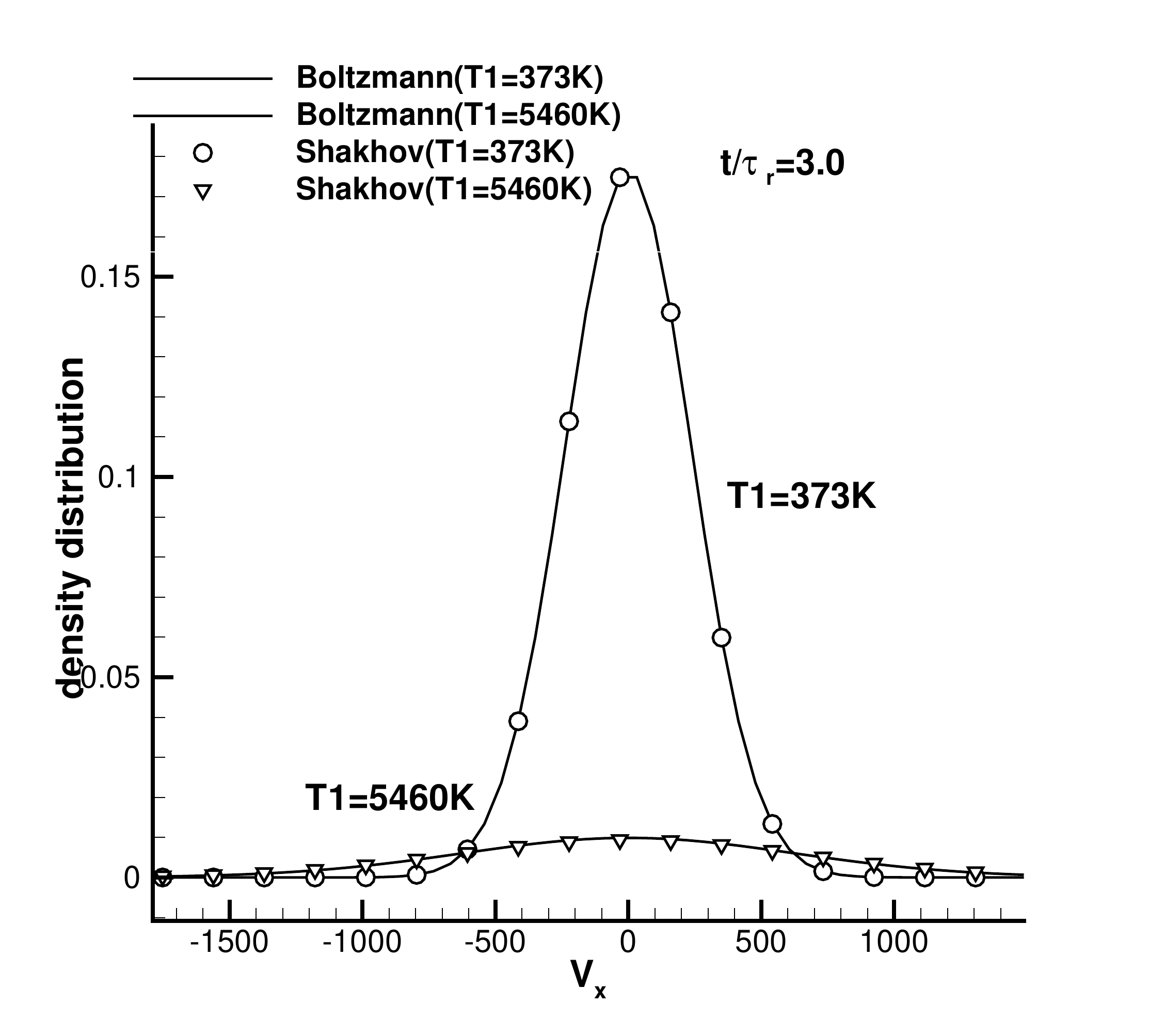}
\includegraphics[width=6cm]{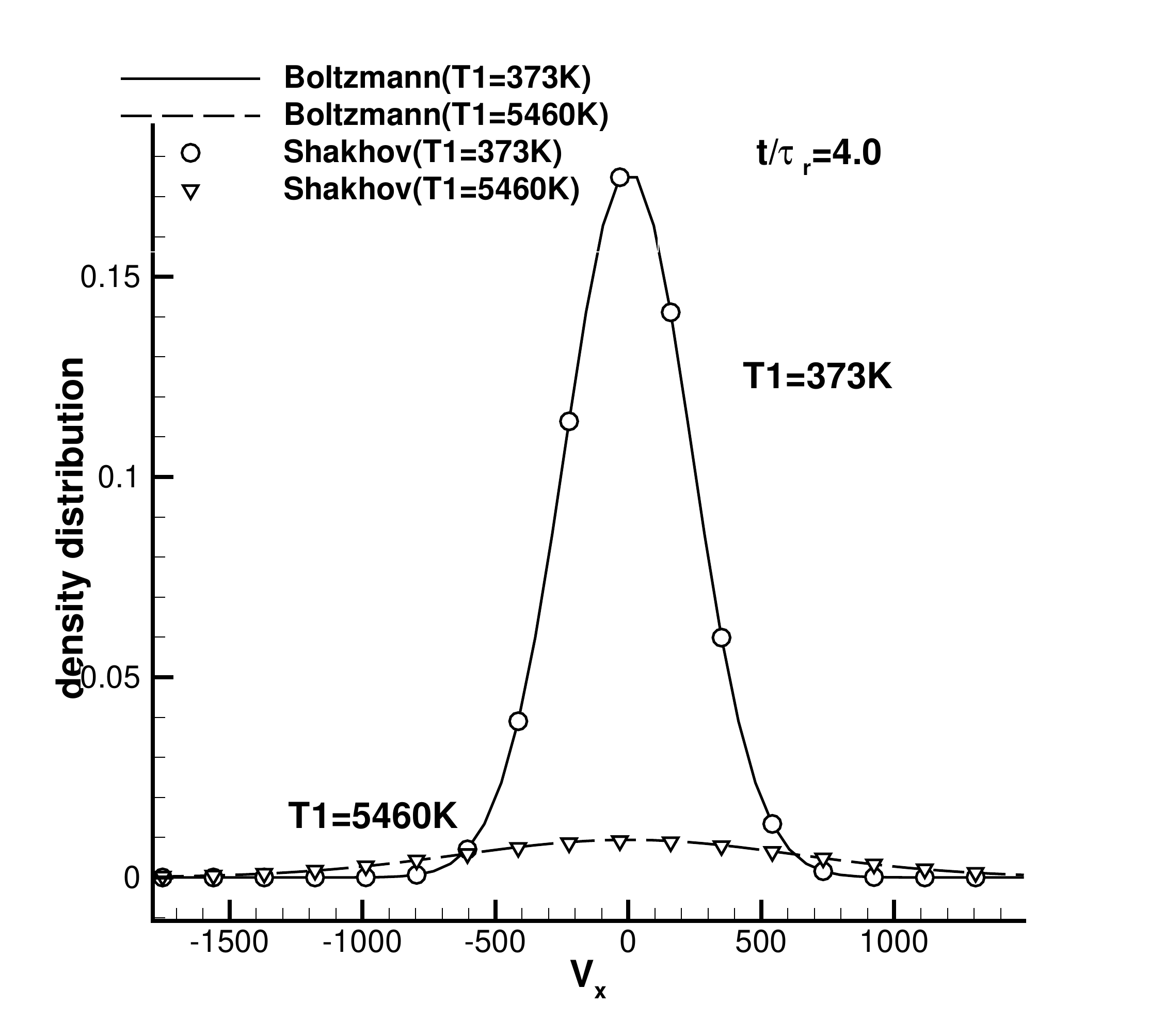}
\caption{Comparison of the x-component distribution functions $f(u,0,0)$ of the relaxation of anisotropic Maxwellian distribution.}
\label{fig:rel1}
\end{figure}

Fig.(\ref{fig:rel1}) shows the x-component distribution functions $f(u,0,0)$  at different output times
where the Shakhov model solution (symbols) and the full Boltzmann results (lines) are compared.
The solutions show that when
{$$t_1>\frac{2}{\rho}\int|f-M| dv\tau_r\approx0.2\tau_{r}$$}
for the first initial condition,
and
{$$t_2>\frac{2}{\rho}\int|f-M| dv\tau_r\approx2\tau_{r}$$}
for the second one,
two solutions agree with each other very well.
This test shows that the larger the temperature difference is, the longer it takes to get the same solution from different collision models.
Even with such a large temperature difference, from $273K$ to $5460K$ in different directions,
four particle collisions are enough to get indistinguishable solutions from the Shakhov model and the full Boltzmann collision term.
Even at $t \simeq 2\tau_r$, the two solutions are close to each other.

\subsection{Relaxation of double half-normal distribution}

Two cases with different initial conditions are tested, and both are related to the double half-normal distributions with discontinuities in the middle,
%
$$  f_0(u,v,w)=\left[\frac{\beta_1}{\sqrt{\pi}}\mathrm{e}^{\beta_1^2u^2}|_{u<0}+\frac{\beta_2}{\sqrt{\pi}}\mathrm{e}^{-\beta_2^2u^2}|_{u\geq0}\right]
  \frac{\beta_2}{\sqrt{\pi}}\mathrm{e}^{-\beta_2^2v^2}\frac{\beta_3}{\sqrt{\pi}}\mathrm{e}^{-\beta_3^2w^2}, $$
%
where $\beta_i =\sqrt{m/2kT_i}$, with $T_1=273K$, $T_2=373K$, $T_3 =273K$ for the first case, and  $T_1=273K$, $T_2=5460K$, $T_3=273K$ for the second case.

\begin{figure}
\center
\includegraphics[width=6cm]{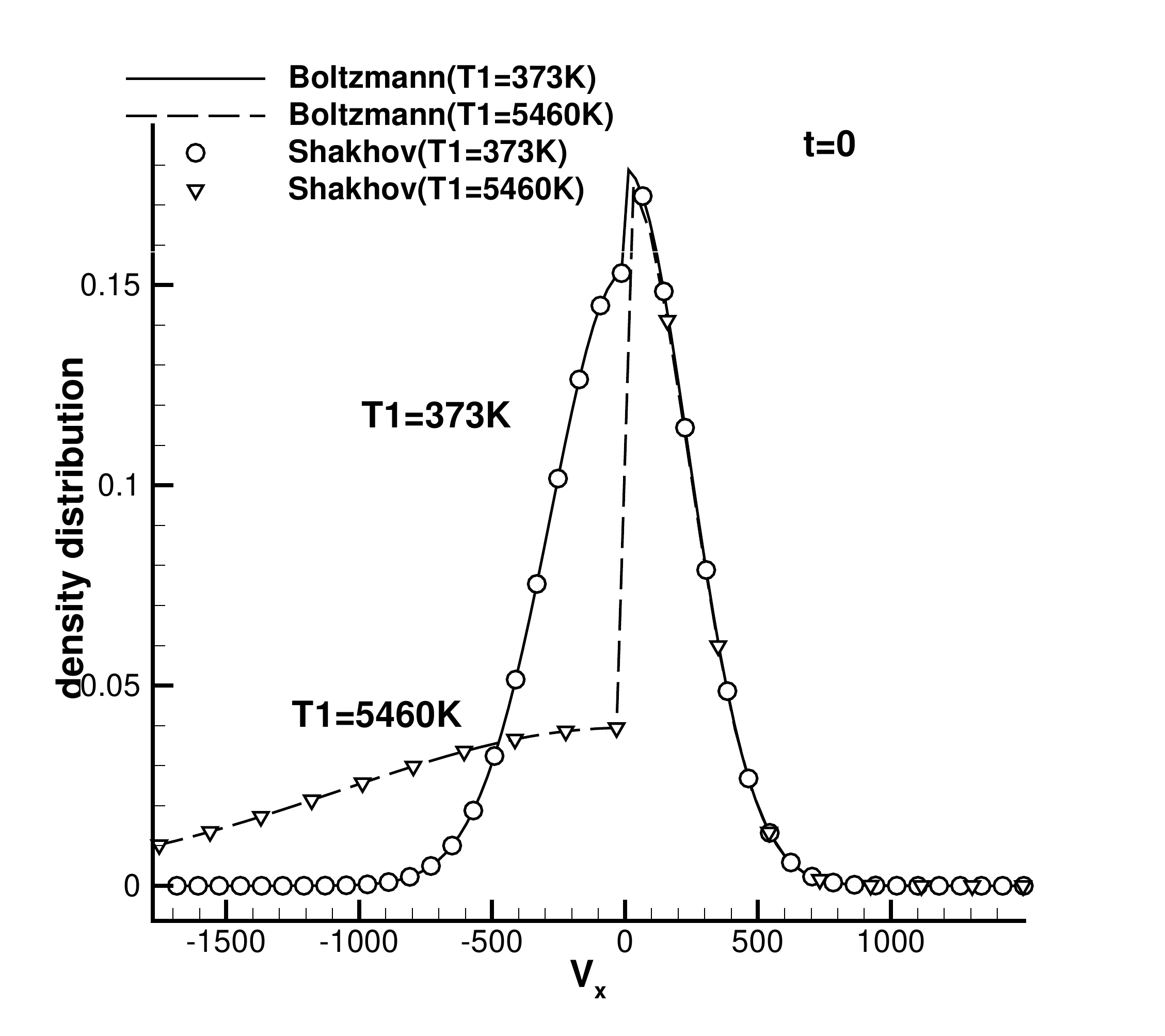}
\includegraphics[width=6cm]{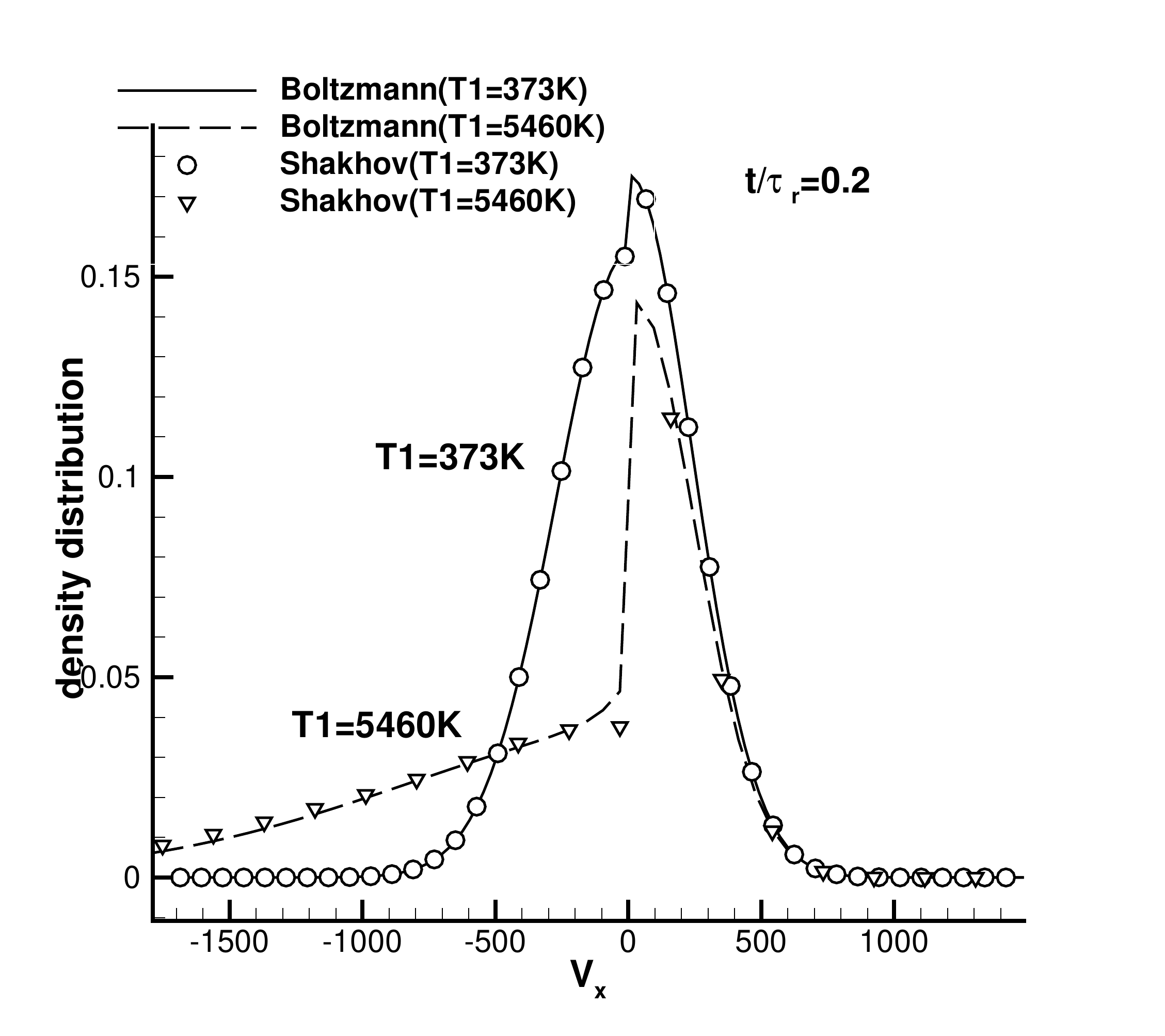}
\includegraphics[width=6cm]{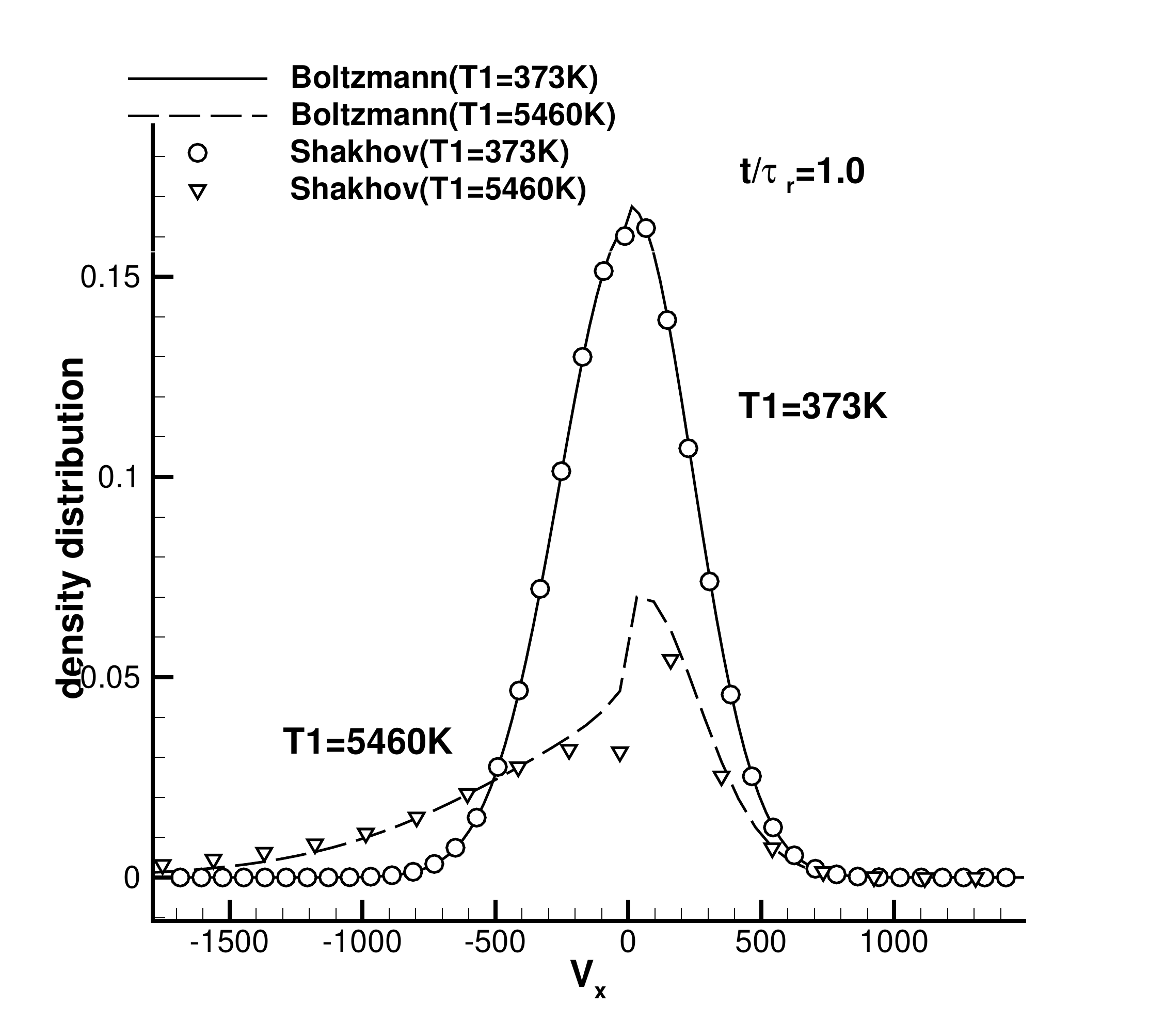}
\includegraphics[width=6cm]{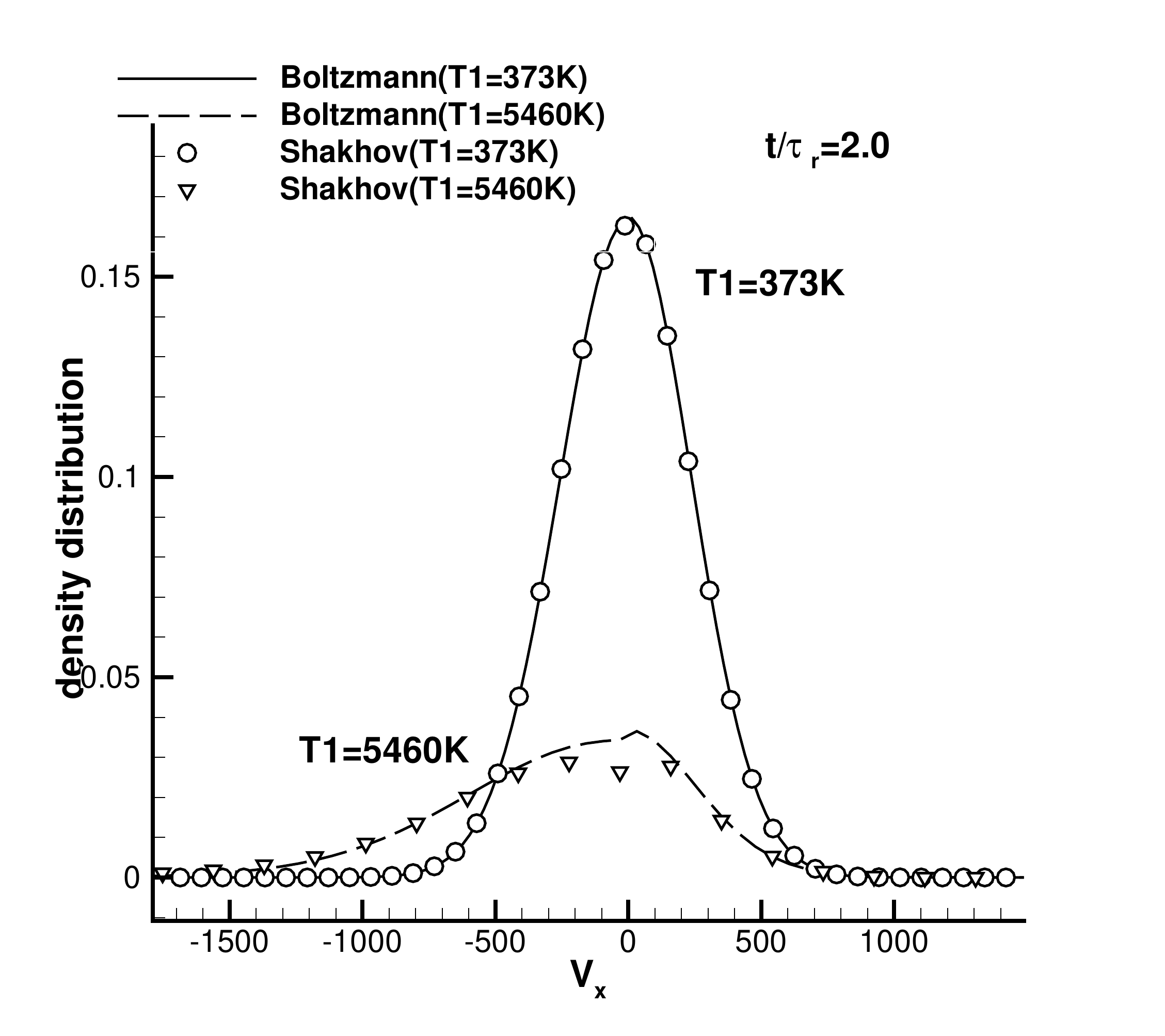}
\includegraphics[width=6cm]{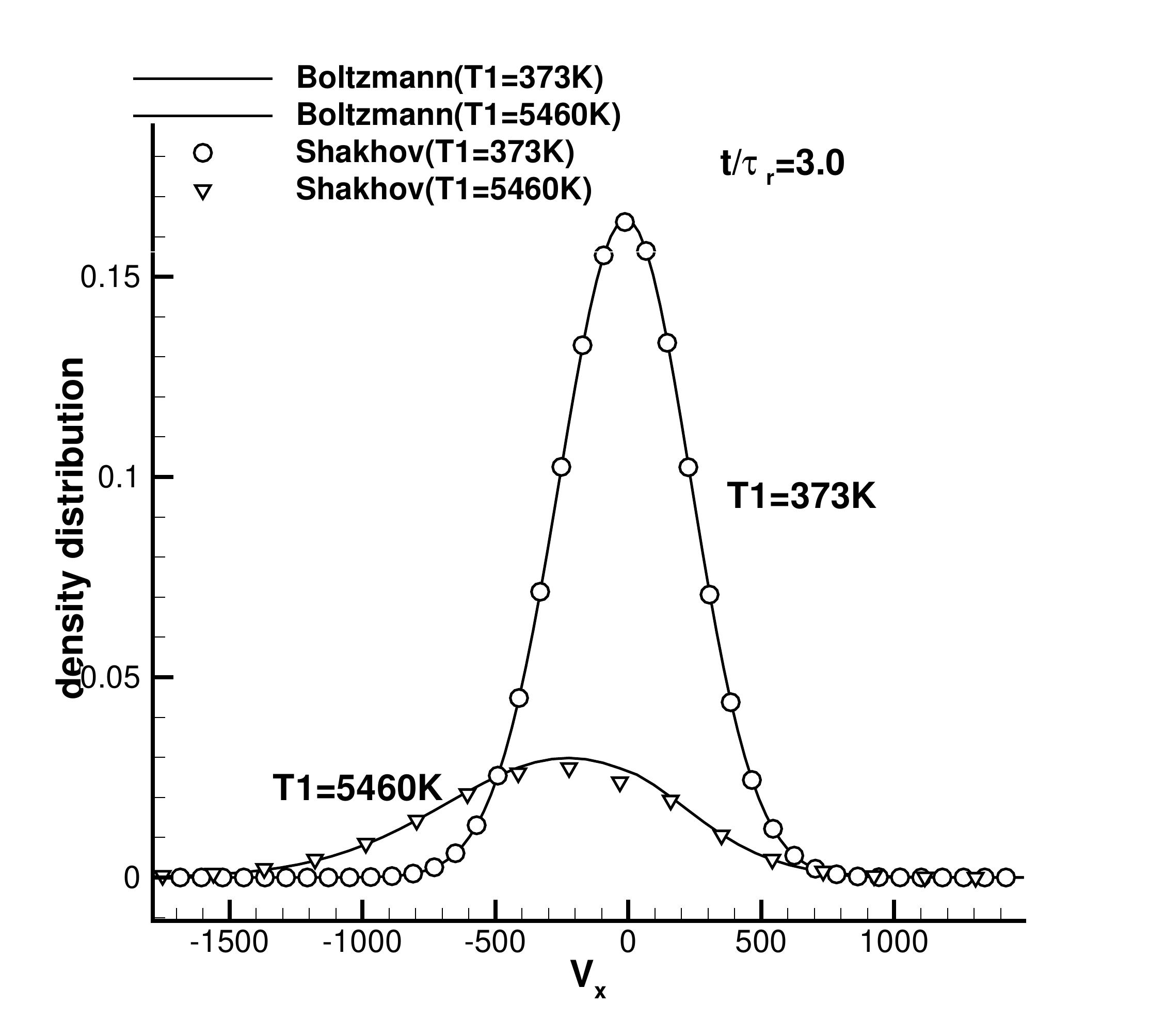}
\includegraphics[width=6cm]{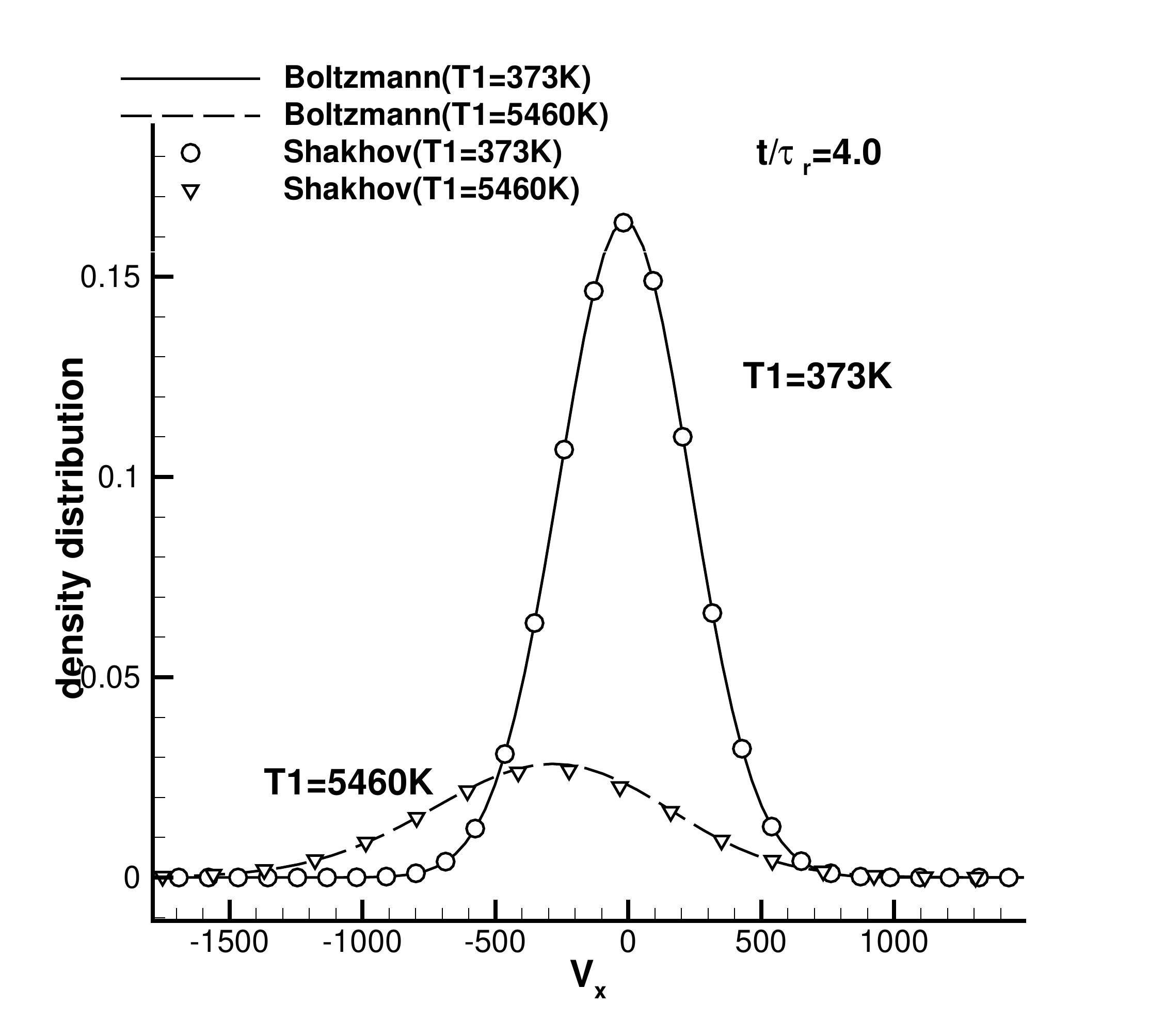}
\caption{Comparison of the x-component distribution functions $f(u,0,0)$ of the relaxation of double half-normal distribution.}
\label{fig:rel2}
\end{figure}

Fig.(\ref{fig:rel2}) shows the x-component distribution functions $f(u,0,0)$
at several output times.
The results show differences near the discontinuity at early times.
However, the deviation in the solutions from Shakhov and Boltzmann decreases with time,
 and becomes negligible when
 {$$t_1>\frac{2}{\rho}\int|f-M| dv\tau_r\approx0.2\tau_{r}$$}
 for the first case, and
 {$$t_2>\frac{2}{\rho}\int|f-M| dv\tau_r\approx2\tau_{r}$$} for the second case.

\subsection{Relaxation of tailored half-Maxwellian distribution}

The tailored half-Maxwellian distribution is designed as follows
$$
f_0(u,v,w)=\frac{2}{\sqrt{\pi}}\frac{\beta_1\beta_2}{\beta_1+\beta_2}\left(\frac{\beta_1}{\sqrt{\pi}}
  \mathrm{e}^{-\beta_1^2u^2}|_{u<0}+\frac{\beta_2}{\sqrt{\pi}}\mathrm{e}^{-\beta_2^2u^2}|_{u\geq0}\right)\frac{\beta_2}{\sqrt{\pi}}\mathrm{e}^{-\beta_2^2v^2}\frac{\beta_3}{\sqrt{\pi}}\mathrm{e}^{-\beta_3^2w^2}, $$
where $\beta_i =\sqrt{m/2kT_i}$, with the initial condition $T_1=273K$, $T_2=373K$, $T_3 =273K$ for the first case,
and $T_1=273K$, $T_2=5460K$, $T_3 =273K$ for the second case.
Here the distribution function is continuous, but with different temperature for the half Maxwellians in the $x$-direction.

\begin{figure}
\center
\includegraphics[width=6cm]{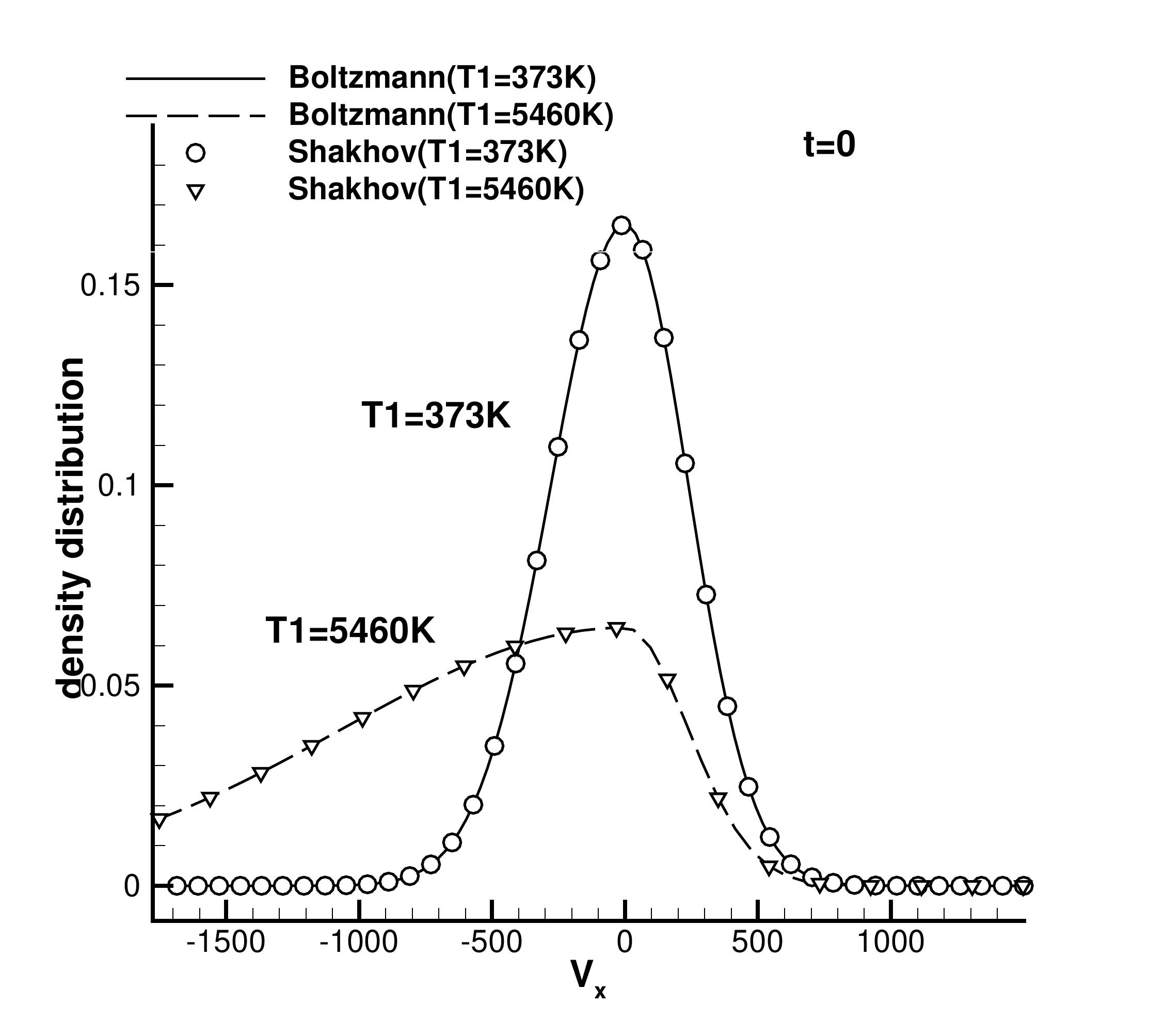}
\includegraphics[width=6cm]{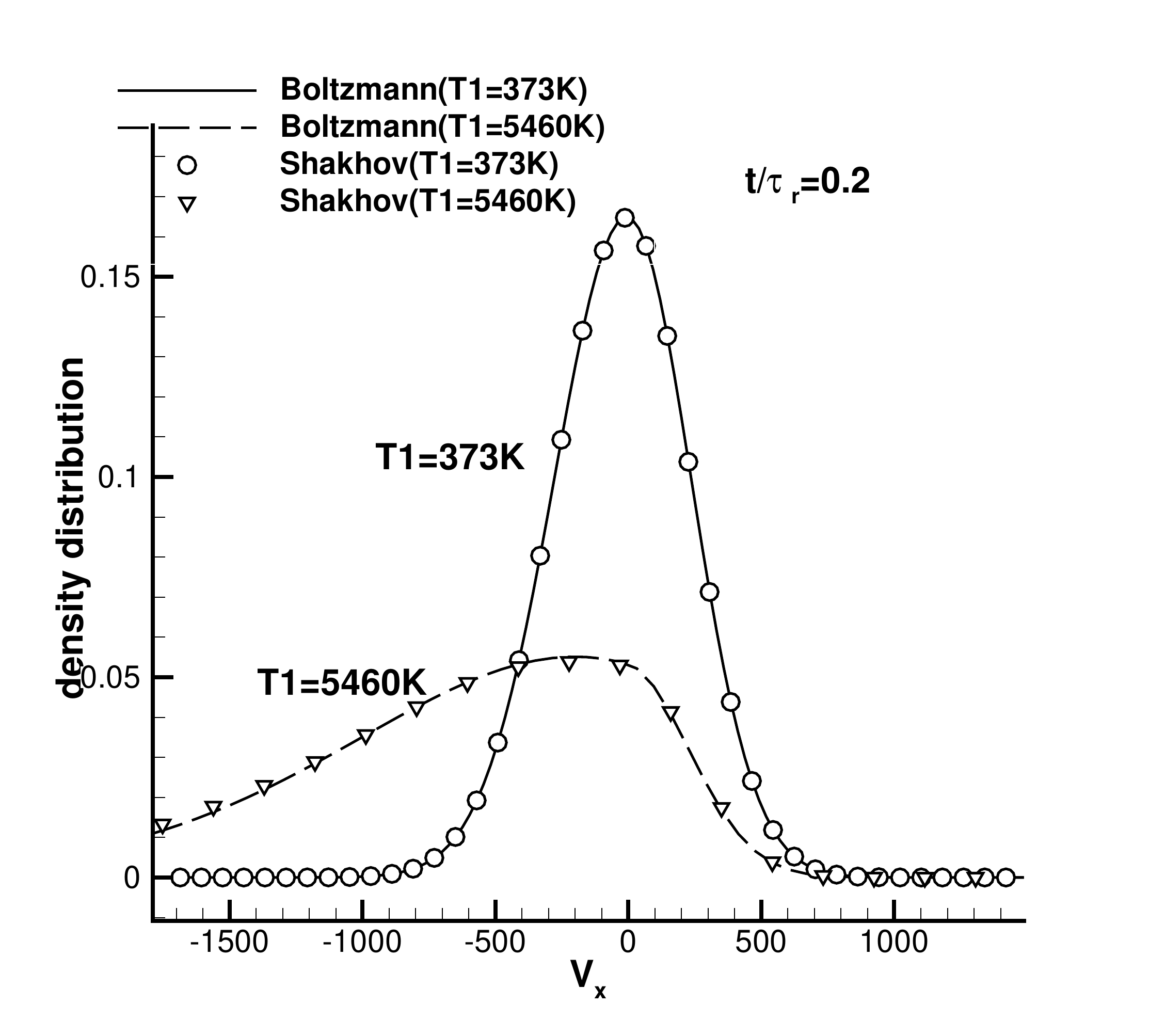}
\includegraphics[width=6cm]{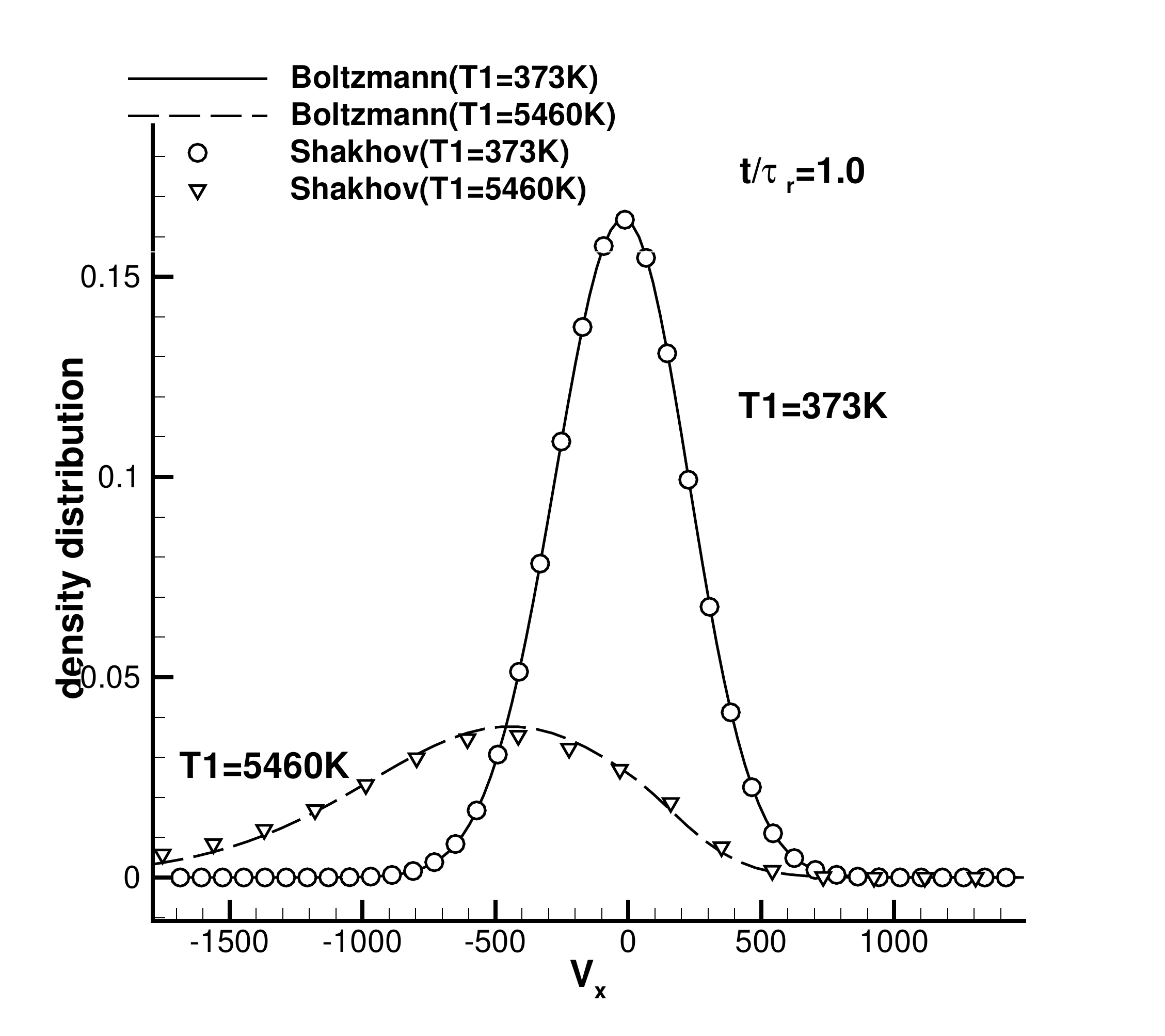}
\includegraphics[width=6cm]{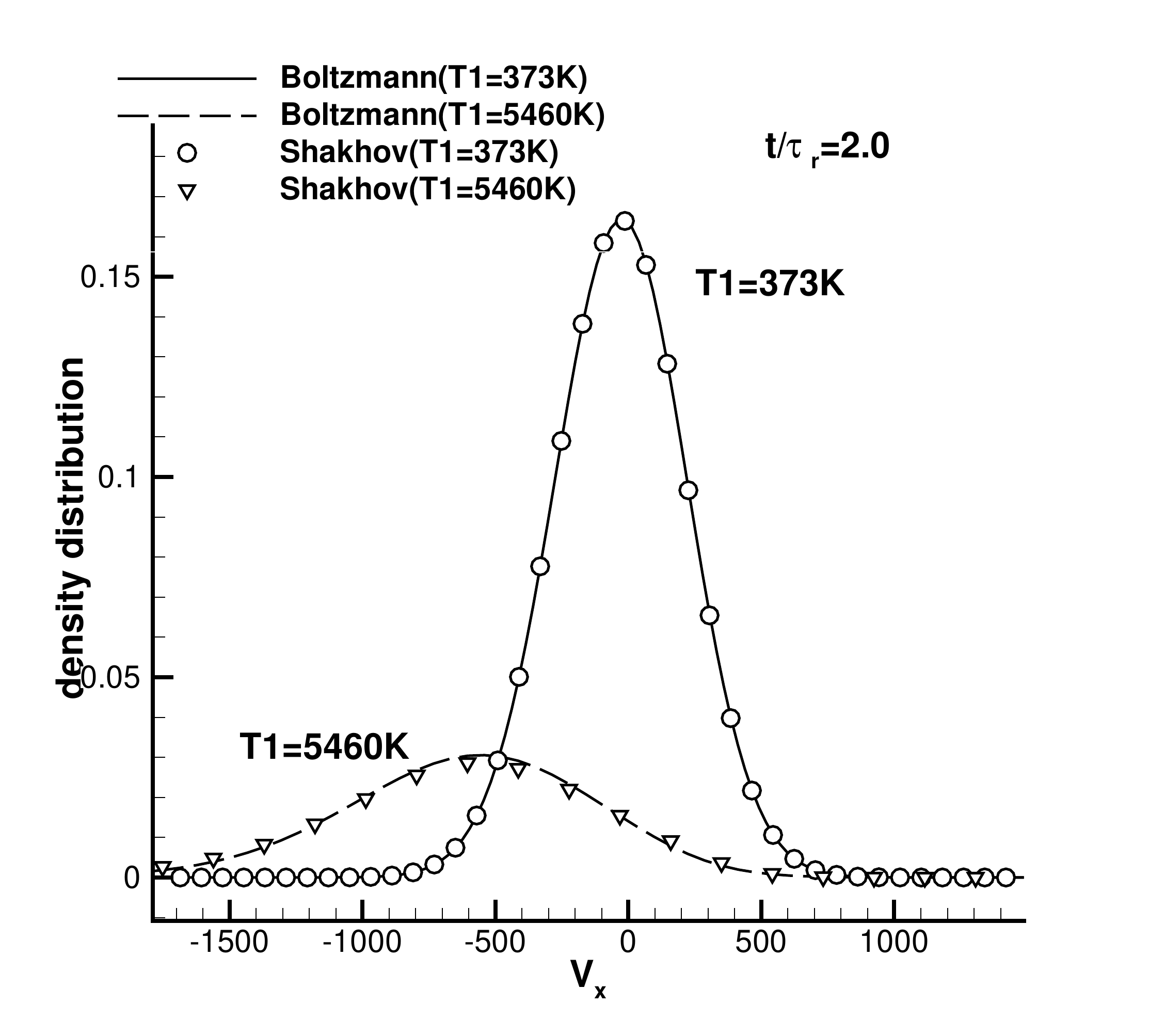}
\includegraphics[width=6cm]{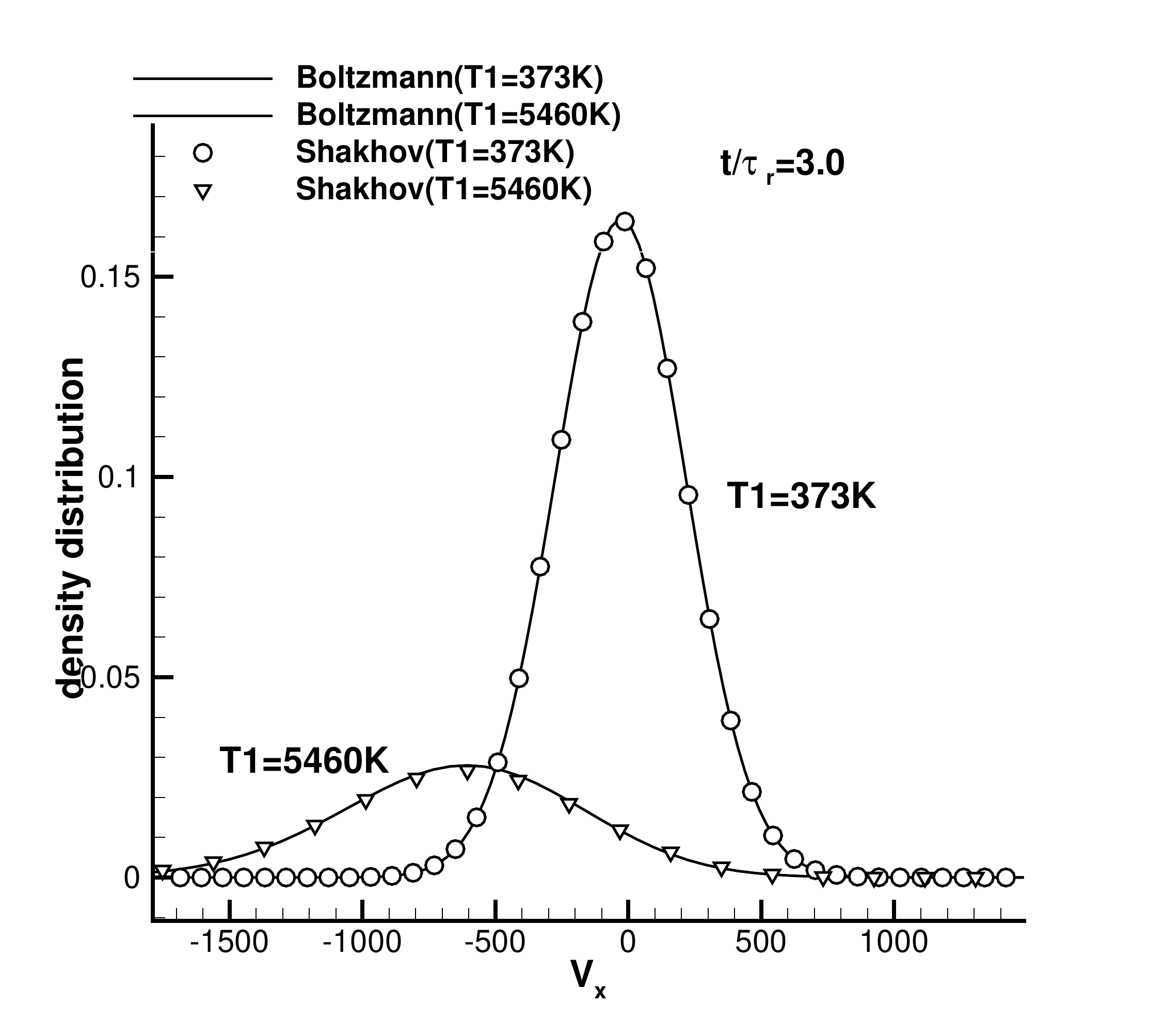}
\includegraphics[width=6cm]{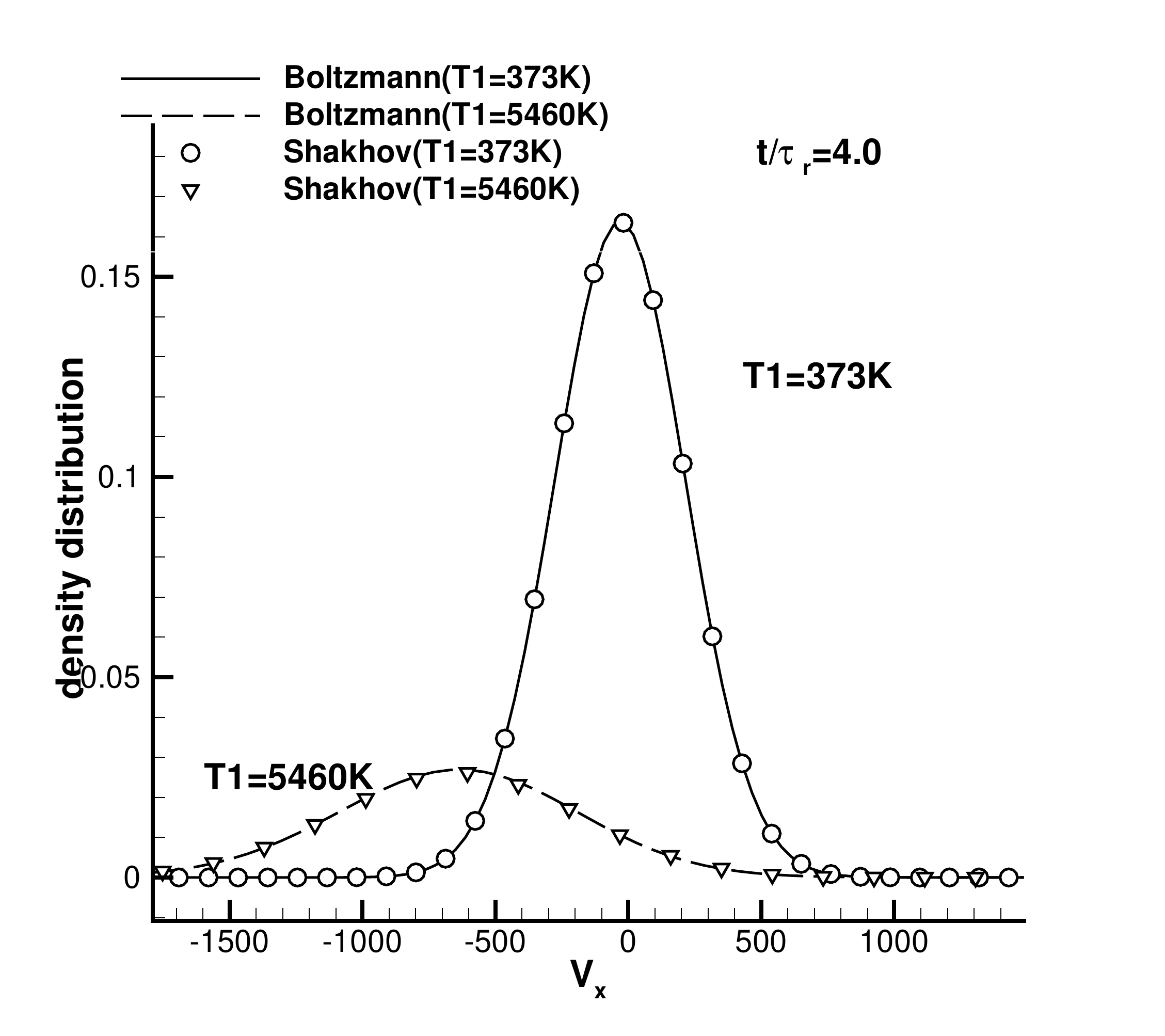}
\caption{Comparison of the x-component distribution functions $f(u,0,0)$ of the relaxation of tailored half-Maxwellian distribution.}
\label{fig:rel3}
\end{figure}

Fig.(\ref{fig:rel3}) shows the time evolution of the x-component distribution functions $f(u,0,0)$ .
The Shakhov (symbols) and Boltzmann (lines) solutions get close after
{$$t_1>\frac{2}{\rho}\int|f-M| dv\tau_r\approx0.2\tau_{r}$$}
for the first case,
and
{$$t_2>\frac{2}{\rho}\int|f-M| dv\tau_r\approx2\tau_{r}$$}
for the second case.

Based on the above observations of all  cases, even for the highly non-equilibrium ones, the Shakhov  and Boltzmann solutions become the same
after $4\tau_{r}$.
Since UGKS is a multiscale method, where the local time step can be varied significantly in terms of the local particle mean free time,
it becomes legitimate to use the kinetic collision model equation if the value of the local time step becomes much larger than
the particle mean free time.
The full Boltzmann collision term is only needed in the highly non-equilibrium region with the time step being less than the mean free time.
The criterion to determine the required region for using the full Boltzmann collision term can be based on the comparison between the local time step
 with the time criterion
 {
 \begin{equation} \label{time-criterion}
    t_{c}=\min(4,D(f,M))\tau_r,
\end{equation}}
 where $D(f,M)=2\int|f-M| dv/\rho$ measures the distance between distribution $f$ and its corresponding equilibrium state $M$.
 The further $f$ deviates from local Maxwellian, the longer it takes to reach the same solution.
 This time criterion will be used in the construction of UGKS with the choices of the full Boltzmann and kinetic model equation.
 The scheme will not be sensitive to such a time criterion
 because even in the cases with $t \leq D(f,M)\tau_r$ the differences in the solutions
 from the full Boltzmann and Shakhov model are not significant.

Past progress on developing asymptotic preserving (AP) schemes \cite{dimarco2011exponential,filbet2010class,jin1999efficient,li2012exponential}
mainly focus on two limiting regimes: the Euler limit and free transport limit.
In the NS regime, although most AP schemes preserve the discrete analogy of the Chapman-Enskog expansion,
viscous effect may not be well resolved due to the large numerical dissipation from the free transport mechanism in the flux evaluation at the cell interface,
or the so-called upwind approach for the transport term across a cell interface \cite{chen-xu}.
In the following, we propose an effective unified scheme, which not only preserves the discrete Chapman-Enskog expansion,
but also leads to the accurate NS solutions.
Basically, in UGKS there is no restriction on the time step in terms of local particle mean free time in the continuum flow regime when a shock capturing approach, without fully resolving the non-equilibrium
shock structure in such a regime, is used and the time step is solely determined by the CFL condition.
Moreover, a local time can be used for the steady state calculation.

\section{Unified gas kinetic scheme with both  full Boltzmann collision term and kinetic model equation}

In this section, we will present the unified gas kinetic scheme (UGKS) in one-dimensional physical space with the inclusion of full Boltzmann collision term.
For two and three-dimensional cases, directional splitting or multidimensional schemes can be constructed accordingly \cite{xu-book}.

\subsection{Unified framework}

The  unified scheme is a direct modeling in the discretized space. It is not targeting to solve any particular partial differential equation, but models
and simulates the flow evolution in the mesh size and time step scales \cite{xu-book,xu2010,Huangxuyu2}.
The physical space is divided into numerical cells with cell size $\Delta x$, and the $j$th-cell is given by $x\in[x_{j-1/2},x_{j+1/2}]$ with cell size $\Delta x=x_{j+1/2}-x_{j-1/2}$.
The temporal discretization is denoted by $t^n$ for the $n$th-time step.
The particle velocity space in $x$-direction is discretized by $2N + 1$ subcells with cell size $\Delta u$,
and the center of $k$th-velocity interval is $u_k = k\Delta u$, and it represents the average velocity $u$ in that interval.
Then, the averaged gas distribution function in cell $j$, at time step $t^n$, and around particle velocity $u_k$, is given by
$$  f_{j,k}^{n}=\frac{1}{\Delta x\Delta u}
  \int_{x_k-1/2}^{x_j+1/2}\int_{u_k-\frac{1}{2}\Delta u}^{u_k+\frac{1}{2}\Delta u} \int
  f(x,t^n,u,\xi)dxdu d\xi, $$
where $\xi$ denotes the freedom in $y,z$ directions with $\xi^2 = w^2 + v^2$ and $d\xi =dvdw$.
The short hand notation $f_{j,k}^n =f(x_j,t^n,u_k)$ will be used for convenience.
The evolution equation for the averaged gas distribution function $f_{j,k}^{n}$ is
\begin{equation}\label{fv}
\begin{aligned}
  f_{j,k}^{n+1}=f_{j,k}^n+
  &\frac{1}{\Delta x}\int_{t^n}^{t^{n+1}}
    (u_{k}f_{j-1/2,k}-u_{k}f_{j+1/2,k})dt\\
  +&\frac{1}{\Delta x}
       \int_{t^n}^{t^{n+1}}\int_{x_{k-1/2}}^{x_{j+1/2}} Q(f,f)_k dxdt,
\end{aligned}
\end{equation}
where $f_{j-1/2,k}=f(x_{j-1/2},t,u_k)$ denotes the time-dependent solution at cell interface $x_{j-1/2}$, the flux transport across the cell interface and the collision term inside each cell need to be modeled.
The above discrete governing equation is actually a physical conservation law, which is valid in all scales.
Here the continuity of the function $f$ is not assumed.
The modeling scale $\Delta x$ and $\Delta t$ in the above equation can be different from the kinetic mean free path and particle mean free time.
So, it is not fully appropriate to state that the above numerical evolution equation is derived from the Boltzmann equation. It is a direct modeling of the physical law.
Instead, the Boltzmann equation can be derived
from the above equation under the constraints on $\Delta x$ and $\Delta t$ to the kinetic scales and with separate consideration of particle transport and collision.
In the above modeling, the dynamics in the flux at cell interface must depend on the ratio of $\Delta t /\tau$, and
a non-splitting treatment of particle transport and collision is needed in the evaluation of flow dynamics around a cell interface.
The non-splitting treatment is critical for the capturing of the NS solutions in the continuum regime
with a relatively coarse mesh with respect to the particle mean free path.
However, many other Boltzmann solvers use the particle free transport in the construction of the cell interface flux, the same as the free transport step in DSMC.
This approach is inconsistent with the physical reality in the continuum regime with $\Delta x \gg \ell$, where the particle will not move freely (straight line) to pass through the cell interface.
This is the main reason for the failure of these flux vector splitting (free transport) modeling schemes in capturing the laminar viscous boundary layer at high Reynolds number.

In  UGKS, besides the evolution equation for $f$ in Eq.(\ref{fv}),
similar to the gas-kinetic scheme (GKS) the update of the conservative variables will be used as well \cite{xu2001},
\begin{equation}\label{fv-w}
  W_j^{n+1}=W^{n}_j+\frac{1}{\Delta x}\int\int_{t^n}^{t^{n+1}}u(f_{j-1/2,k}-f_{j+1/2,k})\psi dt du d\xi,
\end{equation}
for the conservative moments $\psi = (1, u,  \frac{1}{2} (u^2 + \xi^2 ))^T $.
The UGKS is based on the time evolution of two fundamental numerical governing equations (\ref{fv}) and (\ref{fv-w}).
In order to update the gas distribution function and conservative variables, the time-dependent gas distribution function at cell interface and inner cell collision have to be
properly modeled.

\subsection{Gas evolution modeling at a cell interface}

In UGKS, the cell interface flux plays a dominant role to capture the flow dynamics in different scales from kinetic up to the NS ones.
Depending on the scales of $\Delta x$ and $\Delta t$,
the solution at cell interface $f_{j+1/2,k}$ is modeled from an evolution solution of the kinetic model Eq.(\ref{mdl}).
Assume the cell interface is located at $x_{j+1/2}=0$, and the beginning of each time step is $t^n=0$.
The evolution solution is modeled as
\begin{equation}\label{slt}
f(0,t,u_k,\mathbf{\xi})=\frac{1}{\tau_{r}}
\int_{0}^t \tilde{M}(x^\prime,t^\prime,u_k,\mathbf{\xi})\mathrm{e}^{-(t-t^\prime)/\tau_{r}}dt^\prime
+\mathrm{e}^{-t/\tau_{r}}f_0(-u_kt,u_k,\mathbf{\xi}),
\end{equation}
where $x'=-u_k(t-t')$ is the particle trajectory and $f_0(-u_kt,u_k,\xi)$ is the gas distribution function at time $t=0$.
The above solution provides a multiscale modeling from kinetic free transport $f_0$ to the equilibrium realization $\tilde{M}$.
In order to fully determine the evolution solution, the initial condition and the equilibrium states around the cell interface have to be constructed.
Here the conventional reconstruction scheme with nonlinear limiter is used for the initial data reconstruction.
The reconstructed initial condition at time step $t^n$ around the cell interface is
\begin{equation} \label{initial}
f_0(x,u_k,\mathbf{\xi})=\left\{\begin{aligned}
         &f_{j+1/2,k}^L+\sigma_{j,k}x,\quad x\leq 0 ,\\
         &f_{j+1/2,k}^R+\sigma_{j+1,k}x,\quad x>0 .
                          \end{aligned} \right.
                          \end{equation}
In this paper,  the van Leer limiter is used in the reconstruction, where
$$  \sigma_{j,k} = (\mbox{sign}(s_1)+\mbox{sign} (s_2))\frac{|s_1||s_2|}{|s_1|+|s_2|}, $$
with $s_1=(f_{j,k}-f_{j-1,k})/(x_j-x_{j-1})$ and $s_2=(f_{j+1,k}-f_{j,k})/(x_{j+1}-x_j)$.
Certainly, higher-order reconstruction can be used here as well \cite{venugopal}.
The local Maxwellian distribution function $\tilde{M}(f)$ around $(x_{j+1/2},t^n)=(0,0)$ is constructed as,
\begin{equation}\label{eqlm}
\begin{aligned}
\tilde{M}(f)(x,t,u_k,\mathbf{\xi})=&\tilde{M}^{n}_{j+1/2,k}+\partial_x M^{n}_{j+1/2,k}x+\partial_t M^{n}_{j+1/2,k}t\\
=&M^{n}_{j+1/2,k}[1+(1-H(x))a^lx+H(x)a^rx+\tilde{A}t]\\
&+g^{1,n}_{j+1/2,k},
\end{aligned}
\end{equation}
where
$\tilde{M}^{n}_{j+1/2,k} =\tilde{M}(f_0(0,u,\xi))$
and $H(x)$ is Heaviside function defined by
\begin{displaymath}
  H(x)=\left\{\begin{aligned}
  &0, \quad x \leq 0,\\
  &1, \quad x>0.
  \end{aligned} \right.
\end{displaymath}
In 1-D case, the parameters $a^l$,$a^r$ and $\tilde{A}$ depend on the particle velocity in the following form,
$$  a^l=a^l_1+a^l_2u+\frac12 a^l_3(u^2+\xi^2),$$
$$  a^r=a^r_1+a^r_2u+\frac12 a^r_3(u^2+\xi^2), $$
and
 $$ \tilde{A}=\tilde{A}_1+\tilde{A}_2u+\frac12 \tilde{A}_3 (u^2+\xi^2).$$

All parameters $a^l$, $a^r$, and $\tilde{A}$ can be determined based on the correspondence between velocity distribution and conservative variables \cite{xu2010}.
Substituting Eq.(\ref{initial}) and (\ref{eqlm}) into Eq(\ref{slt}), the solution at the cell interface can be expressed as
\begin{equation}\label{integral}
\begin{aligned}
f(x_{j+1/2},t,u_k,\mathbf{\xi})=&(1-\mathrm{e}^{-t/\tau_{r}})(M^{n}_{j+1/2,k}+g^{1,n}_{j+1/2,k})+t\tilde{A}M^{n}_{j+1/2,k}\\
&-(1-\mathrm{e}^{-t/\tau_{r}})\tau_{r}\left((a^r(1-H(u_{k}))+a^lH(u_{k}))u_{k}+\tilde{A}\right)M^{n}_{j+1/2,k}\\
&+\mathrm{e}^{-t/\tau_{r}}[a^r(1-H(u_{k}))+a^lH(u_{k})]u_{k}tM^{n}_{j+1/2,k}\\
&+\mathrm{e}^{-t/\tau_{r}}((f_{i+1/2,k}^L-u_{k}t\sigma_{i,k})H(u_{k})\\
&+(f_{i+1/2,k}^R-u_{k}t\sigma_{i+1,k})(1-H(u_{k}))),
\end{aligned}
\end{equation}
for $t\in[t^n,t^{n+1}]$.
We use $\mathcal{M}_{j+1/2,k}$ and $\mathcal{F}_{j+1/2,k}$ to denote the terms related to the equilibrium distribution and the initial distribution function,
\begin{equation}
  \begin{aligned}
\mathcal{M}_{j+1/2,k}=&(1-\mathrm{e}^{-t/\tau_{r}})(M^{n}_{j+1/2,k}+g^{1,n}_{j+1/2,k})+t\tilde{A}M^{n}_{j+1/2,k}\\
&-(1-\mathrm{e}^{-t/\tau_{r}})\tau_{r}\left((a^r(1-H(u_{k}))+a^lH(u_{k}))u_{k}+\tilde{A}\right)M^{n}_{j+1/2,k}\\
&+\mathrm{e}^{-t/\tau_{r}}[a^r(1-H(u_{k}))+a^lH(u_{k})]u_{k}tM^{n}_{j+1/2,k}\\
    \mathcal{F}_{j+1/2,k}=&\mathrm{e}^{-t/\tau_{r}}((f_{i+1/2,k}^L-u_{k}t\sigma_{i,k})H(u_{k})\\
                         &+(f_{i+1/2,k}^R-u_{k}t\sigma_{i+1,k})(1-H(u_{k}))).
  \end{aligned}
\end{equation}
Based on the distribution function at cell interface, the conservative variables can be updated first  by
\begin{equation}\label{update-macro}
\begin{aligned}
  W_j^{n+1}=W^{n}_j+&\frac{1}{\Delta x}\int_{t^n}^{t^{n+1}}\int u(\mathcal{M}_{j-1/2}-\mathcal{M}_{j+1/2})\psi du d\xi dt\\
                                      +&\frac{1}{\Delta x}\int_{t^n}^{t^{n+1}}\int\sum_k u_{k}(\mathcal{F}_{j-1/2,k}-\mathcal{F}_{j+1/2,k})\psi d\xi dt.
\end{aligned}
\end{equation}

\subsection{Collision term modeling inside each control volume}
Now we have two choices for the collision term modeling inside each control volume, which are
the full Boltzmann collision term $Q(f,f)$ and the Shakhov model $(\tilde{M}(f)-f)/\tau_{r}$.
Depending on the flow regime, the UGKS uses a time step $\Delta t$ which varies significantly relative to the local particle mean free time.
As analyzed in Section 3, starting from a general initial distribution function, the solutions from the full Boltzmann collision term and the kinetic model equation will become the
same after a few mean free times. Therefore, the real place where the full Boltzmann collision term is useful is the region of highly non-equilibrium and with the time step being similar or less than the local mean free time.
As a result, we model the collision term in Eq.(\ref{fv}) as,
\begin{equation}\label{IFSG}
\begin{aligned}
  f_{j,k}^{n+1}=f_{j,k}^n&+\frac{1}{\Delta x}\int_{t^n}^{t^{n+1}}(u_{k} f_{j-1/2,k}-u_{k}f_{j+1/2,k})dt\\
  &+AQ(f_{j}^n,f_{j}^n)_k+B\frac{\tilde{M}(f_j^{n+1})_k-f^{n+1}_{j,k}}{\tau_{r}^{n+1}},
\end{aligned}
\end{equation}
where $\tau_{r}^{n+1}$ denotes the relaxation time at $t^{n+1}$ and the coefficients $A$ and $B$ in the above equation need to satisfy the following constraints:
\begin{description}
  \item[1.] $A+B = \Delta t$ in order to have a consistent collision term treatment.
  \item[2.] The scheme is stable in the whole flow regime.
  \item[3.] In the rarefied flow regime, the scheme gives the Boltzmann solution.
  \item[4.] In continuum regime, the scheme can efficiently recover the Navier-Stokes solutions.
\end{description}
{Based on these constraints, we propose the following choice

\begin{equation}\label{para}
 (A,B)=
 \begin{cases}
 (\beta^n \Delta t,(1-\beta^n) \Delta t), &\Delta t<t_c^n,\\
 (0,\Delta t),  &\Delta t\ge t_c^n,
 \end{cases}
\end{equation}
 with
\begin{equation}\label{beta}
 \beta^n=
 \begin{cases}
 1,&\Delta t<1/\sup_\Upsilon \nu^n,\\
 \exp(1-\Delta t\sup_\Upsilon \nu^n),
 & \Delta t\ge1/\sup_\Upsilon \nu^n.
 \end{cases}
\end{equation}

Here $t_c^n$ is defined by Eq.\eqref{time-criterion} and $\Upsilon$ is the computational domain in velocity space and $\nu^n(u)$ is the collision frequency defined by
\begin{equation}\label{nu}
  \nu^n=\int_\mathrm{R^3}\int_\mathrm{S^2} B(cos\theta, |v-v_*|)f^n(v_*) d\Omega dv_*.
\end{equation}
The collision frequency $\nu^n$ can be calculated using a spectral method \cite{wu} with computational cost $O(N)$, where $N$ denotes the total number of velocity points.
The above choice of parameters presents a continuous transition from the Boltzmann collision term to the kinetic model equation.
The transition parameter $\beta^n$ is proposed based on the following two reasons:
\begin{enumerate}
  \item The Boltzmann collision term is a stiff operator.
        The use of implicit Shakhov model stabilizes the scheme.
        Term $\sup_\Upsilon \nu^n$ in $\beta$ can be viewed as a stiffness indicator.
        When $\sup_\Upsilon \nu^n$ is large, more weight goes to the implicit part.
        This is consistent with the observation in Section 3, where the solutions of Shakhov model and Boltzmann equation get indistinguishable as $\Delta t/\tau_{r}$ increases.
  \item The Boltzmann collision term is physically more accurate than Shakhov model in describing  non-equilibrium flow physics in the kinetic scale,
       and the explicit treatment of the Boltzmann collision term is stable when $\Delta t<1/\sup_\Upsilon \nu^n$.
        In such a case, we can set $\beta=1$ to use the Boltzmann collision term explicitly.
\end{enumerate}

The time criterion $t_c^n$ is proposed to improve the efficiency of the scheme.
When $\Delta t>t_{c}^n$ the solution differences between Boltzmann collision term and Shakhov model are negligible, and a fully implicit Shakhov model can be used
to reduce the computational cost.
The time criterion $t_{c}^n$ is not unique. In many test cases, we can simply assign a fixed value, rather than calculate it in each time step.}

As a result, in the hydrodynamic flow regime the kinetic model will be fully used, and in the kinetic regime
the full Boltzmann equation will be adopted.
Therefore, both the full Boltzmann solution in rarefied regime and the NS solution in continuum regime can be properly obtained.
In the switching regions, the full Boltzmann and kinetic model equations basically present the same solution.
All numerical examples in section 6 show a smooth transition across all regimes.
In the low transition and near continuum flow regime, with the adaptation of large mesh size relative to the local particle mean free path, the kinetic model will be used in most of the domain,
especially when the local time step associated with local mesh size is used for steady state solutions.
In the continuum regime, when the physical shock structure is not necessarily resolved by the numerical cell size,
a shock capturing scheme will be emerged automatically from the the above UGKS.

{In summary, starting from time step $t^n$, the UGKS updates flow variables with four steps:
\begin{description}
  \item[Step 1] Reconstruct macroscopic variables and velocity distribution by Eq.\eqref{initial}; calculate numerical flux based on the time dependent solution Eq.\eqref{integral}; and update the macroscopic variables by Eq.\eqref{update-macro}.
  \item[Step 2] Calculate $t_c^n$ by Eq.\eqref{time-criterion}. If $\Delta t>t_c^n$, go to Step 3; if $\Delta t\le t_c^n$, go to Step 4.
  \item[Step 3] Calculate $\tilde{M}^{n+1}$ from the updated conservative variables, and update velocity distribution function by Eq.\eqref{IFSG} with $(A,B)=(0,\Delta t)$.
  \item[Step 4] Calculate $\tilde{M}^{n+1}$, $Q(f^n,f^n)$, $\nu^n$ and $\beta^n$ in Eq.\eqref{beta} and Eq.\eqref{nu}. Update velocity distribution function by Eq.\eqref{IFSG} with $(A,B)=(\beta^n \Delta t,(1-\beta^n) \Delta t)$.
\end{description}}

\section{Numerical analysis of UGKS}

In this section, we discuss the properties of UGKS.
The stability is discussed for a homogeneous case,
and the AP property is discussed in the Euler, NS, and free molecular regimes.

\subsection{Stability and convergence in a homogeneous case}
{ We show the stability of scheme in a spacial homogeneous case with the full Boltzmann and BGK collision model in both cases of $\Delta t \lessgtr t_c$.

\noindent {\bf Case (1):} $\Delta t < t_c$

Assume $||f_0||_{L^1}=\rho$,
and $\Delta t<t_c$.
By splitting $Q(f,f)$ into gain term $Q_+(f,f)$ and loss term $\nu f$,
we rewrite the UGKS as
\begin{equation}\label{stable1}
\begin{aligned}
f^{n}=&\frac{1-\Delta t \nu^{n-1}\beta^{n-1}}
   {1+\Delta t(1-\beta^{n-1})/\tau_r}f^{n-1}+
   \frac{\Delta t \beta^{n-1}}
   {1+\Delta t(1-\beta^{n-1})/\tau_r}Q_+^{n-1}\\
   &+\frac{\Delta t(1-\beta^{n-1})/\tau_r}{1+\Delta t(1-\beta^{n-1})/\tau_r}M .
\end{aligned}
\end{equation}
Based on the definition of $\beta$ Eq.\eqref{beta},
we have
$$0<\beta^{n-1}<1, \quad 0<\Delta t \nu^{n-1}\beta^{n-1}<1,$$
which show that $f^n$ is a convex combination of $f^{n-1}$, $Q_+^{n-1}$ and $M$.
Hence the numerical solution $f$ keeps positive.
In addition, it is proved that the numerical Boltzmann collision term calculated by the spectral method preserves the total mass \cite{pareschi2000stability}.
By taking $L^1$ norm to Eq.\eqref{stable1}, we have
\begin{equation}\label{L1}
\begin{aligned}
  ||f^{n}||_{L1}=&\frac{1}{1+\Delta t(1-\beta^{n-1})/\tau_r}||f^{n-1}||_{L1}+\frac{\Delta t(1-\beta^{n-1})/\tau_r}{1+\Delta t(1-\beta^{n-1})/\tau_r}||M||_{L1}\\
  =&\rho,
\end{aligned}
\end{equation}
which shows that the UGKS solution is positive with a fixed $L^1$ norm and the solution is stable in such sense.
In addition, from
$$\lim_{\Delta t/\max(\tau,1/\sup \nu) \rightarrow \infty} \beta=0,$$
we get $$\lim_{\Delta t/\max(\tau,1/\sup \nu) \rightarrow \infty} f=M,$$
which implies the L-stable property of the scheme.

By iteration, we can write down the numerical solution of UGKS as,
\begin{equation}\label{convergence1}
   \begin{aligned}
   f^{n}=&\prod_{i=0}^{n-1}\left(1-\frac{(1-\beta^i)/\tau_r+
   \nu^i\beta^i}{1+\Delta t(1-\beta^i)/\tau_r}\Delta t\right)f^0\\
   ~&\\
   &+\sum_{s=0}^{n-1}\frac{\Delta t \beta^s Q^s_+}{1-\Delta t \nu^s \beta^s}
    \prod_{i=s}^{n-1}\frac{1-\Delta t \nu^i\beta^i}
    {1+\Delta t(1-\beta^i)/\tau_r}\\
   ~&\\
    &+\sum_{s=0}^{n-1}\frac{\Delta t(1-\beta^s)/\tau_rM}{1-\Delta t(1-\beta^s)/\tau_r}
    \prod^{n-1}_{i=s}\frac{1-\Delta t \nu^i \beta^i}{1+\Delta t(1-\beta^i)/\tau_r}.
   \end{aligned}
\end{equation}
Let
$$E_{n,\Delta t}=\prod_{i=0}^{n-1}
\left(1-\frac{(1-\beta^i)/\tau_r+\nu^i\beta^i}{1+\Delta t(1-\beta^i)/\tau_r}\Delta t\right).$$
Based on definition of $\beta$ Eq.\eqref{beta}, we have
\begin{equation}\label{stability-condition}
0\le\frac{(1-\beta^i)/\tau_r+\nu^i\beta^i}{1+\Delta t(1-\beta^i)/\tau_r}\Delta t\le1,
\end{equation}
from which we get
\begin{equation}\label{exp}
  \begin{aligned}
  E_{n,\Delta t}
  =&\exp\left\{-\sum_{i=0}^{n-1}\Delta t\frac{(1-\beta^i)/\tau_r+\nu^i\beta^i}{1+\Delta t(1-\beta^i)/\tau_r}\right.\\
   ~\\
  &\left.+\sum_{i=0}^{n-1}\left(\ln\left(1-\frac{(1-\beta^i)/\tau_r+\nu^i\beta^i}{1+\Delta t(1-\beta^i)/\tau_r}\Delta t\right)\right.\right.\\
   ~\\
  &+\left.\left.\Delta t\frac{(1-\beta^i)/\tau_r+\nu^i\beta^i}{1+\Delta t(1-\beta^i)/\tau_r}\right)\right\}\\
   ~\\
  \rightarrow&\exp\left(-\int_0^t\left(\nu(f)\beta+\frac{1-\beta}{\tau}\right)dt\right).
  \end{aligned}
\end{equation}
Based on the regularity property of the $Q^+$ term given by Lions \cite{lions1994}, when $\Delta t$ goes to zero, $\beta$ goes to one and the solution of the UGKS Eq.\eqref{stable1} converges to
\begin{equation}
\begin{aligned}
  f(t)=&f_0\mathrm{e}^{-\int_0^t\left(\nu\beta+\frac{1-\beta}{\tau}\right)ds}
  +\int_0^t \left(\beta Q_+(f,f)(s)
  +\frac{1-\beta}{\tau_r}M\right)\mathrm{e}^{-\int_s^t\left(\nu\beta
  +\frac{1-\beta}{\tau}\right)d\sigma}ds\\
  =&f_0\mathrm{e}^{-\int_0^t\nu ds}
  +\int_0^t Q_+(f,f)(s)\mathrm{e}^{-\int_s^t\nu
  d\sigma}ds,
\end{aligned}
\end{equation}
which is the exact Boltzmann solution.}

\noindent {\bf Case (2):} $\Delta t > t_c$

If $\Delta t>t_c$, the parameters are set as $A=0$ and $B=\Delta t$, and the scheme can be written as
\begin{displaymath}
  \frac{[f^{n+1}-M]-[f^n-M]}{\Delta t}=\frac{{M}(f^{n+1})-f^{n+1}}{\tau_{r}},
\end{displaymath}
which gives,
\begin{displaymath}
  f^{n+1}-M=r_i [f^n-M],
\end{displaymath}
and
\begin{displaymath}
\begin{aligned}
  |r_i|=&\left|\frac{1}{1+\frac{\Delta t}{\tau_{r}}}\right|\\
      <&1,
\end{aligned}
\end{displaymath}
Moreover,
\begin{displaymath}
 \lim_{\Delta t/\tau_{r}\rightarrow \infty} |r_i| =0,
\end{displaymath}
which implies the L-stable property of the scheme as well.

From above discussion, we find the scheme is always stable for a spatially homogeneous case. The exact Boltzmann solution can be obtained from UGKS once the full Boltzmann collision term is used.

\subsection{Asymptotic preserving analysis}

{In this section, we discuss the asymptotic preserving (AP) property of the UGKS (\ref{IFSG}).
AP property has been studied in the past few years \cite{dimarco2011exponential,filbet2010class,jin1999efficient,li2012exponential,liu2010analysis}.
However, the flow regimes considered are mostly two limiting cases, namely the rarefied regime and Euler regime.
Here, we propose a stronger AP property which also includes the NS regime. A scheme for kinetic equation is AP if
\begin{itemize}
  \item Holding the mesh size, it degenerates to collisionless Boltzmann equation as Knudsen number goes to infinity and becomes a suitable scheme for the Euler equation as Knudsen number goes to zero.
  \item Implicit part in the collision term can be effectively calculated explicitly.
  \item When Knudsen number is small, it preserves the discrete Chapman-Enskog expansion. When the local time step is much larger than the local particle collision time,
  i.e. $\Delta t \gg \tau_r$, the scheme is at least a second order time accurate scheme for the NS equations.
\end{itemize}}
For the sake of simple notation, the discussion is in one dimensional physical space.

\subsubsection{Collisionless limit}

In free transport regime, or $\tau_{r}\rightarrow\infty$, scheme  Eq.(\ref{IFSG}) becomes,
\begin{displaymath}
\begin{aligned}
  f_{j,k}^{n+1}=f_{j,k}^{n}&+\frac{1}{\Delta x}((\Delta t f_{j+1/2,k}^L-\frac12\Delta t^2u_k\sigma_{j,k})H[u_k]\\
  &+(\Delta t f_{j+1/2,k}^R-\frac12\Delta t^2u_k\sigma_{j+1,k})(1-H[u_k])),
\end{aligned}
\end{displaymath}
which is a second order upwind scheme for collisionless Boltzmann equation.

\subsubsection{Euler limit}

In the following, we are going to analyze the UGKS in the continuum regime with the conditions $\tau_r \rightarrow 0$ and $\Delta t \gg \tau_r$.

We assume that $f$ has a continuous second order derivative with respect to $x$,
 and the initial condition at the beginning of each time step is reconstructed by central difference. Based on the assumption, the solution at cell interface Eq.\eqref{integral} becomes

\begin{equation}\label{integralconti}
 \begin{aligned}
 f(x_{j+1/2},t,u_k,\mathbf{\xi})=&(1-\mathrm{e}^{-t/\tau_{r}})\left(M^{n}_{j+1/2,k}+g^{1,n}_{j+1/2,k}-\tau_{r}(au_{k}M^{n}_{j+1/2,k}+\tilde{A}M^{n}_{j+1/2,k})\right)\\
&+\tilde{A}M^{n}_{j+1/2,k}t+\mathrm{e}^{-t/\tau_{r}}(au_{k}M^{n}_{j+1/2,k}-u_{k}\sigma_k)t\\
&+\mathrm{e}^{-t/\tau_{r}}f_{i+1/2,k}\\
 \end{aligned}
\end{equation}
where $\sigma_k=\frac{f^n_{j+1,k}-f^n_{j,k}}{\Delta x}$.

In the Euler limit when $\tau_{r}\rightarrow0$, the parameters in scheme (\ref{IFSG}) take values $(A,B)=(0,\Delta t)$.
Taking limit of $\tau_{r}^{n+1}\rightarrow0$ in Eq.(\ref{IFSG}), we have
\begin{equation}\label{ap-euler}
\begin{aligned}
   \lim_{\tau_{r}^{n+1}\rightarrow0}f_{j,k}^{n+1}
  =&\lim_{\tau_{r}^{n+1}\rightarrow0}\frac{f_{j,k}^n
    +\frac{1}{\Delta x}\int_{t^n}^{t^{n+1}}(f_{j-1/2,k}-f_{j+1/2,k})dt
    +\frac{\Delta t }{\tau_{r}^{n+1}}g^1(f^{n+1}_j)_{k}}{1+\frac{\Delta t}{\tau_{r}^{n+1}}}\\
   &+\lim_{\tau_{r}^{n+1}\rightarrow0}\frac{\frac{\Delta t }{\tau_{r}^{n+1}}M(f^{n+1}_j)_{k}}
     {1+\frac{\Delta t}{\tau_{r}^{n+1}}}\\
  =&M(f^{n+1}_{j})_{k}.
\end{aligned}
\end{equation}
Note that based on the definition of $g^1$ in Eq.\eqref{shakhov-M}, we have $g^1(f^{n+1}) \sim O(\tau^{n+1})$,
thus, $$\lim_{\tau_{r}^{n+1}\rightarrow0}
 \frac{\frac{\Delta t }{\tau_{r}^{n+1}}
 g^1(f^{n+1}_j)_{k}}{1+\frac{\Delta t}{\tau_r^{n+1}}}=0.$$
Therefore, as $\tau_{r}$ approaching to zero, the numerical solution $f$ converges to $M(f)$.
By linear interpolation, the cell interface distribution function $f_{j+1/2}(t^n,u_k,\mathbf{\xi})$ is
\begin{equation}\label{init1}
  \begin{aligned}
  f_{j+1/2}(t^n,u_k,\mathbf{\xi})
  = &f_{j,k}^n+\frac{f_{j+1,k}^n-f_{j,k}^n}{\Delta x}\frac12
     \Delta x\\
  = &M^{n}_{j+1/2,k}+O(\Delta x ^2),
  \end{aligned}
\end{equation}
where
\begin{align*}
M^{n}_{j,k}=M(f_j(t^n))(u_k), ~~M^{n}_{j+1/2,k}=M(f_{j+1/2}(t^n))(u_k).
\end{align*}
Substituting initial condition Eq.\eqref{init1} into the integral solution Eq.\eqref{integralconti}, we have for $t\in[t^n,t^{n+1}]$,
\begin{equation}\label{1order}
\begin{aligned}
 f(x_{j+1/2},t,u_k,\mathbf{\xi})=&M^n_{j+1/2,k}+t\tilde{A}M^n_{j+1/2,k}\\
 &+t\mathrm{e}^{-t/\tau_{r}}(au_kM^n_{j+1/2,k}-\sigma_ku_k)\\
 &-\tau_r(1-\mathrm{e}^{-t/\tau_{r}})(au_kM^n_{j+1/2,k}+\tilde{A}M^n_{j+1/2,k})\\
 &+(1-\mathrm{e}^{-t/\tau_{r}})g^{1,n}_{j+1/2,k}\\
 =&M(f_{j+1/2,k}(t))+O(\Delta x^2, \Delta t^2),
\end{aligned}
\end{equation}
as $\tau_r\rightarrow 0$.
By taking conservative moments $\psi$ to Eq.(\ref{1order}), one can get the cell interface flux for conservative variables,
\begin{equation}
\begin{aligned}
  F_{w}
  = \left(
  \begin{array}{c}
   \rho U \\
   \rho U^2+P\\
    (\rho E+P)U \\
     \end{array}
      \right)_{j+1/2}+O(\Delta t^2, \Delta x^2).
\end{aligned}
\end{equation}
Substituting the microscopic flux Eq.(\ref{1order}) into Eq.(\ref{fv-w}), taking conservative moments, and keeping $O(1)$ terms, we get the discrete Euler system,
 \begin{equation}
     \left\{
     \begin{aligned}
       &\frac{\rho^{n+1}-\rho^{n}}{\Delta t}+\frac{1}{\Delta t\Delta x}\int_{t^{n}}^{t^{n+1}} \left[(\rho U)_{j+1/2}- (\rho U)_{j-1/2}\right] dt=O(\Delta t^2,\Delta x^2),\\
       &\frac{\rho^{n+1}U^{n+1}-\rho^{n}U^{n}}{\Delta t}+\frac{1}{\Delta t\Delta x}\int_{t^{n}}^{t^{n+1}}
           \left[(\rho U^2+P)_{j+1/2}-(\rho U^2+P)_{j-1/2}\right] dt \\
           & =O(\Delta t^2,\Delta x^2),\\
       &\frac{(\rho E)^{n+1}-(\rho E)^{n}}{\Delta t}+\frac{1}{\Delta t\Delta x}\int_{t^{n}}^{t^{n+1}}
       \left[((\rho E+P)U)_{j+1/2}-((\rho E+P)U)_{j-1/2}\right]dt \\
       & =O(\Delta t^2,\Delta x^2),
     \end{aligned}
     \right.
   \end{equation}
from which we can observe that UGKS is a second order scheme for the Euler equations as Knudsen number goes to zero,
   \begin{equation}
     \left\{
     \begin{aligned}
       &\frac{\partial \rho}{\partial t}+\frac{\partial \rho U}{\partial x}=O(\Delta t^2, \Delta x^2),\\
       &\frac{\partial (\rho U ) }{\partial t}+\frac{\partial (\rho U^2+P)}{\partial x}=O(\Delta t^2, \Delta x^2),\\
       &\frac{\partial (\rho E)}{\partial t}+\frac{\partial  (\rho E+P)U}{\partial x}=O(\Delta t^2, \Delta x^2).
     \end{aligned}
     \right.
   \end{equation}

\subsubsection{Navier-Stokes limit}

Next, we analyze the asymptotic property of the scheme (\ref{fv}) in the Navier-Stokes regime with
small $\tau_r $ and $\Delta t$, under the condition $\Delta t \gg \tau_r$.

The following analysis is given for a well resolved flow region with $\Delta t > t_{c}$ and the initial condition at the beginning of each time step is reconstructed by central difference.
When $\tau_r$ is small, the parameters in the scheme take the
values $(A,B)=(0,\Delta t)$.
The initial condition is assumed to be the form $f_0=M(f_0)+O(\tau_{r})$.

Following the Chapman-Enskog theory, when $\tau_{r}\ll 1$,
the cell averaged solution $f_j(t,u,\mathbf{\xi})$ can be formally written as an asymptotic expansion of small parameter $\tau_{r}$,
\begin{equation}\label{ce}
  f_j=f_j^0+\tau_{r} f_j^1+O(\tau_{r}^2).
\end{equation}
The modified equilibrium distribution function $\tilde{M}(f_j)$ can be expanded as
\begin{displaymath}
  \tilde{M}(f_j)=M(f_j)+g^1(f_j),
\end{displaymath}
where $g^1$ is of order $\tau_r$.
The stress tensor and heat flux can be expanded as
\begin{displaymath}
  \begin{aligned}
    \theta &= TI+\tau_{r} \theta^1+O(\tau_{r}^2),\\
    q&=0+\tau_{r} q^1+O(\tau_{r}^2),
  \end{aligned}
\end{displaymath}
where
\begin{displaymath}
  \begin{aligned}
    \theta^1&=\frac{1}{\rho}\int_{\mathbb{R}^3} (v-U)^2 f_j^1(t,u,\mathbf{\xi})du d\xi,\\
    q^1&=\frac12\int_{\mathbb{R}^3}((u-U)^2+\xi^2)(u-U) f_j^1 (t,u,\mathbf{\xi})du d\xi .
  \end{aligned}
\end{displaymath}

Assuming that scheme (\ref{IFSG}) depends continuously on $t\in[t^n,t^{n+1}]$, by taking time derivative, we have
{
\begin{equation}\label{CTDS}
   \partial_{t} f_{j,k}(t)+u_k\frac{f_{i+1/2,k}(t)-f_{i-1/2,k}(t)}{\Delta x}
   =\frac{\tilde{M}(f_{j}(t))_k-f_{j,k}(t)}{\tau_{r}}+O(t).
\end{equation}}
From Eq.\eqref{ap-euler}, we have $f_j^0(t^n)=M(f_j(t^n))$.
Balancing the $O(1)$ terms in Eq.\eqref{CTDS}, we have
\begin{displaymath}
\begin{aligned}
  f_{j,k}^1(t^n)=&g^1(f_{j}(t^n))_{k}/\tau_r-\partial_{t} M(f_j(t^n))_{k}\\
                            &-u_{k}\frac{M(f_{j+1/2}(t^n))_{k}-M(f_{j-1/2}(t^n))_{k}}{\Delta x}+O(\Delta t)\\
                          =&g^1(f_{j}(t^n))_{k}/\tau_r-\partial_{t} M(f_j(t^n))_{k}-u_{k}\partial_x M(f_j(t^n))_{k}+O(\Delta t,\Delta x^2).
\end{aligned}
\end{displaymath}
Then, based on the cell center values, by interpolation the distribution function at cell interface at $t^n$ up to order $O(\tau_{r})$ is
\begin{equation}\label{init2}
  \begin{aligned}
  &f_{j+1/2}(t^n,u_k,\mathbf{\xi})\\
  =&f_{j,k}^{n}+\frac{f_{j+1,k}^{n}-f_{j,k}^{n}}{\Delta x}\frac12 \Delta x\\
  =&M^{n}_{j,k}+g^{1,n}_{j,k}-\tau_{r}(\partial_tM^{n}_{j,k}+u_k\partial_xM^{n}_{j,k}
  +O(\Delta x^2,\Delta t))\\
  &+\partial_x (M^{n}_{j,k}+g^{1,n}_{j,k}-\tau_{r}(\partial_tM^{n}_{j,k}+u_k\partial_xM^{n}_{j,k}))
  \frac12 \Delta x\\
  &+O(\tau_{r} \Delta t, \tau_{r} \Delta x^2)+O(\Delta x^3)\\
  =&M^{n}_{j+1/2,k}+g^{1,n}_{j+1/2,k}-\tau_{r}(\partial_tM^{n}_{j+1/2,k}+u_k\partial_xM^{n}_{j+1/2,k})\\
   &+O(\tau_{r}\Delta t, \Delta x^2),
  \end{aligned}
\end{equation}
where
\begin{displaymath}
\begin{aligned}
&g^{1,n}_{j,k}=g^1(f_j(t^n))(u_k), g^{1,n}_{j+1/2,k}=g^1(f_{j+1/2}(t^n))(u_k),\\
&M^{n}_{j,k}=M(f_j(t^n))(u_k), M^{n}_{j+1/2,k}=M(f_{j+1/2}(t^n))(u_k).\\
\end{aligned}
\end{displaymath}
Substituting initial condition (\ref{init2}) into the integral solution (\ref{integral}),  we have for $t\in[t^n,t^{n+1}]$,
  \begin{equation}\label{tauorder}
    \begin{aligned}
    &f_{j+1/2}(t,u_k,\mathbf{\xi})\\
    =&M^{n}_{j+1/2,k}+g^{1,n}_{j+1/2,k}-\tau_{r}(\tilde{A}M^{n}_{j+1/2,k}+u_kaM^{n}_{j+1/2,k})\\
    &+t\tilde{A}M^{n}_{j+1/2,k}\\
    &-\tau_{r} t\mathrm{e}^{-t/\tau_{r}}(u_k\partial_x(\tilde{A}M^{n}_{j+1/2,k}+u_kaM^{n}_{j+1/2,k})+O(\tau_{r}\Delta t, \tau_{r}\Delta x^2))\\
    &+O(\Delta x^2)\\
    =&M_{j+1/2,k}(t)+g^1_{j+1/2,k}(t)-\tau_{r}(\partial_t M_{j+1/2,k}(t)+u_k\partial_x M_{j+1/2,k}(t))\\
    &+O(\tau_{r} \Delta t,\Delta t^2,\Delta x^2),
    \end{aligned}
   \end{equation}
where
\begin{displaymath}
\begin{aligned}
&M_{j+1/2,k}(t)=M(f_{j+1/2}(t))(u_k), ~~g^{1}_{j+1/2,k}(t)=g^1(f_{j+1/2}(t))(u_k).\\
\end{aligned}
\end{displaymath}
The non-equilibrium part in Eq.(\ref{tauorder}) is
   \begin{displaymath}
  \begin{aligned}
    f_{j+1/2}^1(t,u_k,\mathbf{\xi})=&g^1_{j+1/2}(t,u_k)-\tau_{r}(\partial_t M_{j+1/2}(t,u_k)+u_k\partial_x M_{j+1/2}(t,u_k))\\
                                       =&M_{j+1/2}(t,u_k)\frac{4(1-\mbox{Pr})\lambda^2}{5\rho}c'_{k}\cdot q(2\lambda c'^2_{k}-5)\\
                                          & -M_{j+1/2}(t,u_k)\left(\tau_{r}\left(\lambda c'^2_k-\frac52\right)c'_k\frac{\partial}
                                            {\partial x}\ln T+\frac43\tau_{r}\lambda c'^2_{k}\frac{\partial U}{\partial x}\right),\\
  \end{aligned}
   \end{displaymath}
   where $c'_k$ denotes the $k$-th peculiar velocity.

By taking conservative moments $\psi$ to Eq.(\ref{tauorder}), one can get the cell interface flux for conservative variables up to $O(\tau_{r})$,
\begin{equation}
\begin{aligned}
  F_{w}=&\int u \mathcal{M}_{j+1/2}(t,u,\mathbf{\xi}) \mathbf{\psi} du d\xi +\sum_k \int u_k \mathcal{F}_{j+1/2}(t,u_k,\mathbf{\xi}) \mathbf{\psi} d\mathbf{\xi}  \\
  =& \left(
  \begin{array}{c}
   \rho U \\
   \rho U^2+P-\frac43\mu U_x\\
    (\rho E+P)U-\frac43\mu U_xU- \kappa T_x \\
     \end{array}
      \right)_{j+1/2}+{O(\tau\Delta t, \Delta t^2, \Delta x^2)},
\end{aligned}
\end{equation}
 with the viscosity coefficient $\mu=\tau_{r}P$ and heat conduction coefficient $\kappa =C_p\tau_{r} P/\mbox{Pr}$.
From the above formulation, we can observe that UGKS is a suitable scheme for the NS equations,
   \begin{equation}
     \left\{
     \begin{aligned}
       &\frac{\partial \rho}{\partial t}+\frac{\partial (\rho U) }{\partial x}=
       O(\tau\Delta t, \Delta t^2, \Delta x^2),\\
       &\frac{\partial (\rho U) }{\partial t}+\frac{\partial}{\partial x} (\rho U^2+P-\frac43\mu U_x)=
       O(\tau\Delta t, \Delta t^2, \Delta x^2),\\
       &\frac{ \partial (\rho E)}{\partial t}+\frac{\partial}{\partial x} ((\rho E+P)U-\frac43\mu UU_x-\kappa T_x)=
       O(\tau\Delta t, \Delta t^2, \Delta x^2).
     \end{aligned}
     \right.
   \end{equation}

{In summary, when it comes to the hydrodynamic regime, the UGKS has the following properties:
   \begin{enumerate}
     \item In hydrodynamic regime, if the flow is well resolved under a fine mesh and $\Delta t$ is comparable to $\tau_{r}$, the scheme approximates the NS equations with a truncation error of $O(\tau_{r}\Delta t)$.
         If $\Delta t \gg \tau_r$, the scheme approximates the NS solution with a dominating error of $O(\Delta t^2)$.
         The time step of UGKS is not limited to be smaller than $\tau_r$ when it is applied in the NS regime.\\
     \item When $\Delta t \gg \tau_r$, the numerical flux of UGKS will not be sensitive to the initial distribution function at the beginning of each time step as well as the time discretization of the collision term.
         Even if we use an explicit or implicit method to calculate the collision term or use nonlinear limiter to reconstruct the flow field, the dominant numerical errors in solving NS equation keep second order.
         More specifically, let's investigate the time averaged cell interface flux of UGKS,
         \begin{displaymath}
           \begin{aligned}
            \tilde{F}_{ugks}=&\frac{1}{\Delta t}\int_0^{\Delta t} u f_{j+1/2}(t) dt\\
            =&\frac{1}{\Delta t}u[\tau_{r}(1-\mathrm{e}^{\Delta t/\tau_{r}})(H(u)f^l_{i+1/2}+(1-H(u))f^{r}_{i+1/2})\\
            &+\tau_{r} (\tau_{r} (\mathrm{e}^{-\Delta t/\tau_{r} -1}-1)+\Delta t \mathrm{e}^{-\Delta t/\tau_{r}})u(H(u)\sigma^l+(1-H(u))\sigma^r)\\
            &+(\Delta t -\tau_{r} (1-\mathrm{e}^{-\Delta t /\tau_{r}}))M_0\\
            &+\tau_{r} (2\tau_{r} (1-\mathrm{e}^{-\Delta t/\tau_{r}})-\Delta t (1+\mathrm{e}^{-\Delta t/\tau_{r}}))au M_0\\
            &+(\frac12 (\Delta t)^2+\tau_{r}(\tau_{r}(1-\mathrm{e}^{-\Delta t/\tau_{r}}-\Delta t))){\tilde A} M_0].
          \end{aligned}
    \end{displaymath}
     When Knudsen number approaches to zero with $\Delta t\gg\tau_{r}$, the above numerical flux goes to
      \begin{displaymath}
          \begin{aligned}
               \tilde{F}_{ugks}=& u \left[ M_0(1-\tau_{r}(au+{\tilde A})+\frac12\Delta t {\tilde A})  \right.\\
               & \left. +\frac{\tau_{r}}{\Delta t}[H(u)f_i+(1-H(u))f_{i+1}-M_0]+O(\tau_{r}^2) \right],
      \end{aligned}
      \end{displaymath}
   which shows that the numerical flux of the UGKS will not be sensitive to the initial reconstruction when time step is much larger than the local mean free time.
   The numerical flux is mainly contributed from the integration of the equilibrium state, which presents a NS flux. The initial term decays with $\tau_r /\Delta t$.
   The nonlinear limiter is to introduce a kinematic dissipation of $O(\Delta x^2)$ in the initial flow reconstruction \cite{xu2001}.
   Unfortunately, in the hydrodynamic limiting case, many other AP schemes will evaluate the cell interface flux from
   $[H(u)f_i+(1-H(u))f_{i+1}]$ only, where large numerical dissipation is intrinsically rooted from this flux vector splitting mechanism \cite{xu2001dissipative}.
   The last numerical example in the next section is about the laminar boundary layer at $\mbox{Re} = 10^5$ which
   is basically under such a situation. Accurate solution can be obtained from UGKS, but many other AP schemes with the upwind treatment for the interface flux will have difficulties here.
\end{enumerate}}

\subsubsection{Dynamics in the transition regime}

The above analysis presents the limiting governing equations of UGKS, i.e., the collisionless Boltzmann equation, the Euler equations, and the NS ones.
Between these limits, a smooth dynamical transition is practically obtained in UGKS with respect to the modeling scales of $\Delta x $ and $\Delta t$ and their ratios to the particle mean free
path and mean collision time.
Since the identified physics is closely related to the mesh size scale, it is difficult to figure out the underlying governing equations in the transition regime.
On the other hand, with the increasing of rarefaction, the degrees of freedom increases dramatically from a few conservative flow variables in the NS system to the infinite number of individual particle
movement. How to describe such a system theoretically and what kind of flow variables can be used to describe the non-equilibrium flow
with a continuous variation of degrees of freedom are basically unclear.
Actually, from theoretical point of view, the basic problem in any physical modeling is about the scale to identify the flow physics. Unfortunately, this has never been fully answered
in the traditional non-equilibrium thermodynamics research.
That is why all theoretical analysis of irreversible or extended thermodynamics are only limited to the near equilibrium regime.
However,
the UGKS, as a direct modeling, has no much difficulty to present a complete description of flow physics from equilibrium to non-equilibrium ones, because there is a specific modeling scale in UGKS.
If the algorithm itself can be considered as a kind of governing equation for the description of physical laws,
to understand and conduct theoretical analysis of the algorithm is as importance as the theoretical study of  partial differential equations.
Different from the derivation of the Boltzmann equation, the modeling scale for the validity of the NS equations has never been explicitly stated in the fluid mechanics research.
The current analysis only shows the
convergence of the UGKS to the NS equations under the limiting conditions, which may not be satisfied in practical computations.
The above analysis of so-called numerical error in terms of the NS equations can only provide a reference. Theoretically,
it is very hard to make any judgement about which one presents a more physically accurate solution
in the "continuum" flow regime, the NS or the UGKS.
The numerical examples in the next section demonstrate the dynamic differences from these two models in the near continuum regime.

\section{Numerical experiments}

\subsection{Sod shock tube problem}
\begin{figure}
\center
\includegraphics[width=6cm]{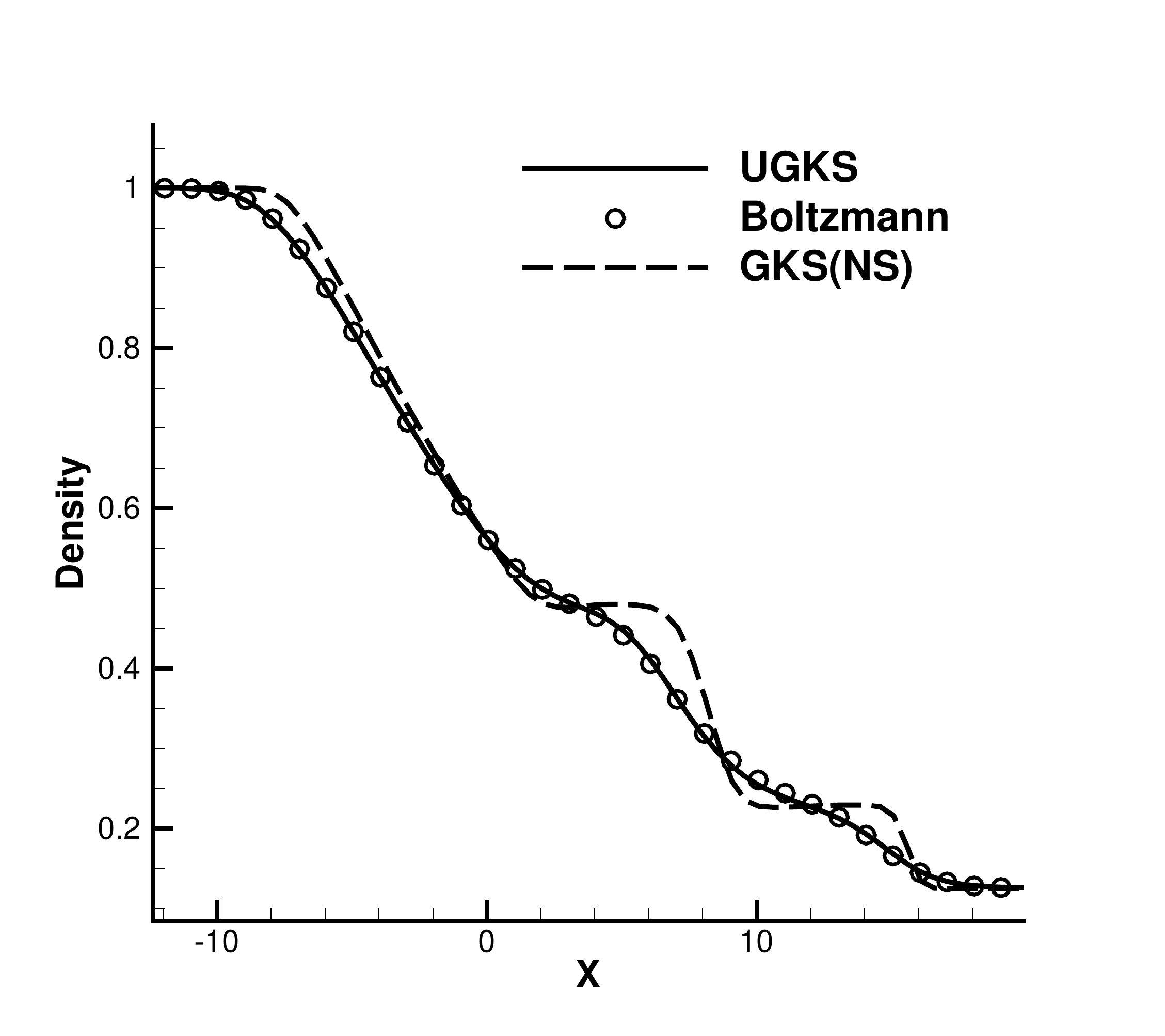}{a}
\includegraphics[width=6cm]{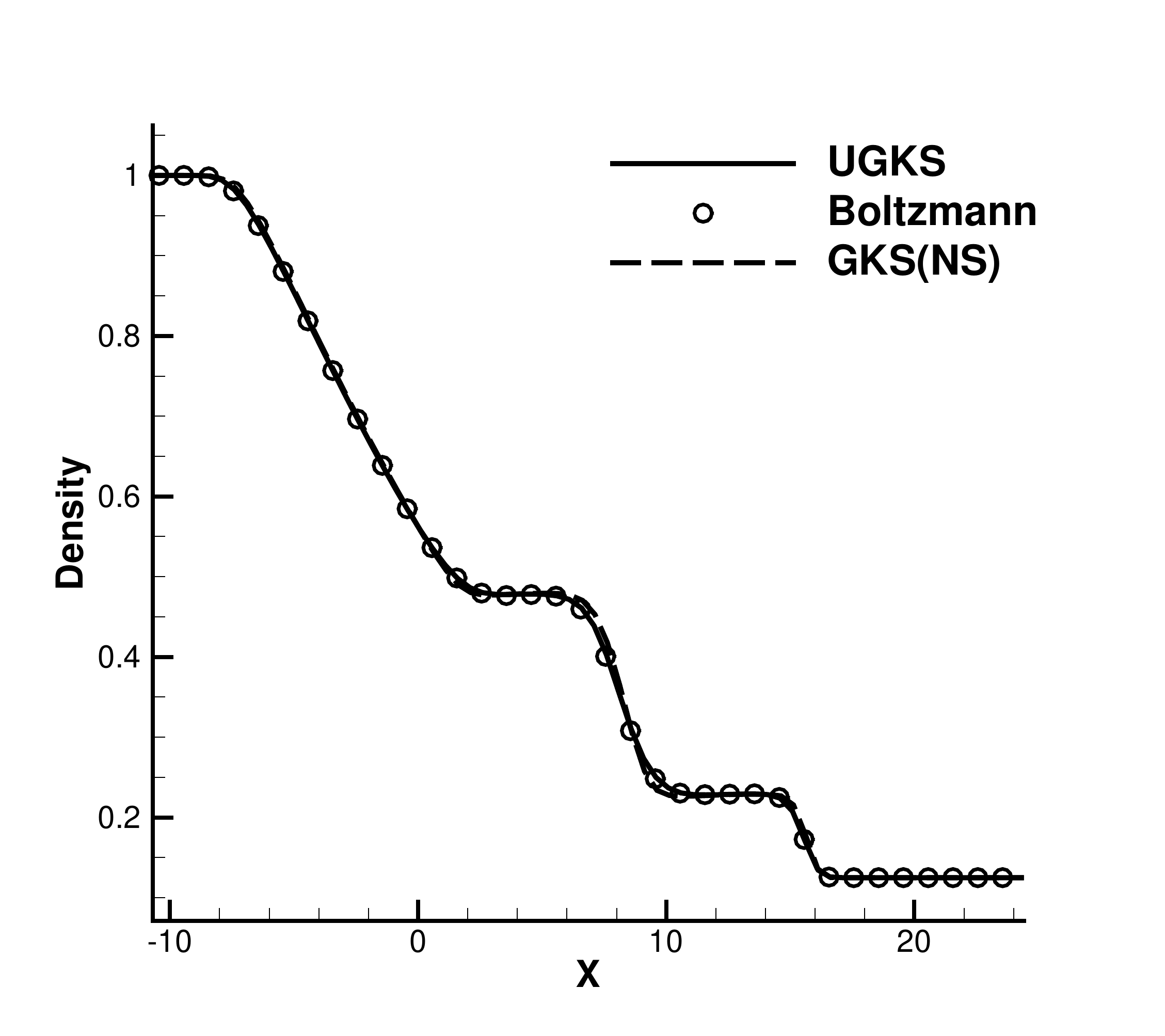}{b}
\includegraphics[width=6cm]{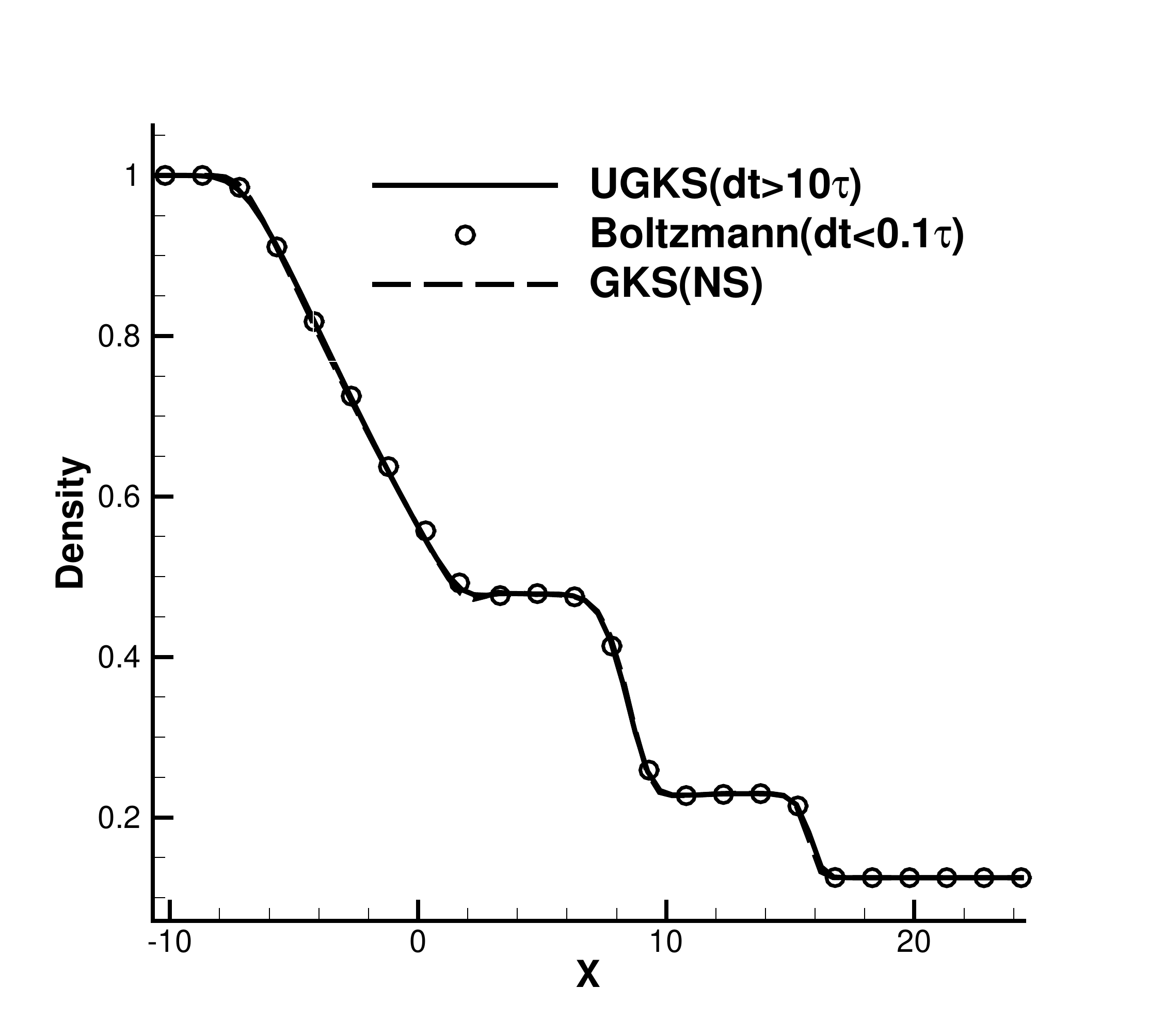}{c}
\caption{Density profiles of Sod test case. Line: UGKS solution; Circles: direct full Boltzmann solver; dash lines: GKS (NS solutions).
 (a) Kn$=0.1$; (b) Kn$=10^{-2}$; (c) Kn$=10^{-3}$. }
 \label{sod}
\end{figure}

We first calculate one dimensional Sod shock tube problem with Knudsen number in the range from $0.1$ to $10^{-3}$, to test the performance of UGKS in different flow regimes.
The gas medium is argon modeled by VHS model,
and the dimensionless initial condition is given by $\rho_l=1.0$, $U_l=0$, $T_l=1.0$, and $\rho_r=0.125$, $U_r=0$, $T_r=1.25$.
The results in Fig.(\ref{sod}), show that the UGKS is consistent with numerical Boltzmann solution
in rarefied regime and will converge to NS solution in hydrodynamic regime
with a time step being much larger than the particle mean free time ($\Delta t \geq 10 \tau$).
For the single scale methods, the NS equations cannot give physically consistent solution when Knudsen number becomes
large, and the direct Boltzmann solvers always require a small time step ($\Delta t \leq 0.1 \tau$) for a physical solution, even in the continuum regime.

\subsection{Normal shock structures}

The shock structure is one of the most important test case for the non-equilibrium flow.
In this calculation, we use non-uniform mesh in physical domain, such as a fine mesh in the upstream and a relative coarse mesh in the downstream.
In addition, local time step is used to get stationary solution.
In previous studies, shock structures have been calculated by UGKS with the Shakhov collision model only \cite{xuhuangshock}. The major difference between the
previous UGKS and DSMC solution is in the temperature profile around the upstream region, where the temperature from the UGKS rises earlier than that in the DSMC
for high Mach number shock wave.
For the density profiles, perfect match has been obtained between UGKS and DSMC solution.
Here in order to further improve the UGKS, the UGKS with the inclusion of the full Boltzmann collision term is tested.
The parameters to determine the switching between full Boltzmann and Shakhov model in the current UGKS
depends on the relative values of the the local time step and particle local mean free time.
In all shock structure calculations,
we set $t_{c}=0.4\tau_{r}$.
The test cases are mostly chosen from a recent paper about the full Boltzmann solver \cite{wu}, which provides easy comparison between the UGKS results and the full Boltzmann solutions.

\begin{figure}
\center
\includegraphics[width=6cm]{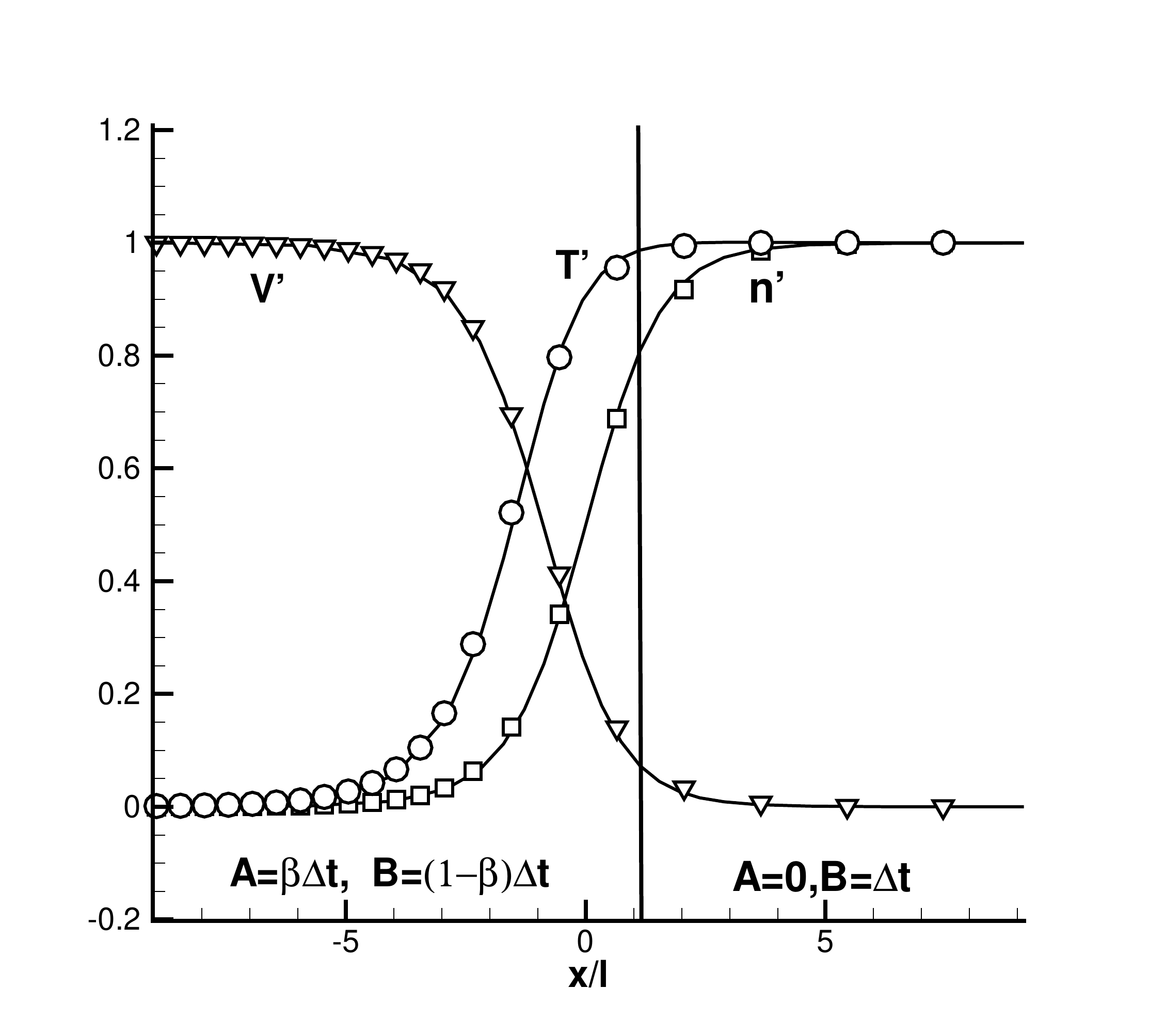}{a}
\includegraphics[width=6cm]{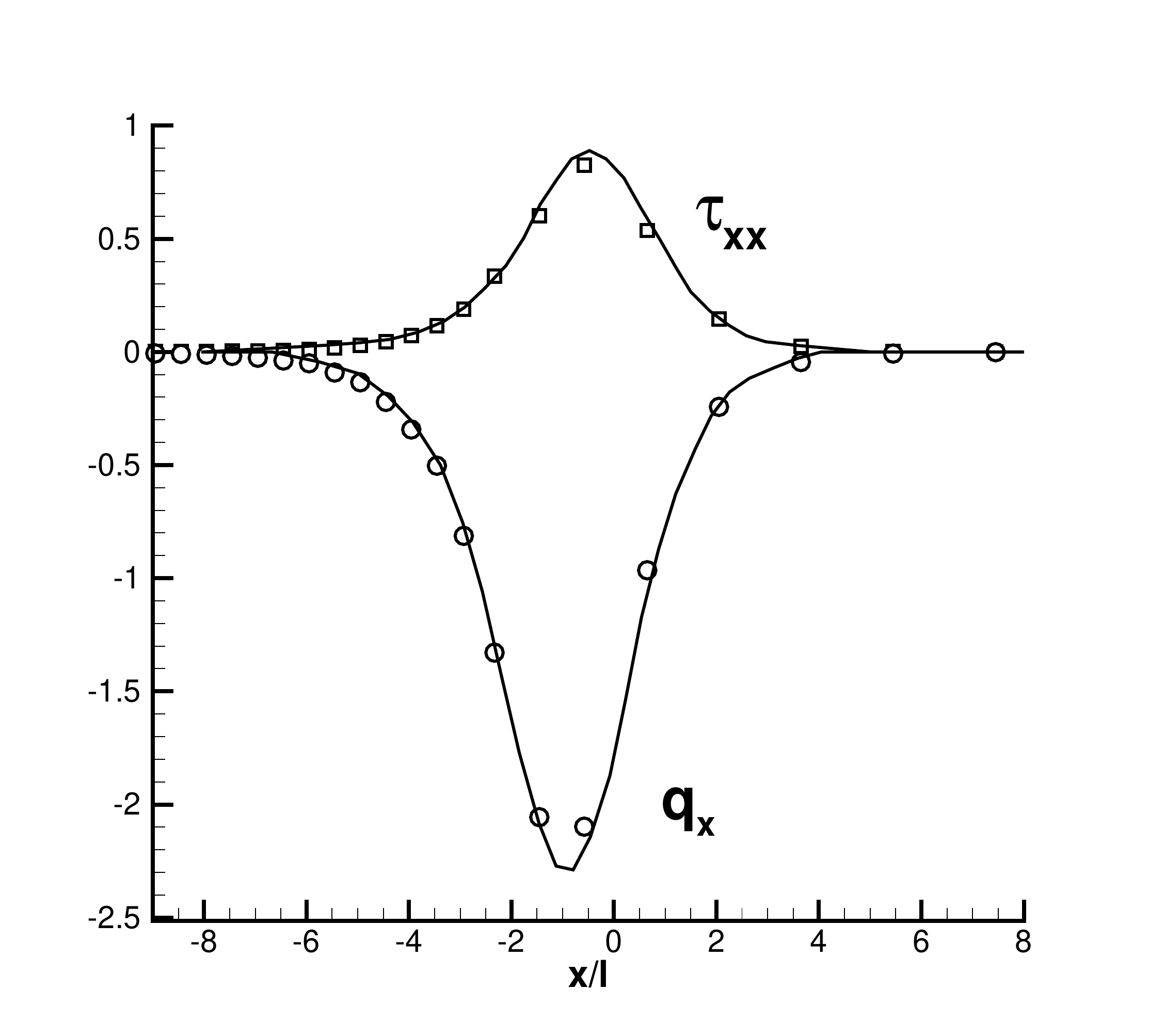}{b}
\label{ma3}
\caption{Shock structure computations with non-uniform mesh in the physical space at $\mbox{M}=3$ from the UGKS (symbols) and finite difference Boltzmann solution (lines) of \cite{ohwada1993structure}.
a: density, temperature and velocity; b: shear stress and heat flux. The vertical lines show the domain where the full Boltzmann collision term and Shakhov model are used.
}
\label{fig:shock-ohwada}
\end{figure}

We first consider the shock wave computation of hard sphere molecules.
Ohwada solved this problem by means of a finite difference method \cite{ohwada1993structure}.
Fig.(\ref{fig:shock-ohwada}) shows the shock structure, i.e., density, shear stress, and heat flux,  at Mach number $3$  from the UGKS (symbols) and reference solutions (lines).
The UGKS results get perfect match with the full Boltzmann solutions.
The vertical line in Fig.(\ref{fig:shock-ohwada}-a) shows the location for the switching between the combined full Boltzmann-Shakhov models and purely Shakhov model.
Based on this test, we can realize that the UGKS can use a large cell size in computation, especially in the downstream region. Even with the stretched cell size, accurate solutions can be obtained.

\begin{figure}
\centering
\includegraphics[width=8cm]
{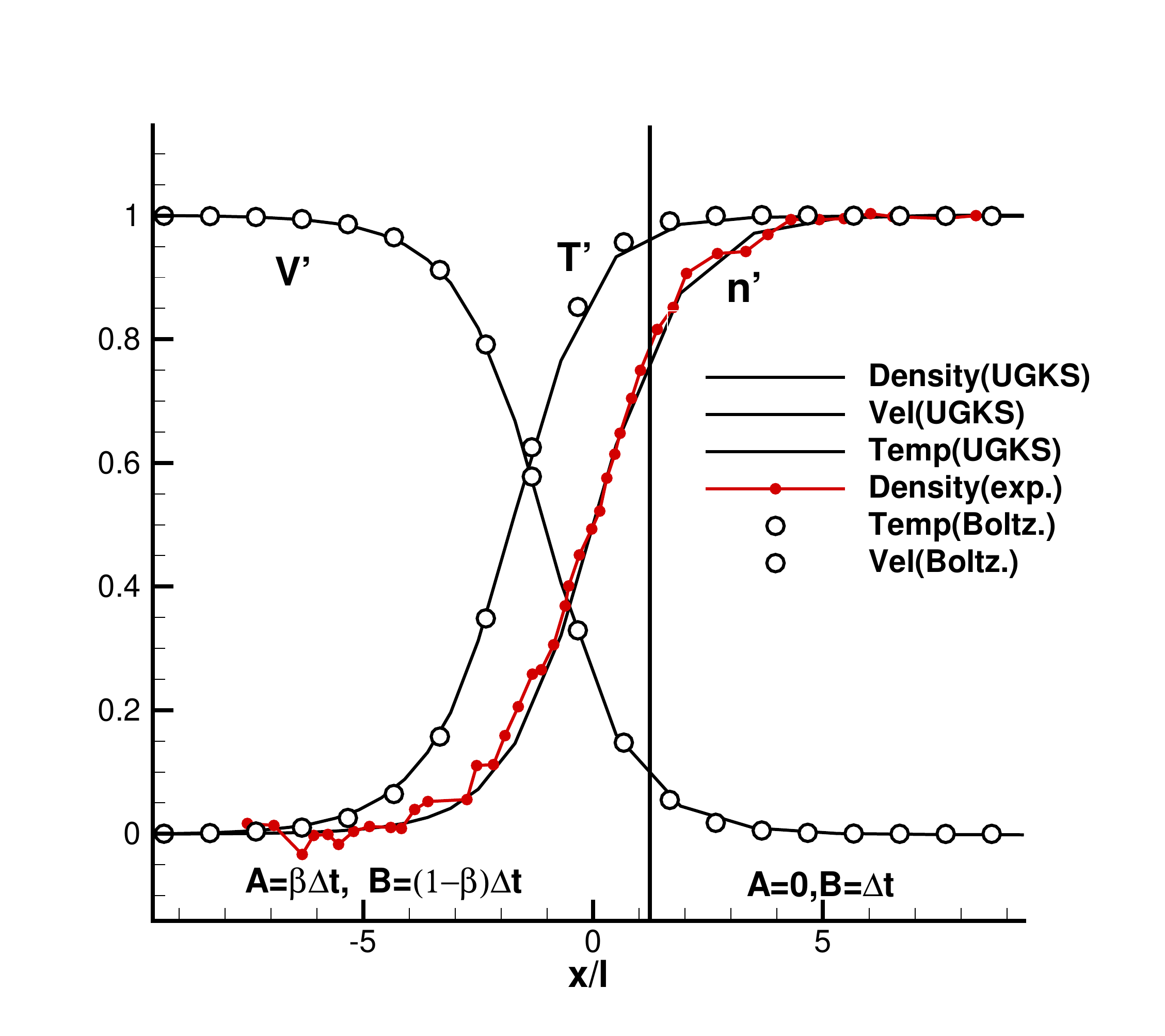}
\caption{Argon shock structure at $\mbox{M}=2.8$ from UGKS and experiment measurements \cite{kowalczyk2008numerical}.}
\label{ma2.8}
\end{figure}

\begin{figure}
\centering
\includegraphics[width=8cm]
{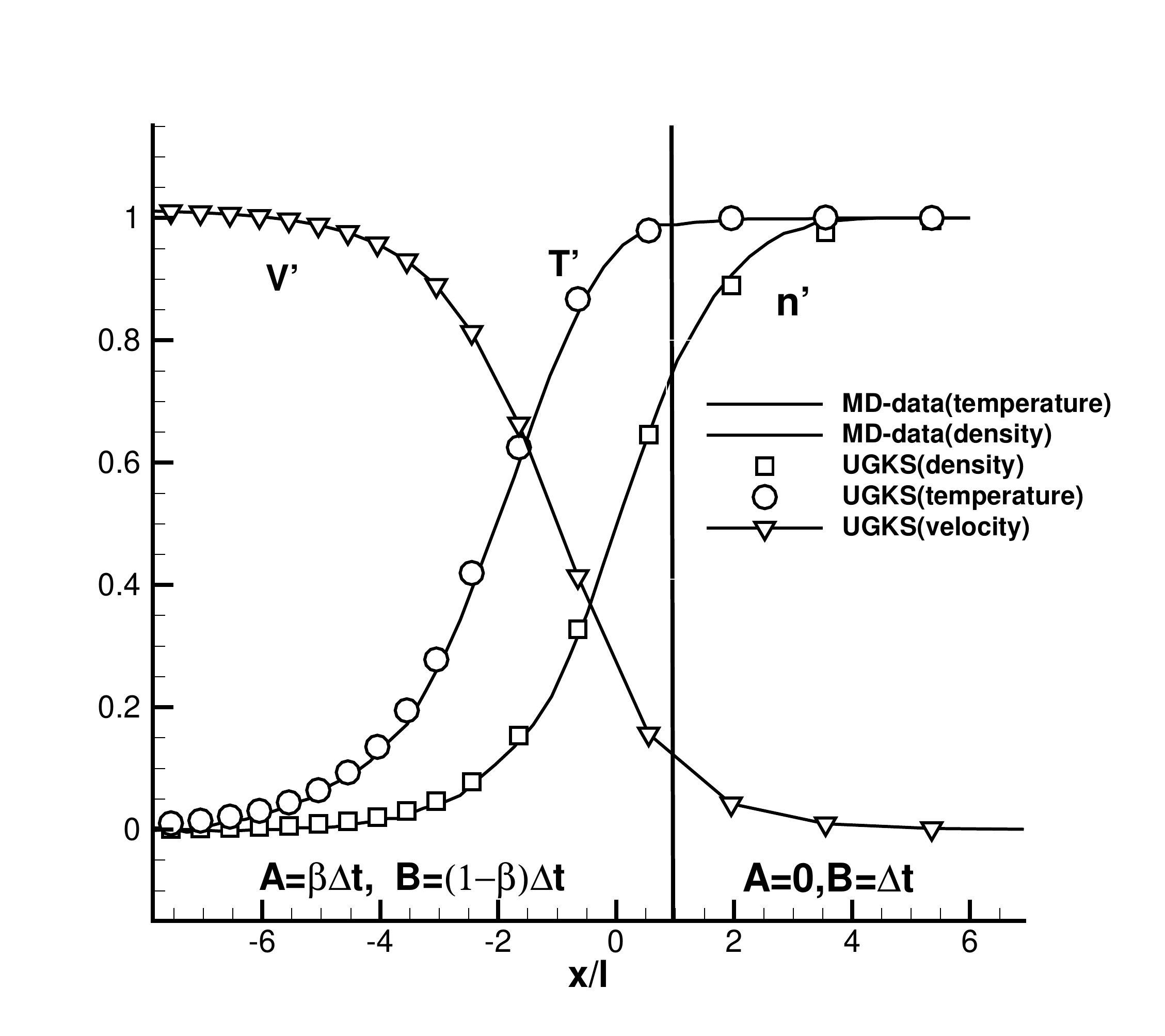}
\caption{Normalized number density, temperature and velocity distributions from UGKS (symbols) and MD solutions (lines) \cite{valentini2009large}.
} \label{ma5a}
\end{figure}

\begin{figure}
\center
\includegraphics[width=6cm]{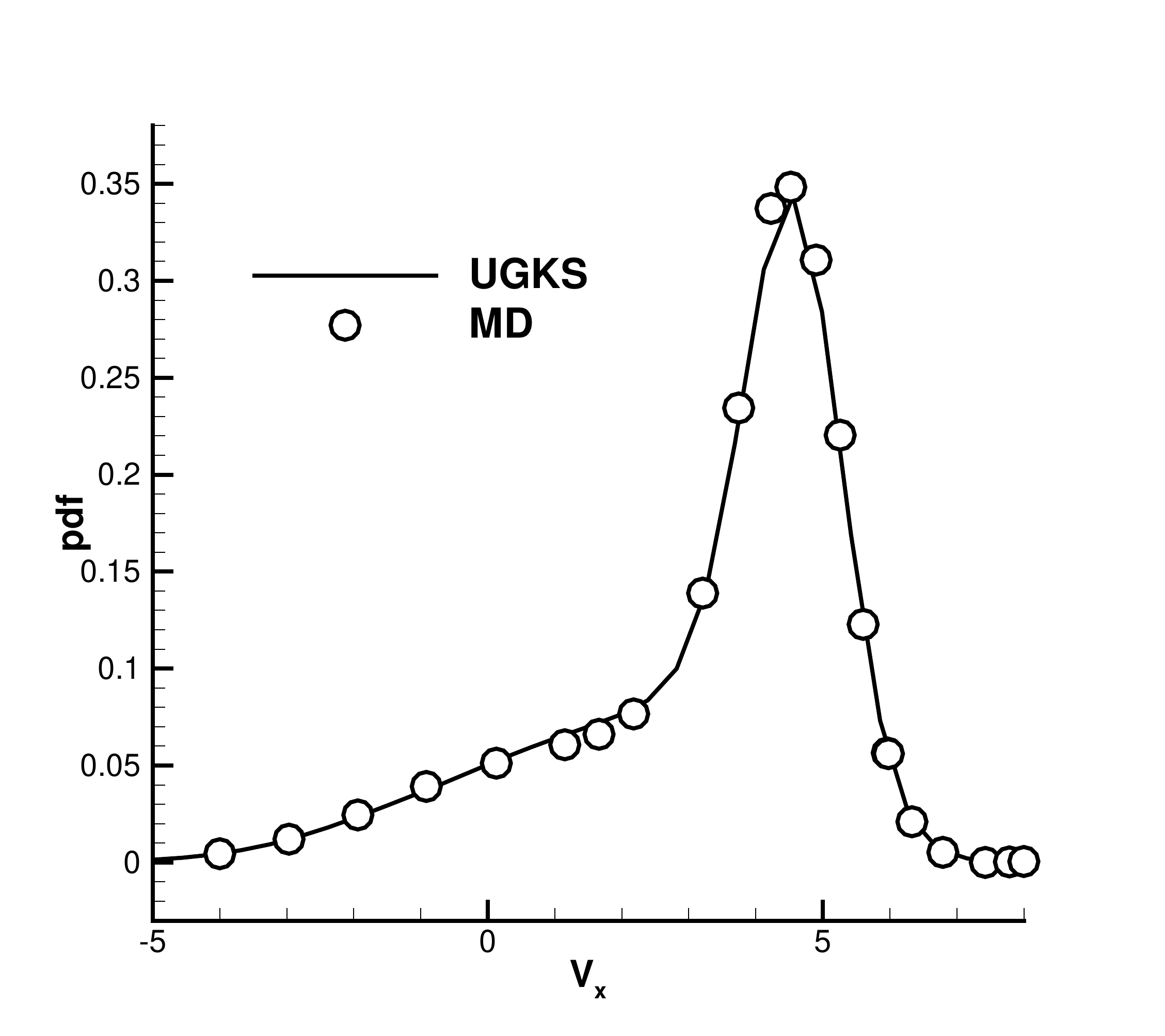}
\includegraphics[width=6cm]{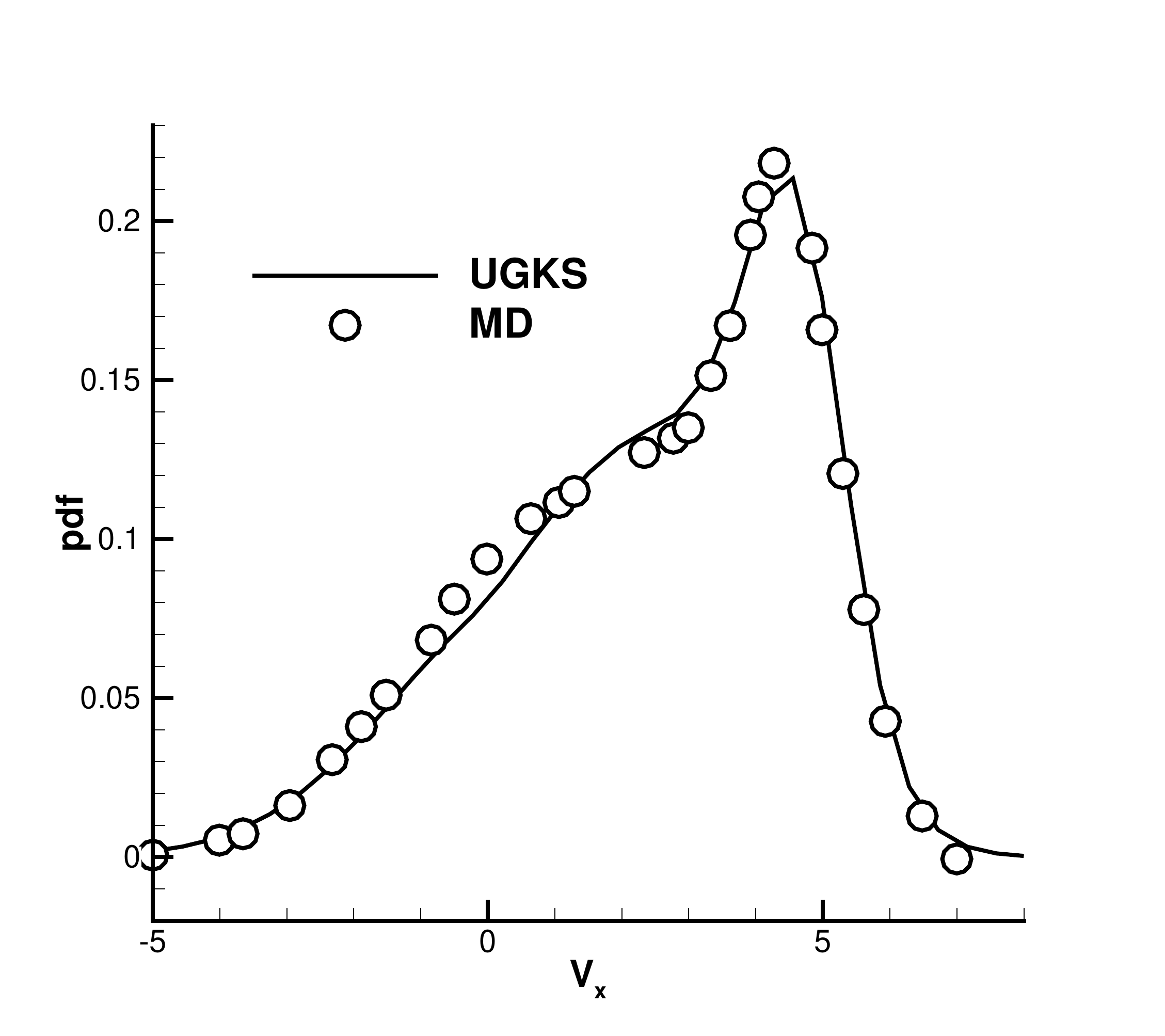}\\
\includegraphics[width=6cm]{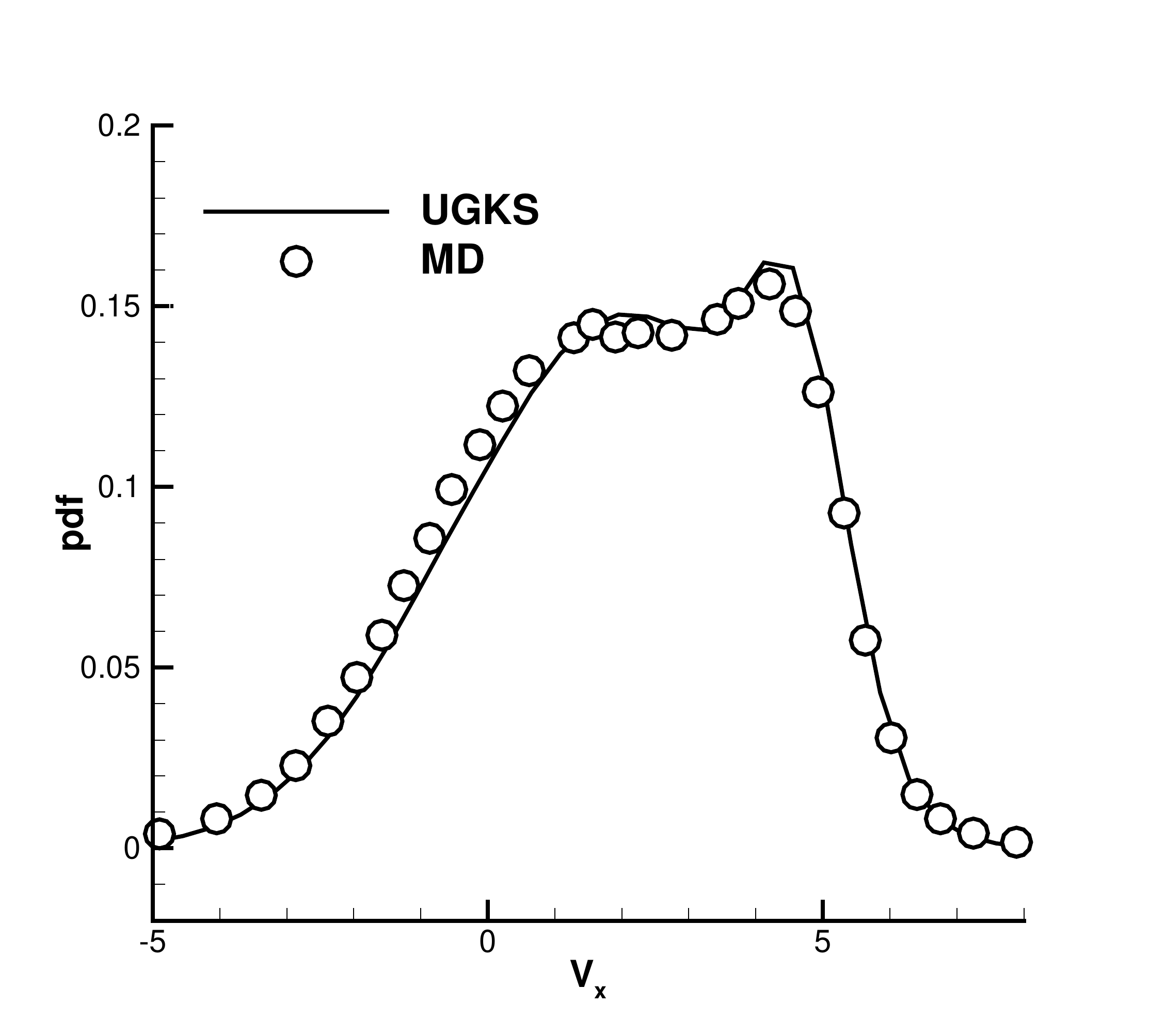}
\includegraphics[width=6cm]{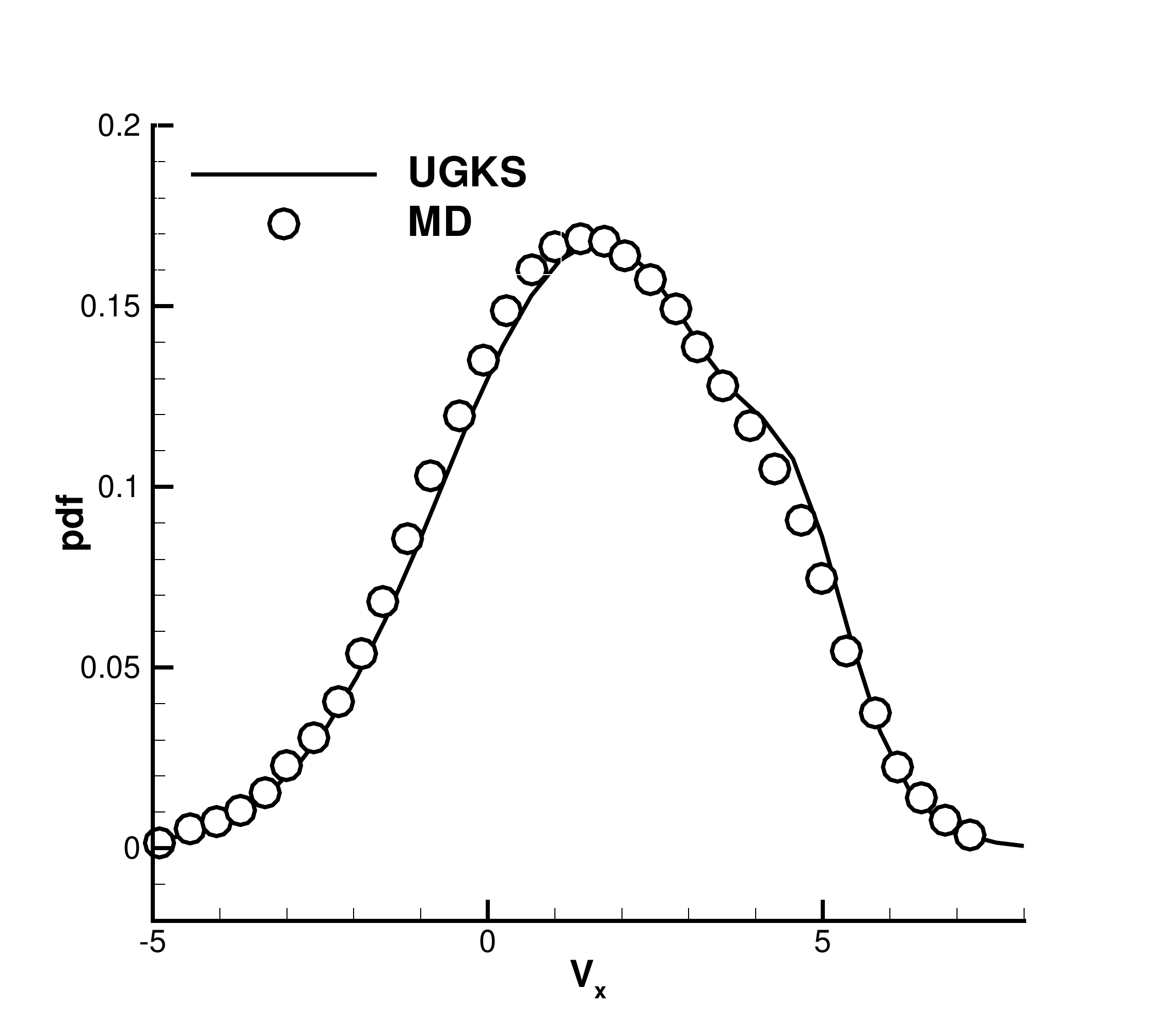}
\caption{
The distribution function $\int\int f dudw/n$ at locations of density $n'=0.151, 0.350, 0.511$, and $0.759$.
UGKS solutions (lines) and MD solution (symbols) \cite{valentini2009large}.
} \label{ma5b}
\end{figure}

Next we consider argon gas with L-J potential. We use the generalized anisotropic collision model to recover L-J potential which is expressed in Eq.\eqref{LJ}.
Fig.(\ref{ma2.8}) shows the shock structure of argon gas with L-J potential at Mach number $2.8$ from the UGKS and experiment measurement \cite{kowalczyk2008numerical}. For the shock wave of argon gas with L-J potential at $\mbox{M}=5$, we compare the UGKS solution with molecular dynamics simulation of \cite{valentini2009large}.
Fig.(\ref{ma5a}) and (\ref{ma5b}) present the shock wave structure and the distribution functions inside the shock layer.

\begin{figure}
\centering
\includegraphics[width=8cm]
{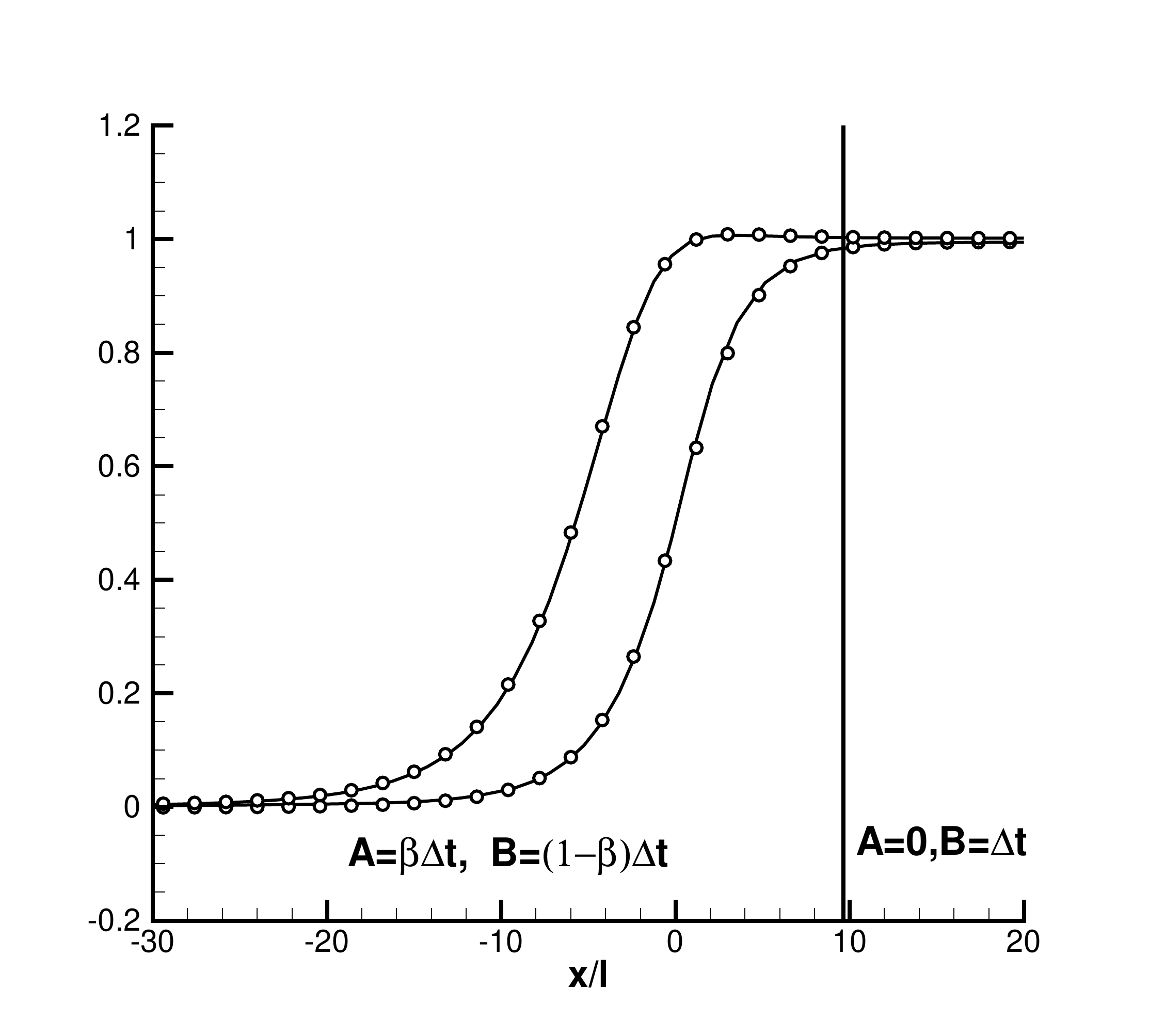}
\caption{Shock structure calculations at $\mbox{M}=6$ from the UGKS (symbols) and the full Boltzmann solution (lines). In UGKS, a non-uniform mesh and local time step are used.}
\label{ma6}
\end{figure}

The last shock structure calculation is the argon gas at $\mbox{M}=6$ from UGKS (symbols) with non-uniform mesh and
the full Boltzmann solution (lines) with a much refined mesh. The results are shown in Fig.(\ref{ma6}).
This shows that even with the variation of mesh size, the physical solution can be always captured by UGKS.

\subsection{Flow passing through a circular cylinder}

In order to test the performance of UGKS for two-dimensional high speed flow in various flow regimes,we calculate a flow of argon gas passing through a circular cylinder at Mach number 2 and Knudsen number $\text{Kn}=1.0$, $0.1$, $10^{-2}$, and $10^{-3}$ relative to cylinder radius. The gas medium is argon modeled by VHS model.

\begin{figure}
\center
\includegraphics[width=6cm]{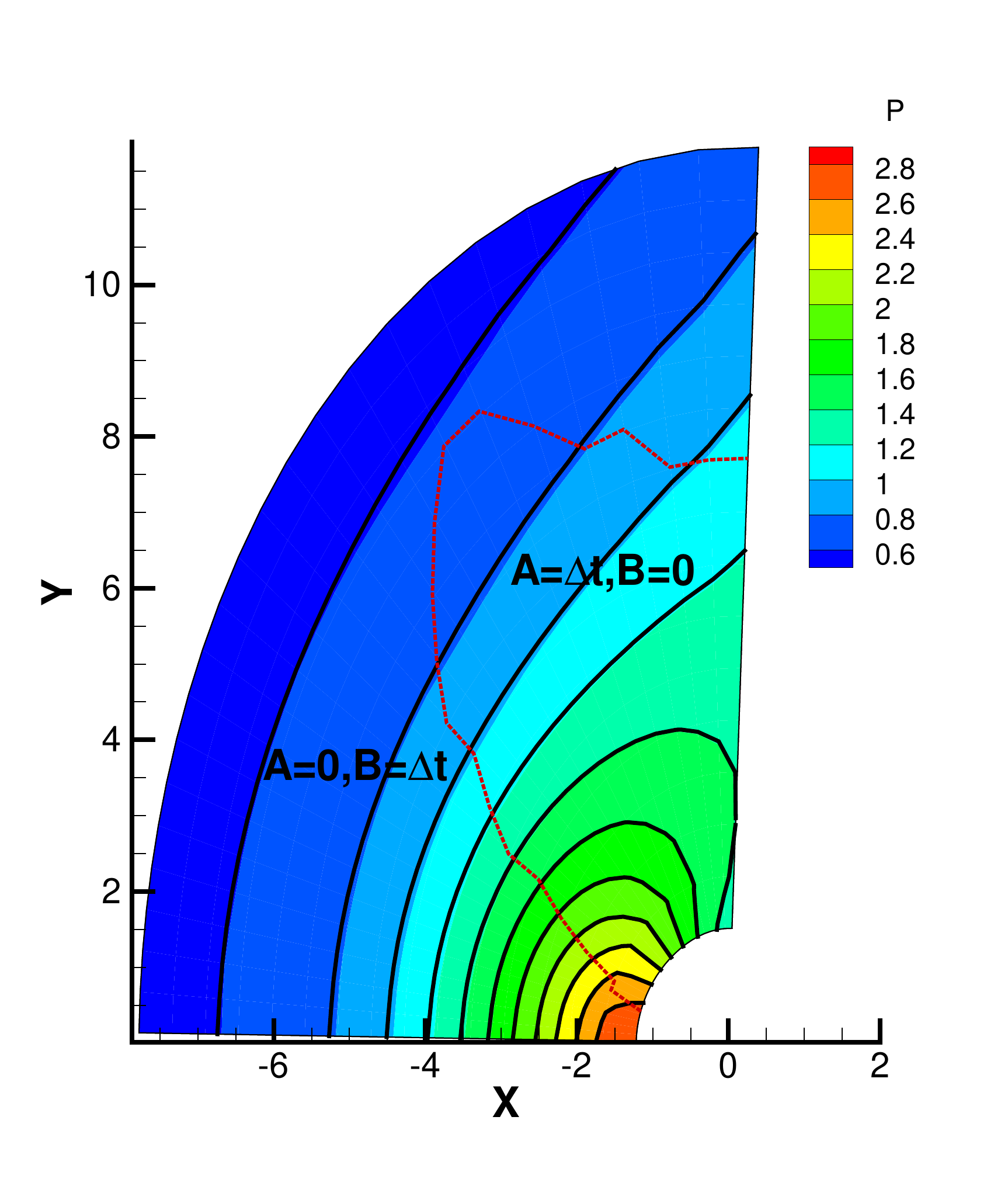}{a}
\includegraphics[width=6cm]{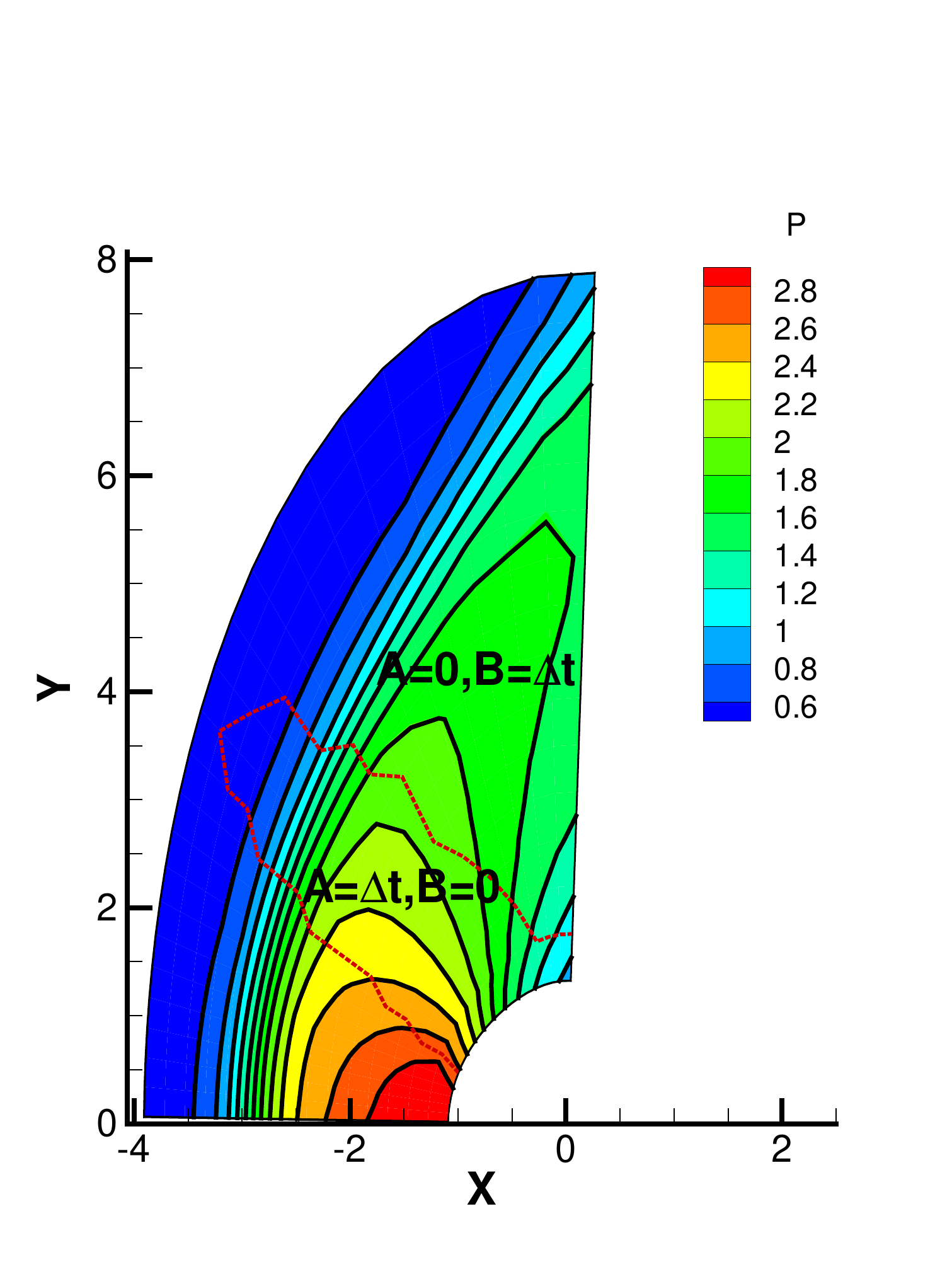}{b}\\
\includegraphics[width=6cm]{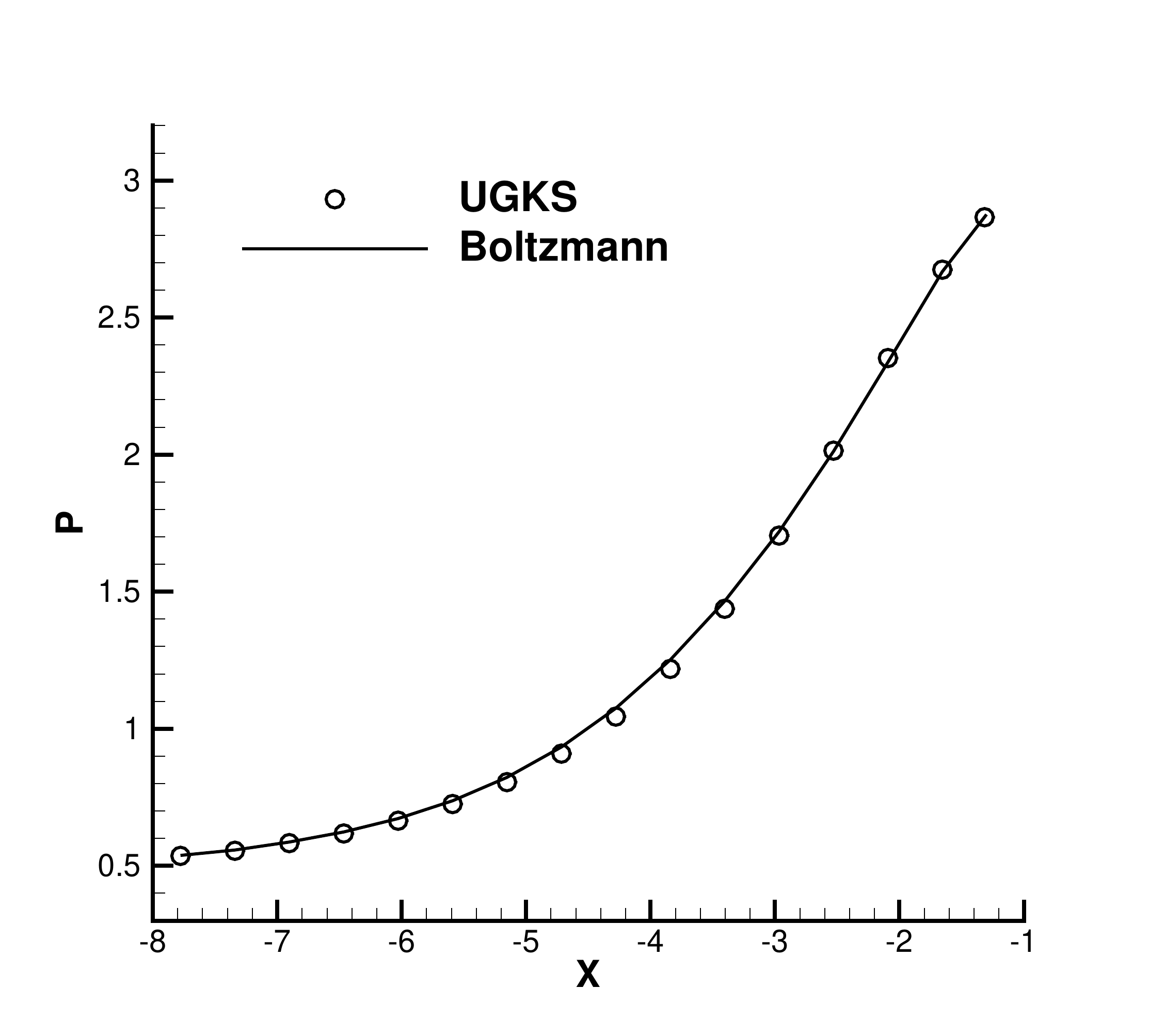}{c}
\includegraphics[width=6cm]{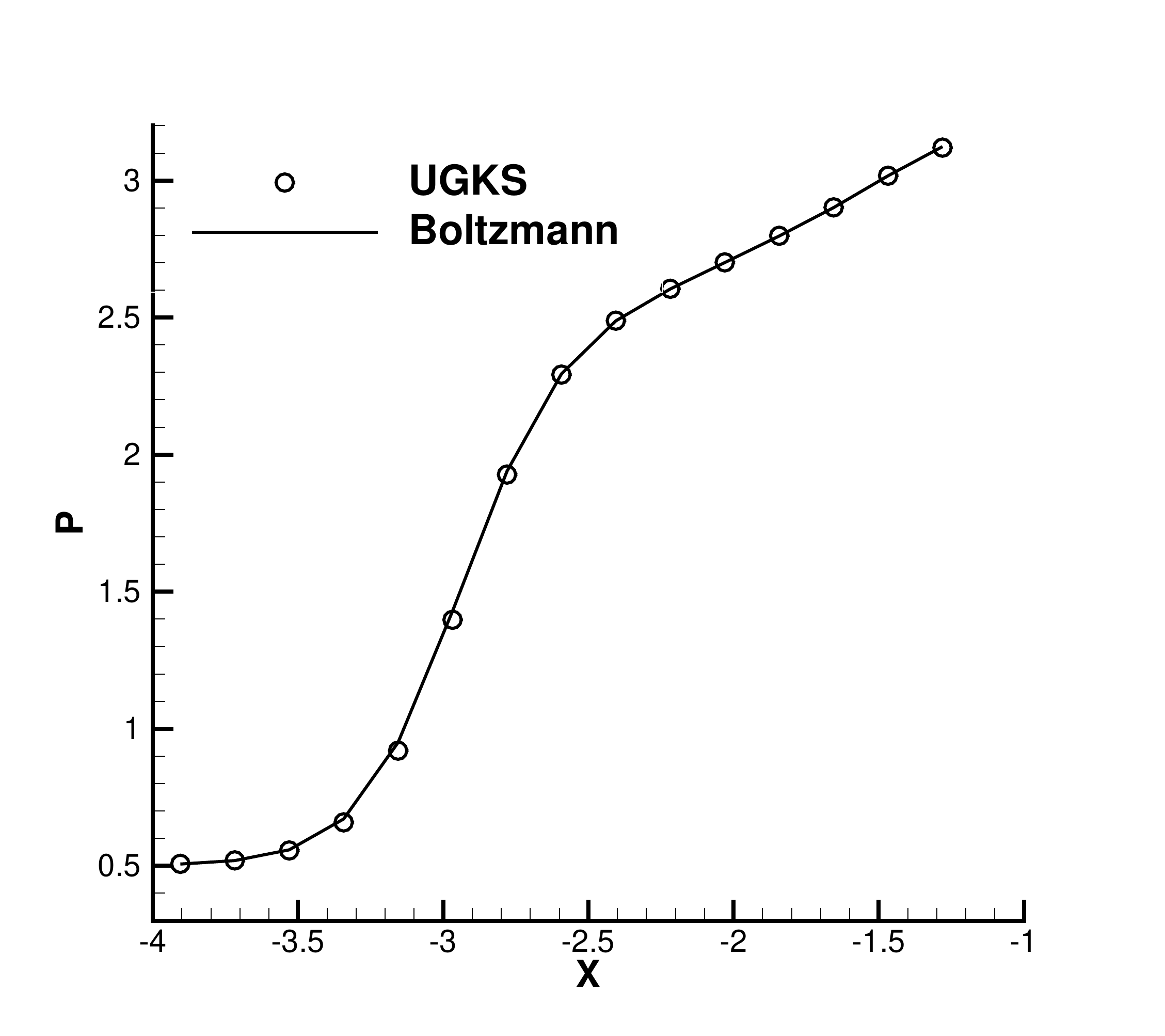}{d}
\caption{Flow passing through a circular cylinder at $M=2.0$.
a: pressure contours at Kn$=1.0$. Background: UGKS solution; Solid lines: direct Boltzmann solution;
b: pressure contours at Kn$=0.1$. Background: UGKS solution; Solid lines: direct Boltzmann solution;
c: pressure along central line in front of cylinder at Kn$=1.0$;
d: pressure along central line in front of cylinder at Kn$=0.1$.}
\label{cylinderb}
\end{figure}

For the case with $\text{Kn}=1.0$ and $0.1$, the velocity space $v\in[-8,8]^3$ is divided into $32$ equally spaced velocity grids in each direction and $24\times16$ grids are used in physical space. We fix the CFL number to be $0.5$ and use local time to reach steady state. {We set $t_{c}=1.2\times 10^{-3}\tau_r$ for $\text{Kn}=1.0$ and $t_{c}=3.4\times 10^{-3}\tau_r$ for $\text{Kn}=0.1$. Because $t_c<1/\sup\nu$, the parameters $A$ and $B$ are chosen as
\begin{equation}\nonumber
 (A,B)=
 \begin{cases}
 (\Delta t,0), &\Delta t<t_c^n,\\
 (0,\Delta t), &\Delta t\ge t_c^n.
 \end{cases}
\end{equation}}
The UGKS solutions are compared with the direct Boltzmann solver, as shown in Fig.(\ref{cylinderb}).

\begin{figure}
\center
\includegraphics[width=6cm]{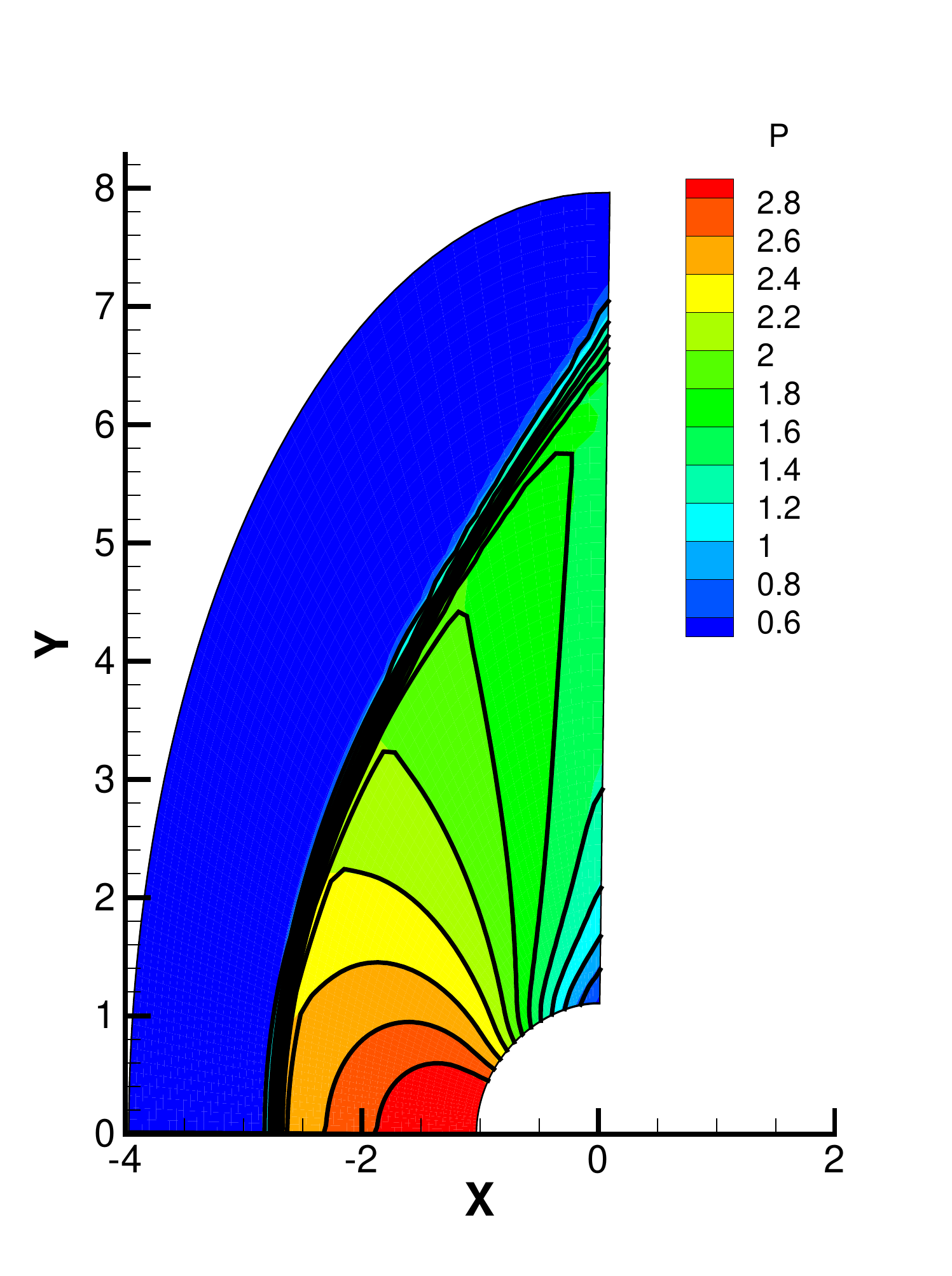}{a}
\includegraphics[width=6cm]{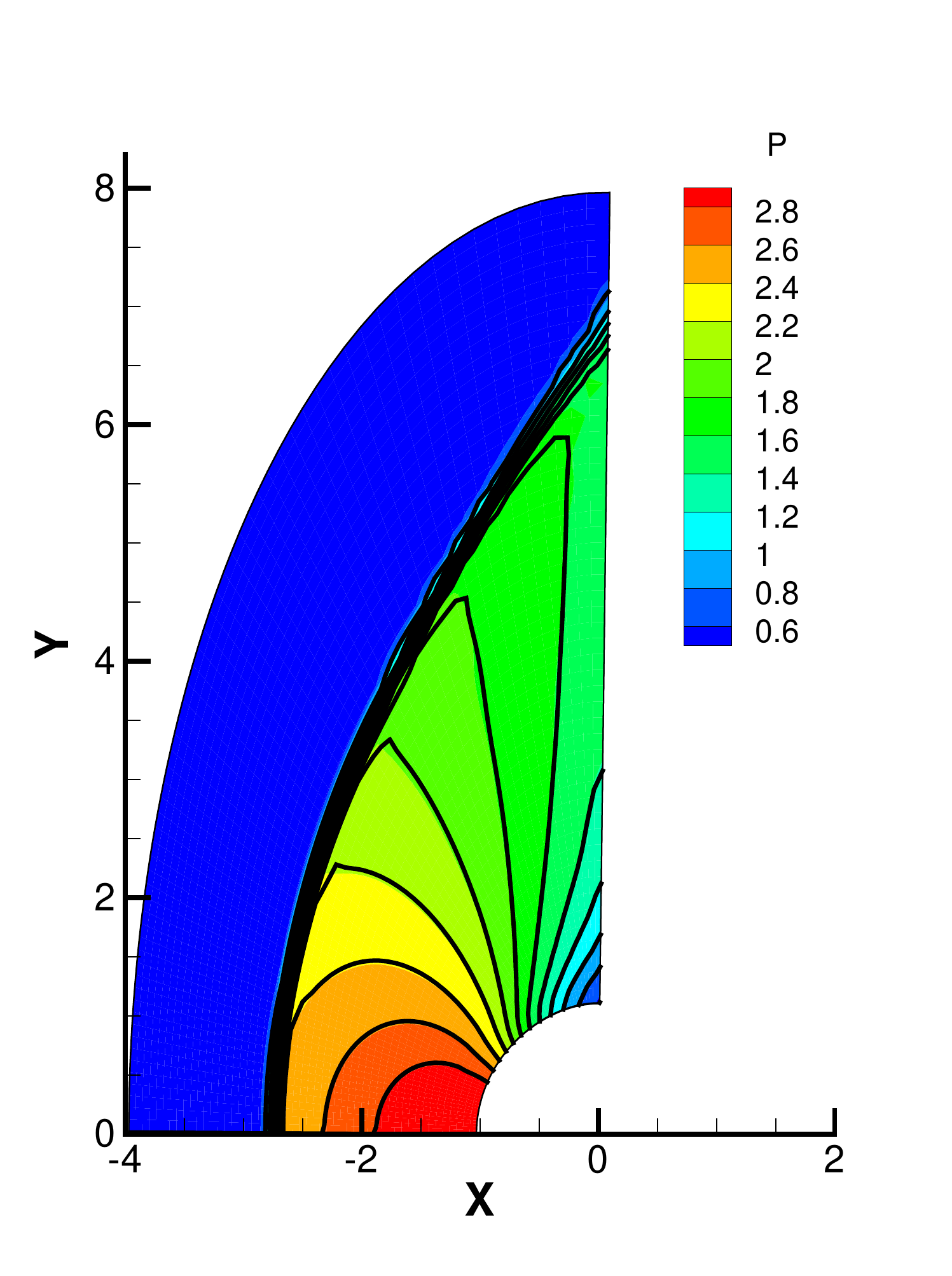}{b}\\
\includegraphics[width=6cm]{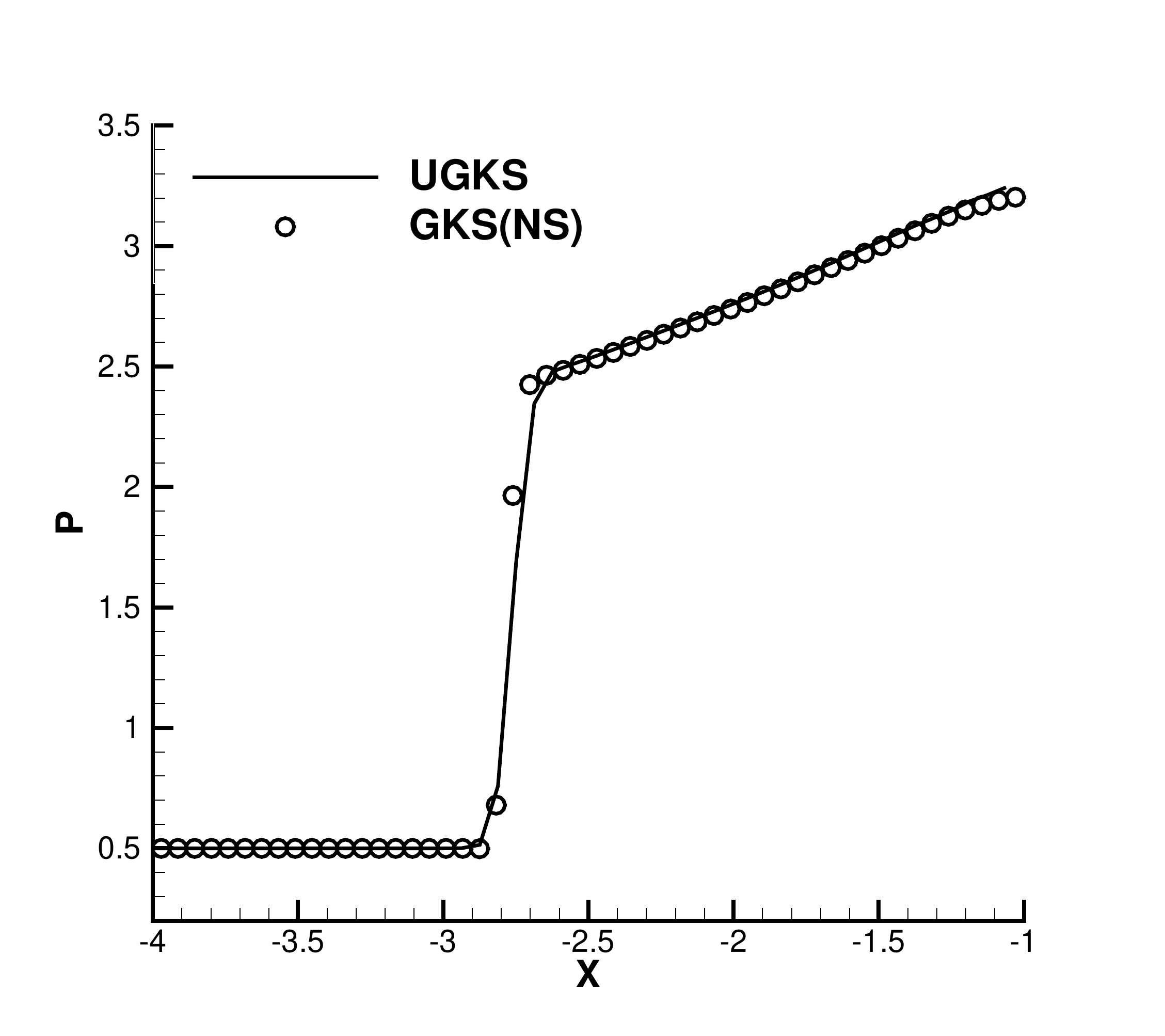}{c}
\includegraphics[width=6cm]{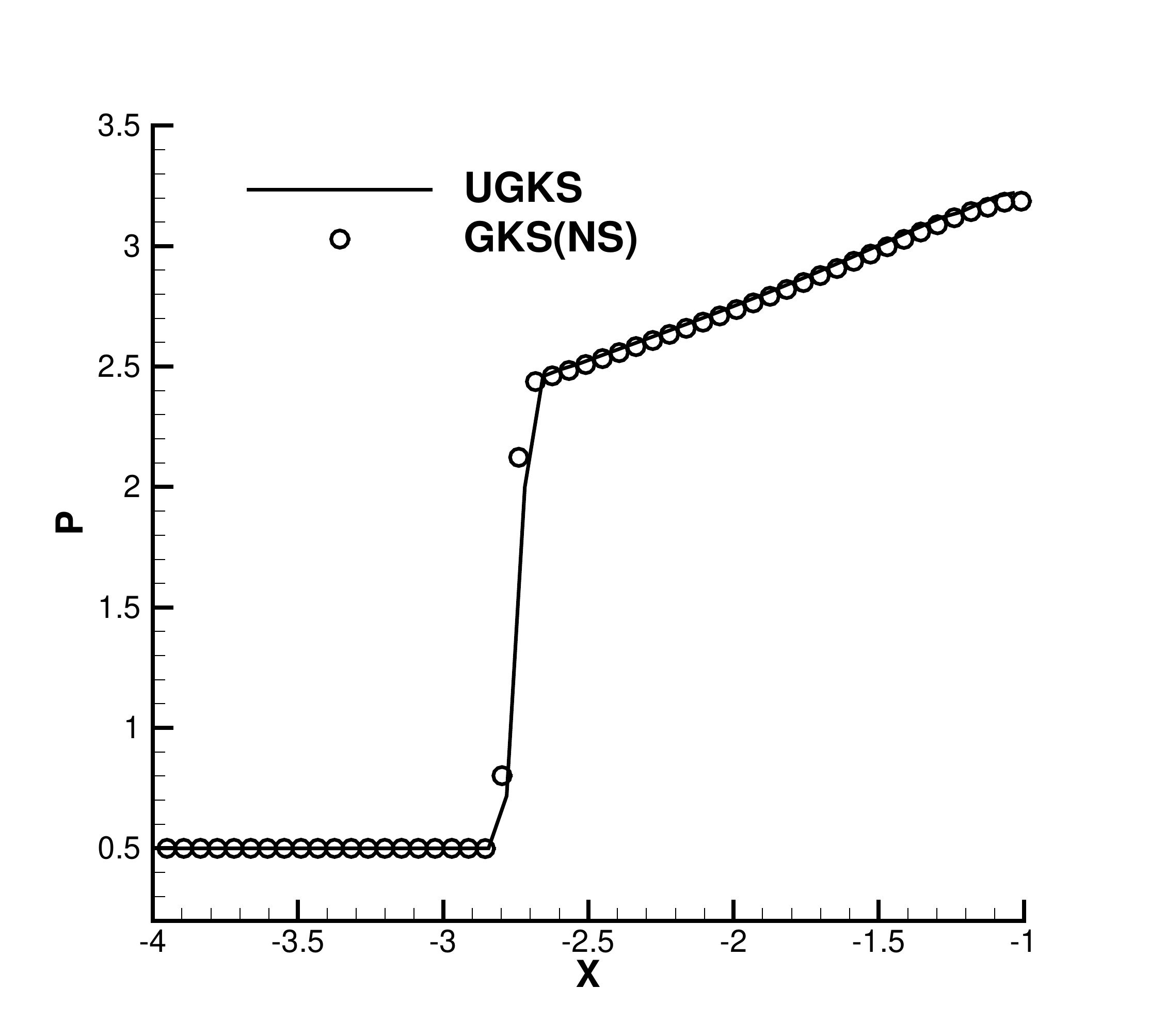}{d}
\caption{Flow passing through a circular cylinder at $M=2.0$.
a: pressure contours at Kn$=10^{-2}$. Background: UGKS solution; Solid lines: GKS (NS) solution;
b: pressure contours at Kn$=10^{-3}$. Background: UGKS solution; Solid lines: GKS (NS) solution;
c: pressure along central line in front of cylinder at Kn$=10^{-2}$;
d: pressure along central line in front of cylinder at Kn$=10^{-3}$.}
\label{cylindern}
\end{figure}

For the cases with $\text{Kn}=10^{-2}$, and $10^{-3}$, the velocity space $v\in[-8,8]^3$ is divided into $42$ equally spaced velocity grids in each direction and $64\times48$ grids are used in physical space. For these two cases, the implicit Shakhov collision operator is used in UGKS and the UGKS solutions are compared with the NS solutions by the GKS \cite{xu2001}, which are shown in Fig(\ref{cylindern}).
In these small Knudsen number cases, the UGKS becomes a shock capturing scheme for the NS solutions.

\subsection{Lid-driven cavity flow}

The 2-d lid-driven cavity flow is used for the study of flow physics in the whole flow regimes.
In the following calculations, the gas medium consists of argon modeled by VHS.
The wall temperature is kept the same as reference temperature of $T_w=T_0=273K$, and the up wall velocity is kept fixed at
$U_w=50m/s$. Maxwell's diffusion boundary condition  with full accommodation is used at the boundaries.
In the physical space, a non-uniform mesh is used in
order to identify the flow structure with different resolution.
The grid point follows
\begin{equation}
  x=(10-15s+6s^2)s^3-0.5, ~~ s=(0,1,...,N)/2N,
\end{equation}
 in $x$-direction and similar formula is used in the $y$-direction.

\begin{figure}
\center
\includegraphics[width=6.0cm]{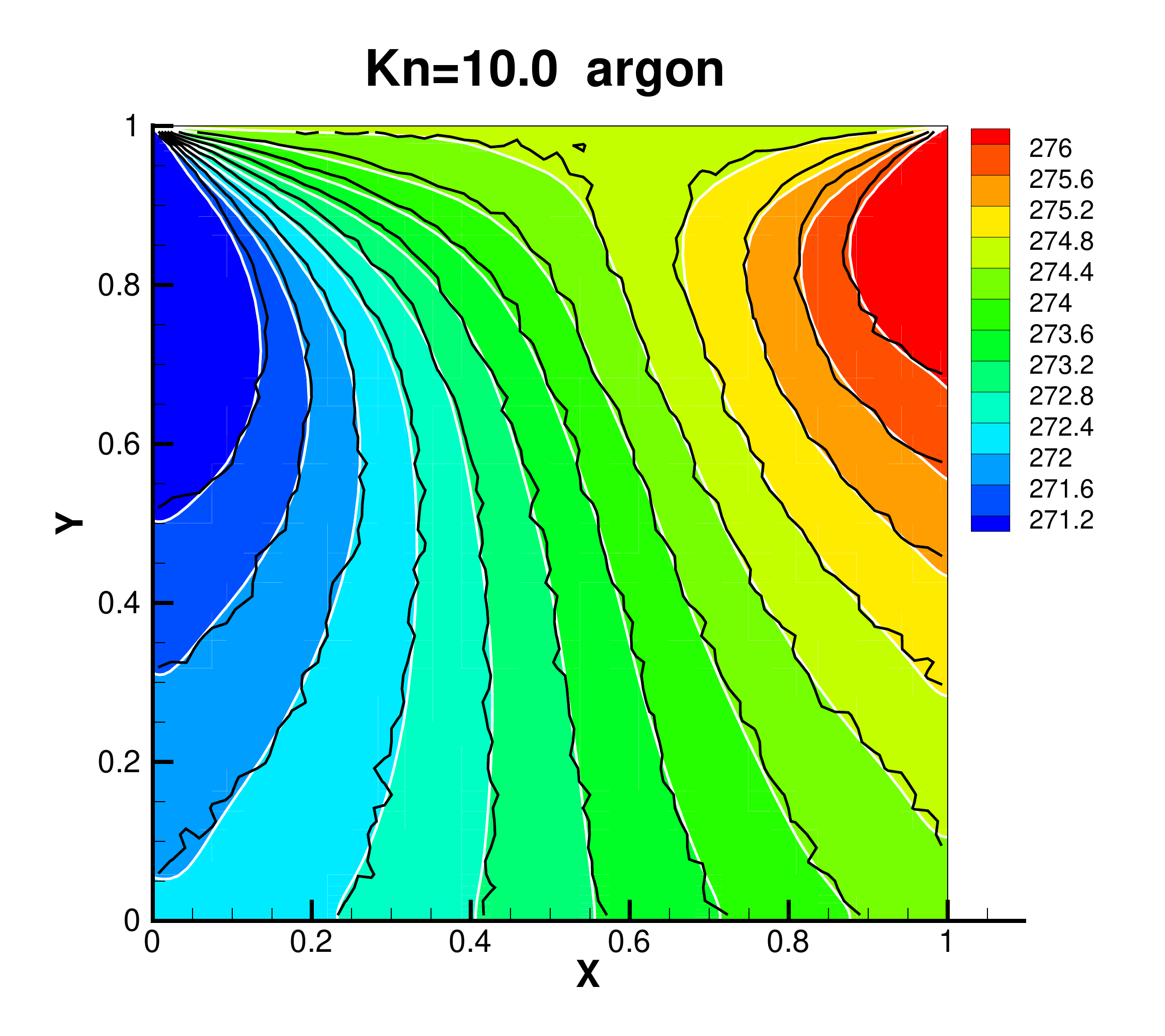}{a}
\includegraphics[width=6.0cm]{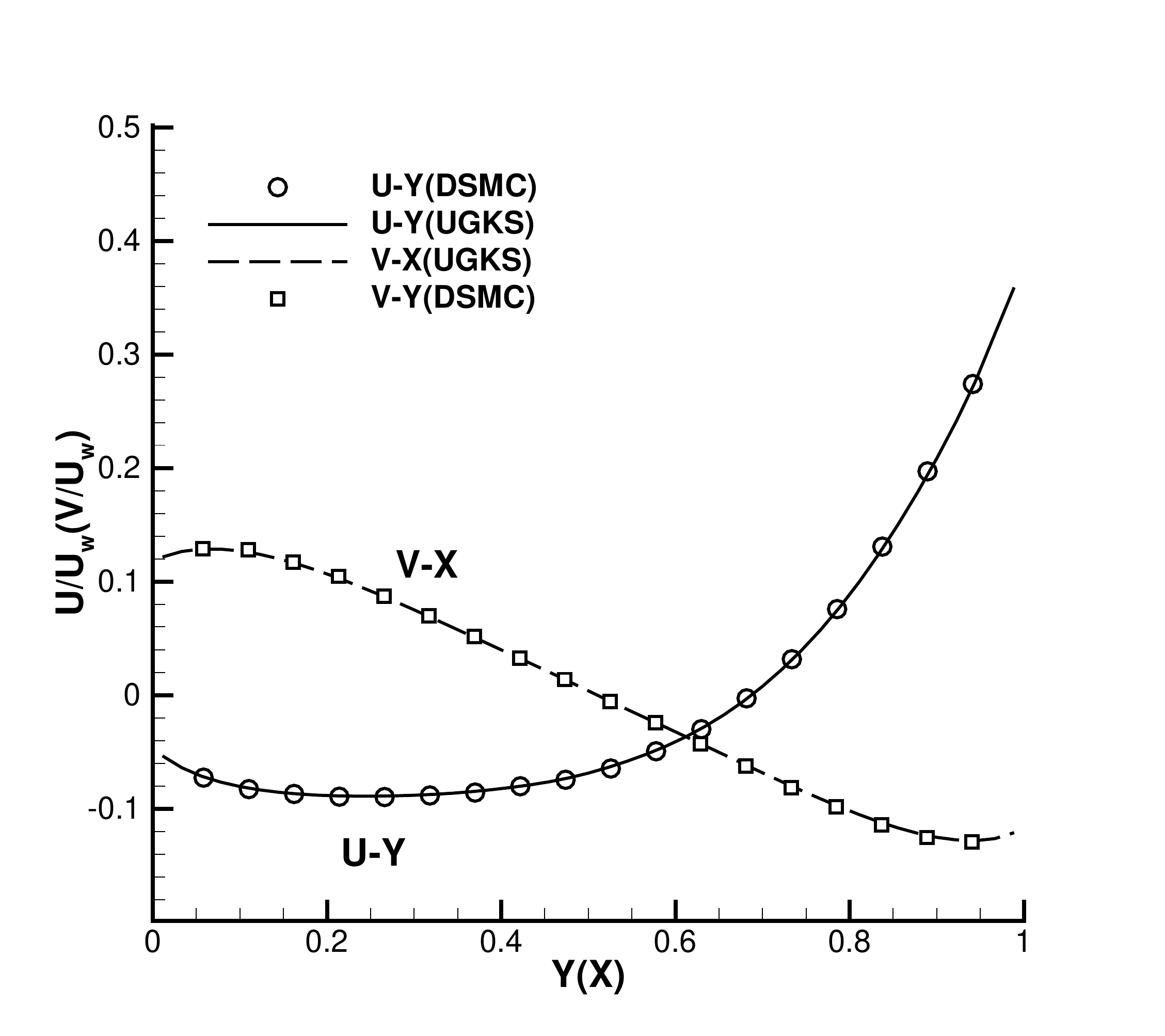}{b}
\caption{Cavity flow at \mbox{Kn}=10.
 (a) temperature contours, black lines: DSMC, white lines and background: UGKS;
 (b) U-velocity along the central vertical line and V-velocity along the central horizontal line, circles: DSMC, line:UGKS.}
 \label{k10}
\end{figure}

\begin{figure}
\center
\includegraphics[width=6.0cm]{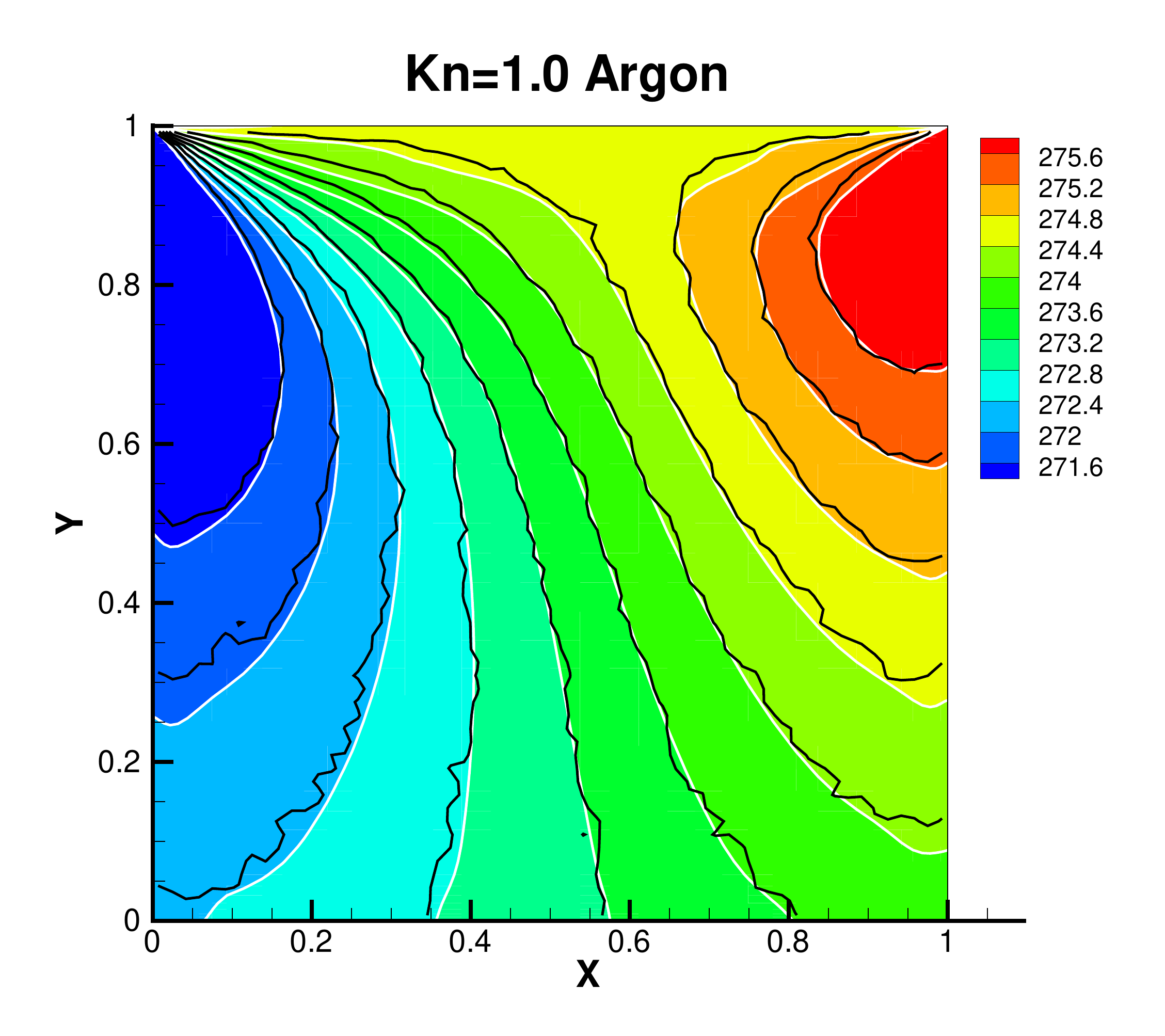}{a}
\includegraphics[width=6.0cm]{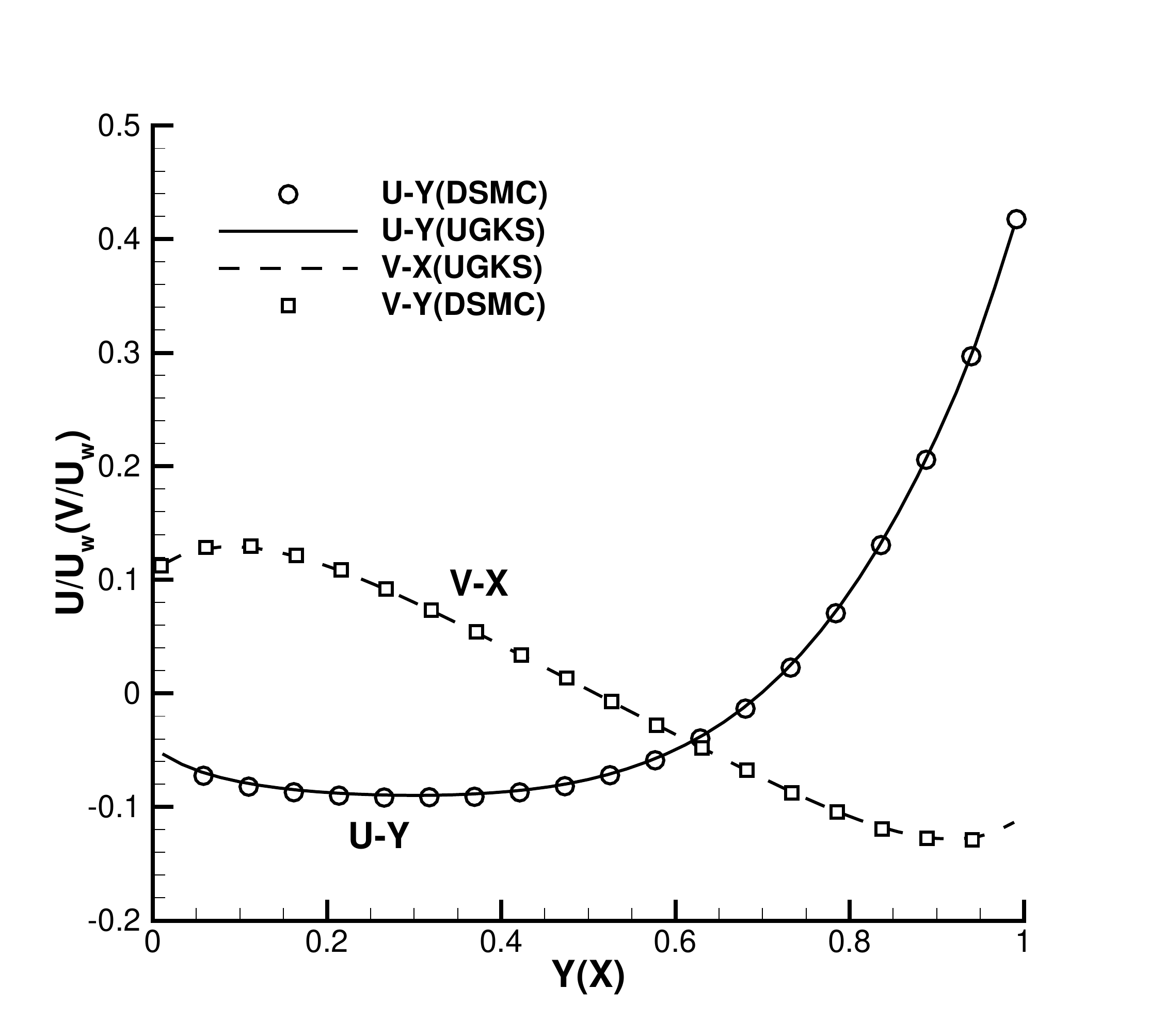}{b}
\caption{Cavity flow at \mbox{Kn}=1.
 (a) temperature contours, black lines: DSMC, white lines and background: UGKS;
 (b) U-velocity along the central vertical line and V-velocity along the central horizontal line, circles: DSMC, line:UGKS.}
 \label{k1}
\end{figure}

\begin{figure}
\center
\includegraphics[width=6.0cm]{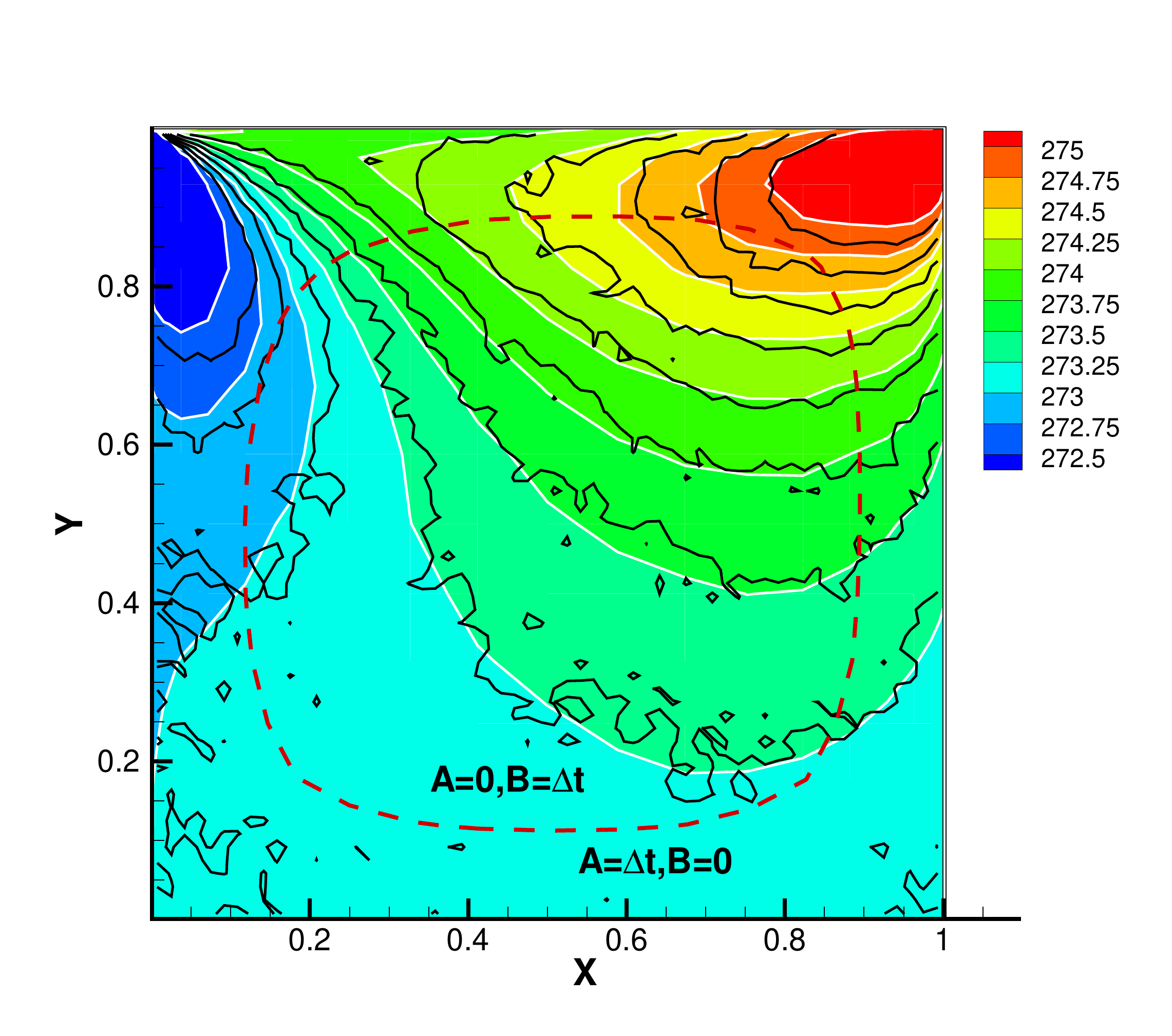}{a}
\includegraphics[width=6.0cm]{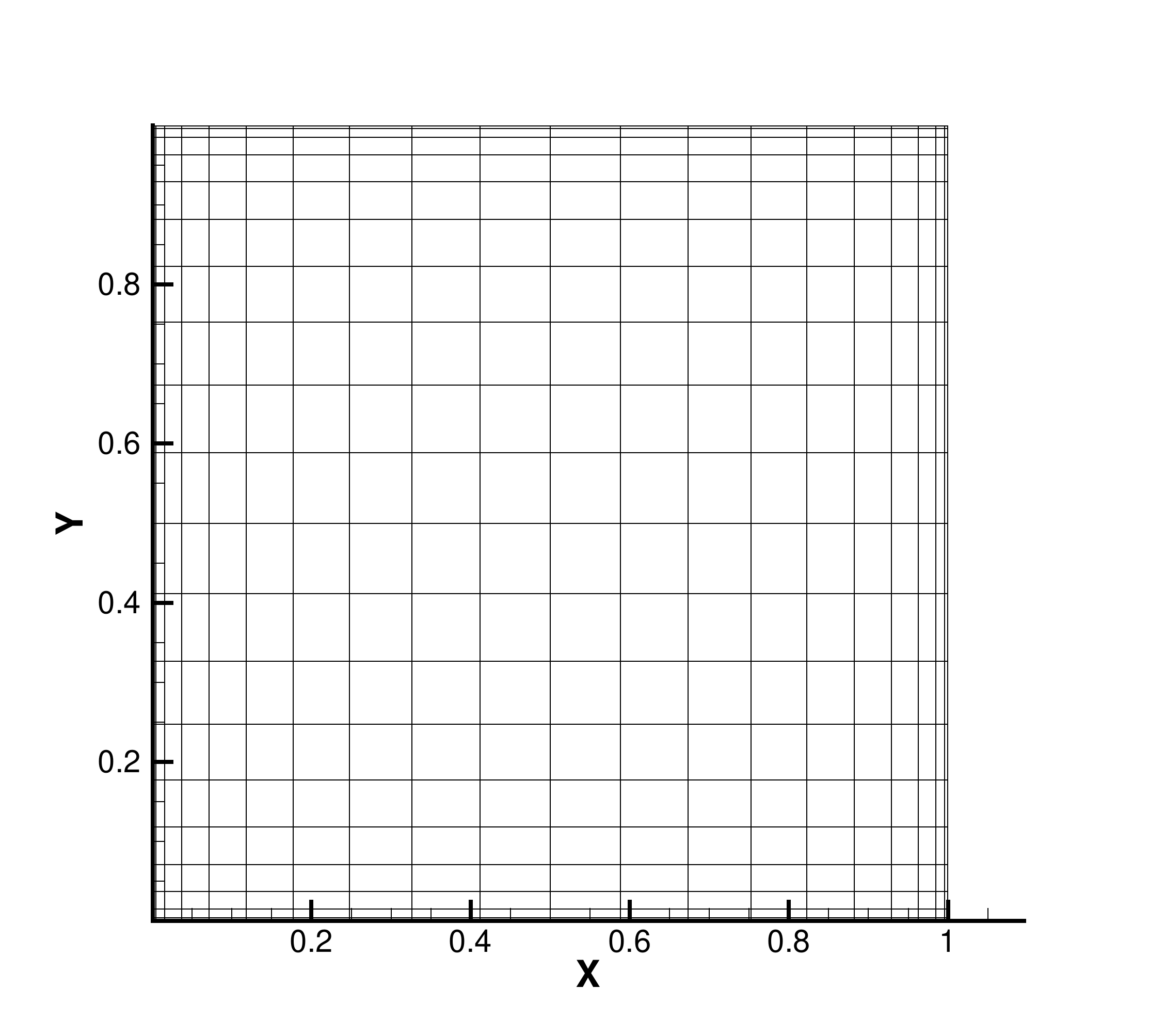}{b}
\includegraphics[width=6.0cm]{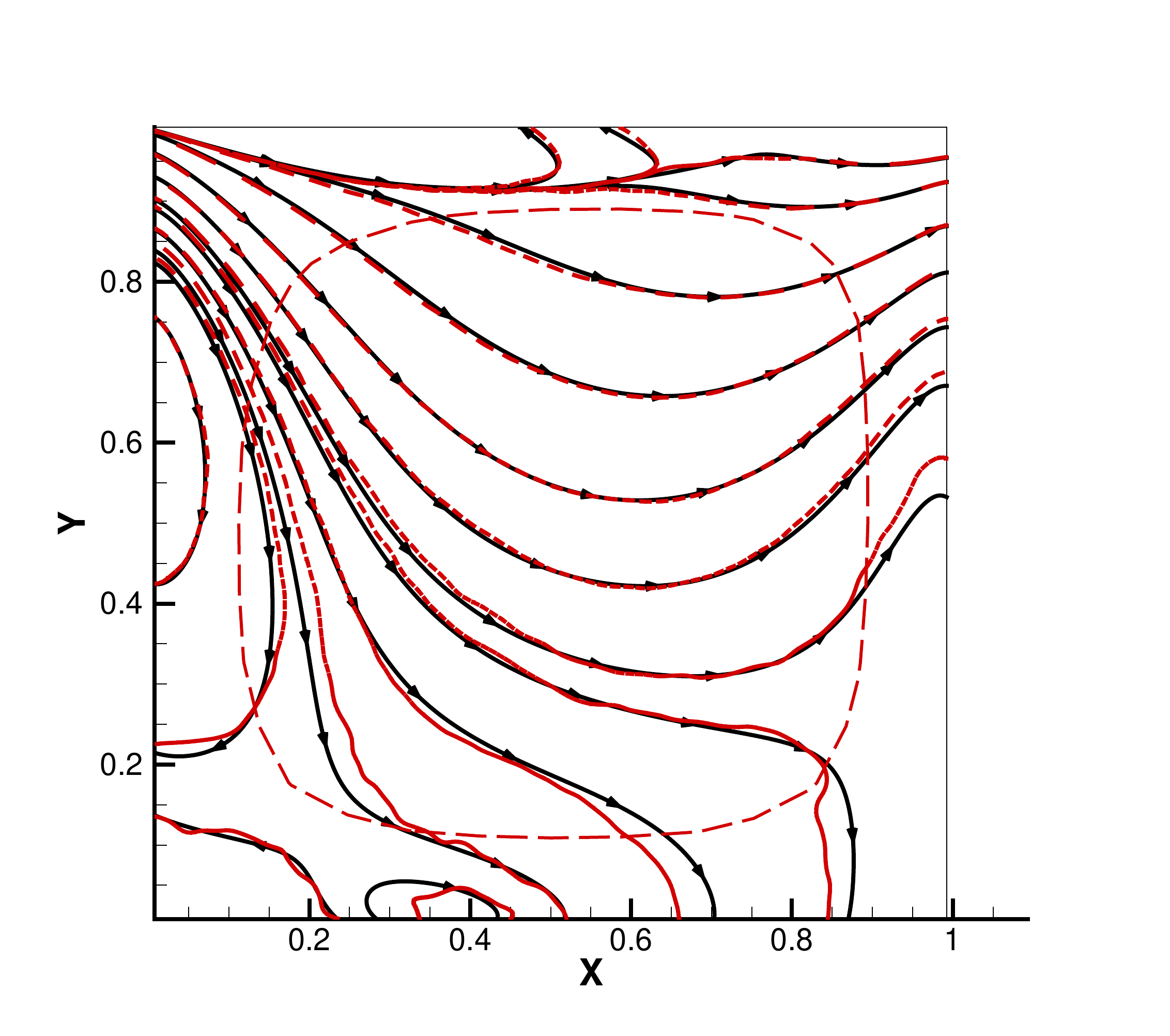}{c}
\includegraphics[width=6.0cm]{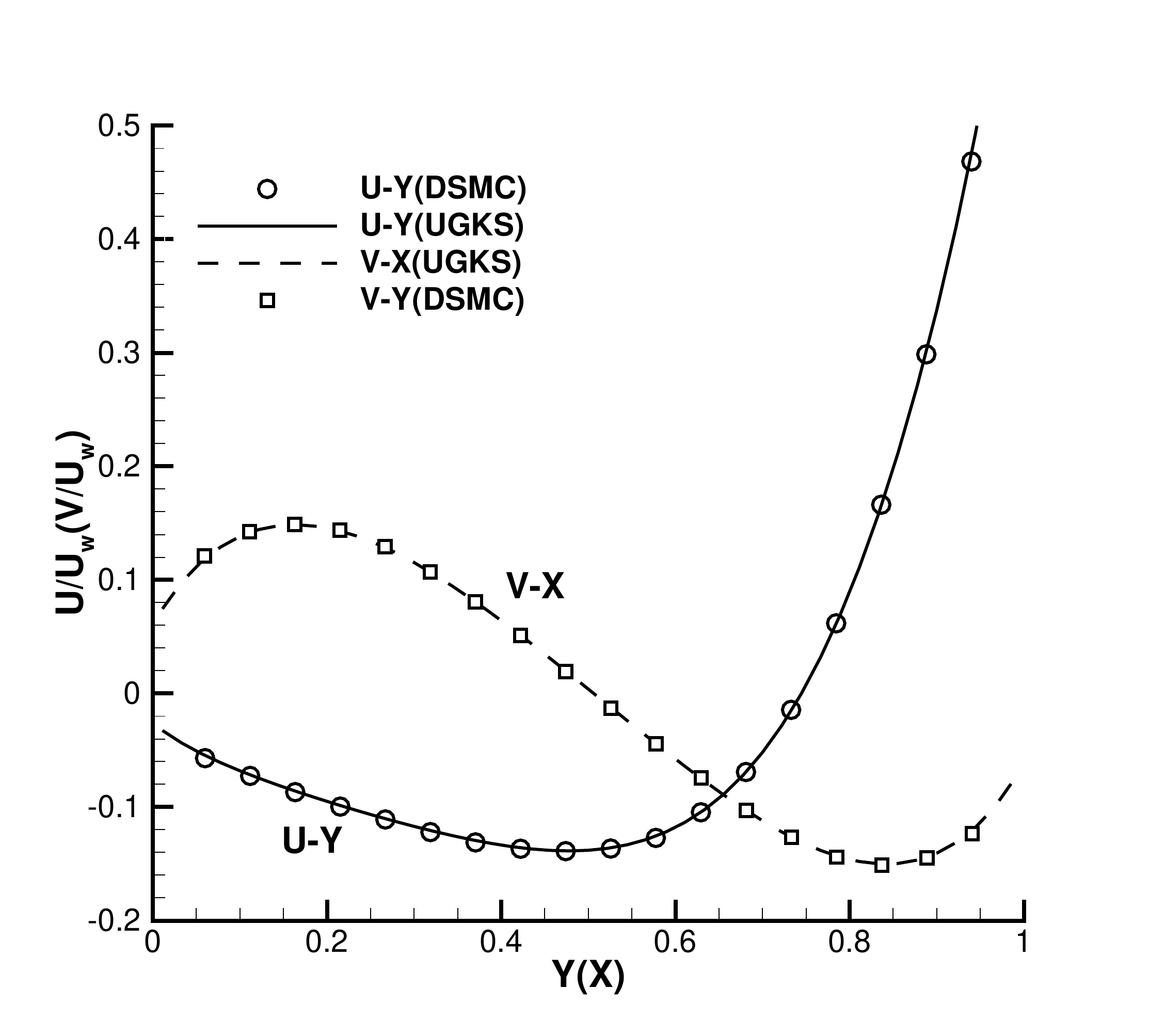}{d}
\caption{Cavity flow at \mbox{Kn}=0.075.
 (a) temperature contours with domain interface for different collision models, black lines: DSMC, white lines and background: UGKS;
 (b) Computational mesh in physical space;
 (c) heat flux, dash lines: DSMC, solid lines: UGKS;
 (d) U-velocity along the central vertical line and V-velocity along the central horizontal line, circles: DSMC, line:UGKS.}
 \label{k75}
\end{figure}

The first few tests are in the rarefied and transitional regime, where the UGKS solutions are compared with DSMC ones.
 Fig.(\ref{k10})-(\ref{k75}) show the results from UGKS and DSMC solutions of \cite{john2011effects} at Knudsen numbers $10$, $1$, and $0.075$.
The computational domain for $\mbox{Kn} =10$ and $\mbox{Kn} =1$ cases is composed of $50\times50$ non-uniform mesh in physical space
and $72\times72\times24$ points in  the velocity space.
Due to the reduction of Knudsen number, the mesh size over the particle mean free path can be changed significantly.
The computational domain for $\mbox{Kn} =0.075$ case is composed of $23\times23$ non-uniform mesh in physical space
and $32\times32\times12$ points in  the velocity space.
Due to the use of non-uniform of mesh and the local time step,
 Fig.(\ref{k75}) includes the switching interface between the use of the full Boltzmann collision term and the
Shakhov model. Even with the hybrid collision models, a smooth transition is obtained  in the solutions.
Same as the previous calculation \cite{Huangxuyu2}, perfect match with DSMC results has been obtained from the current UGKS.

\begin{figure}
\center
\includegraphics[width=6cm]{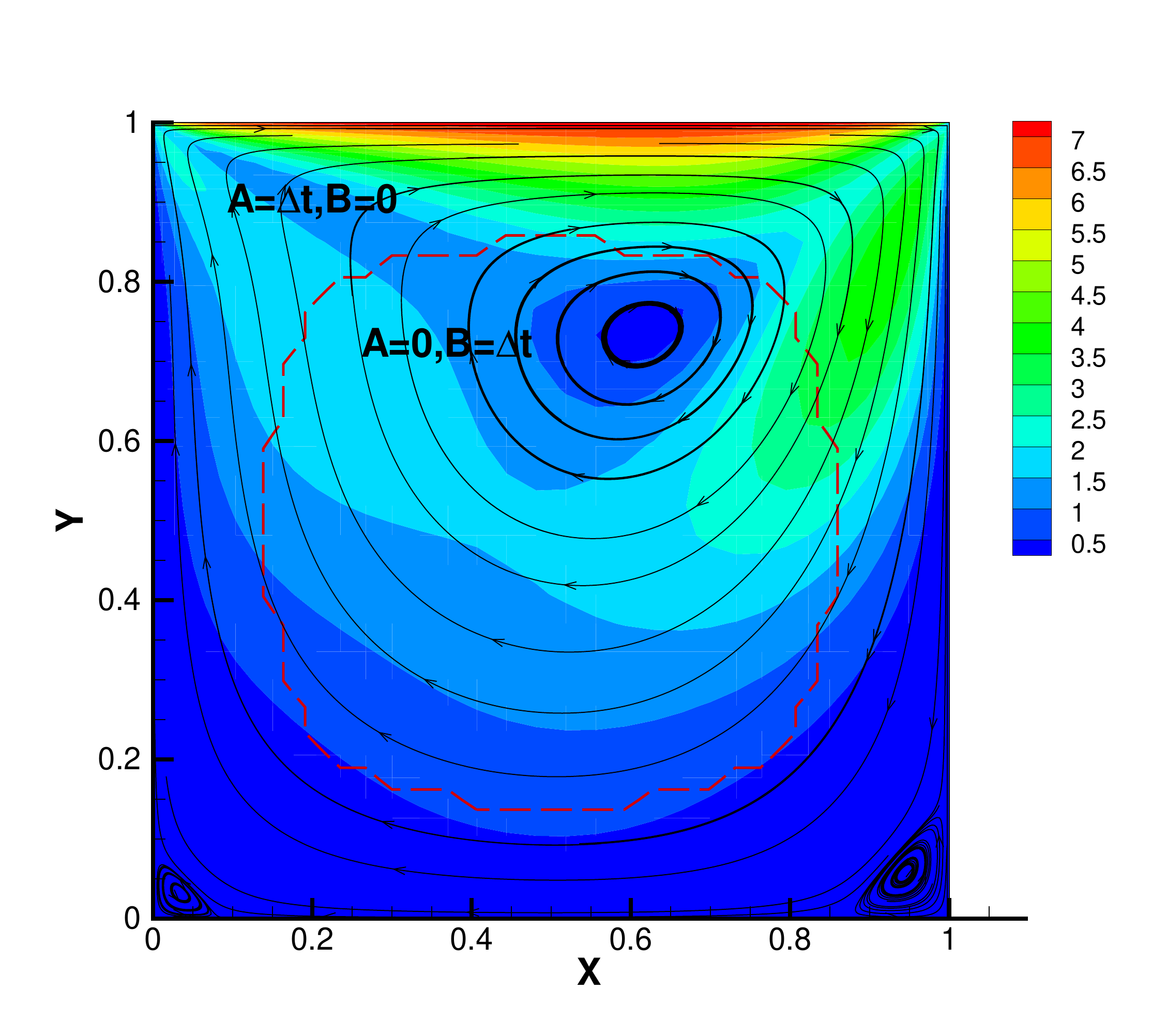}{a}
\includegraphics[width=6cm]{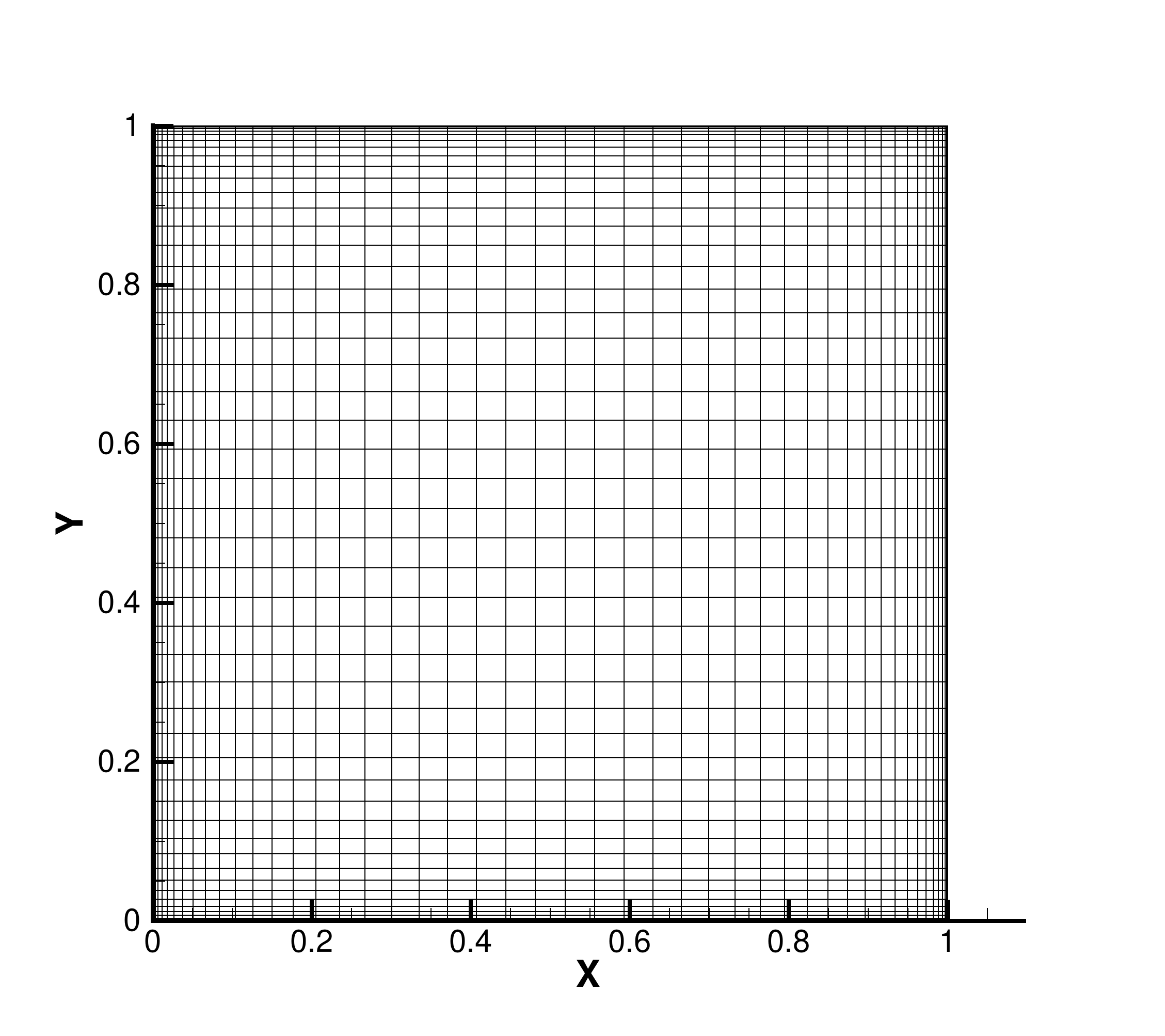}{b}
\includegraphics[width=6cm]{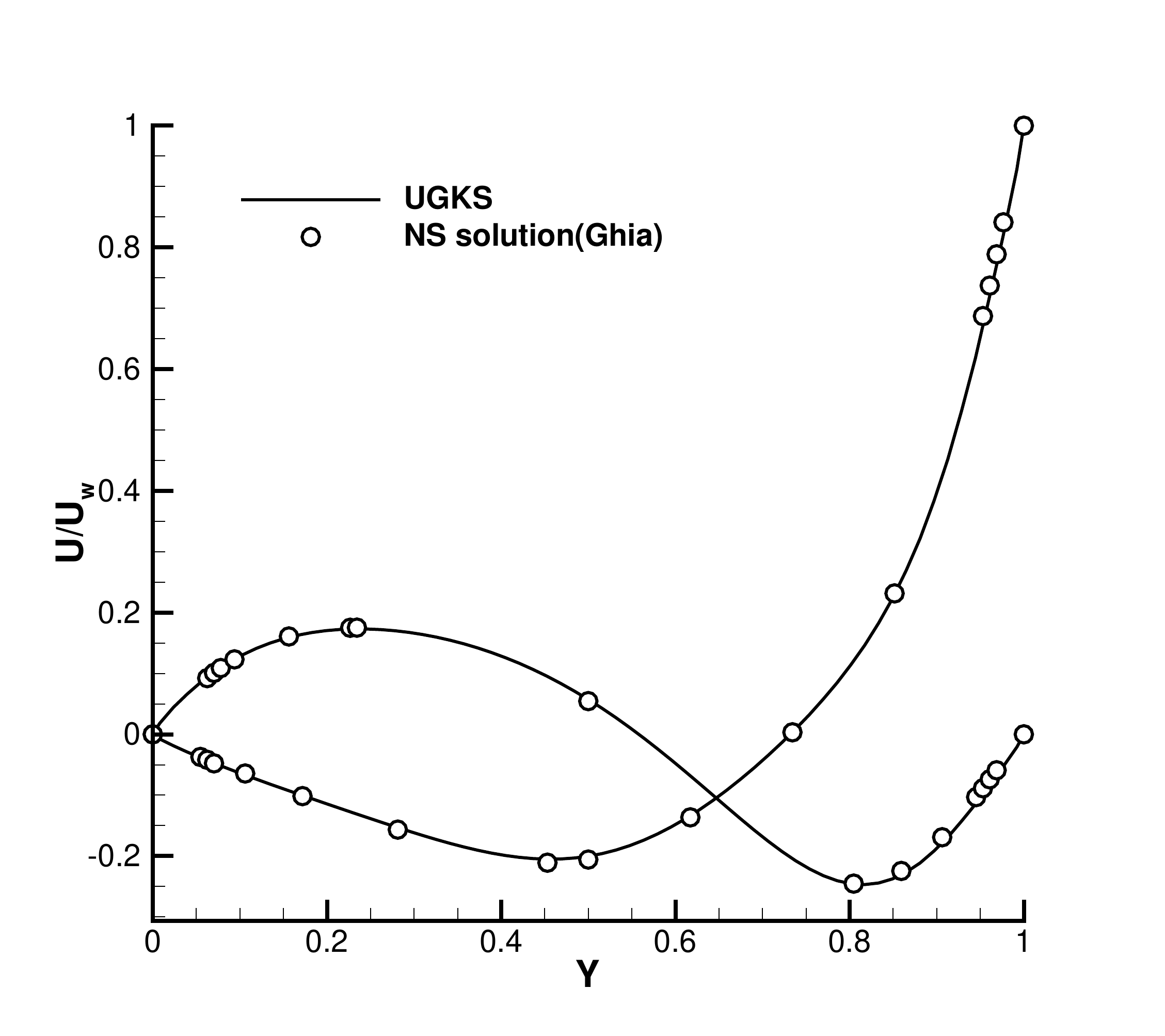}{c}
\caption{Cavity flow at $\mbox{Kn} = 1.42\times10^{-3}$ and $\mbox{Re} = 100$.
 (a) stream lines with velocity contour background and domain interface for different collision models;
 (b) Computational mesh in physical space;
 (c) U-velocity along the central vertical line and V-velocity along the central horizontal line, circles: NS solution, line: UGKS.}
 \label{r100}
\end{figure}

\begin{figure}
\begin{minipage}{0.1\linewidth}
  \includegraphics[width=8.0cm]{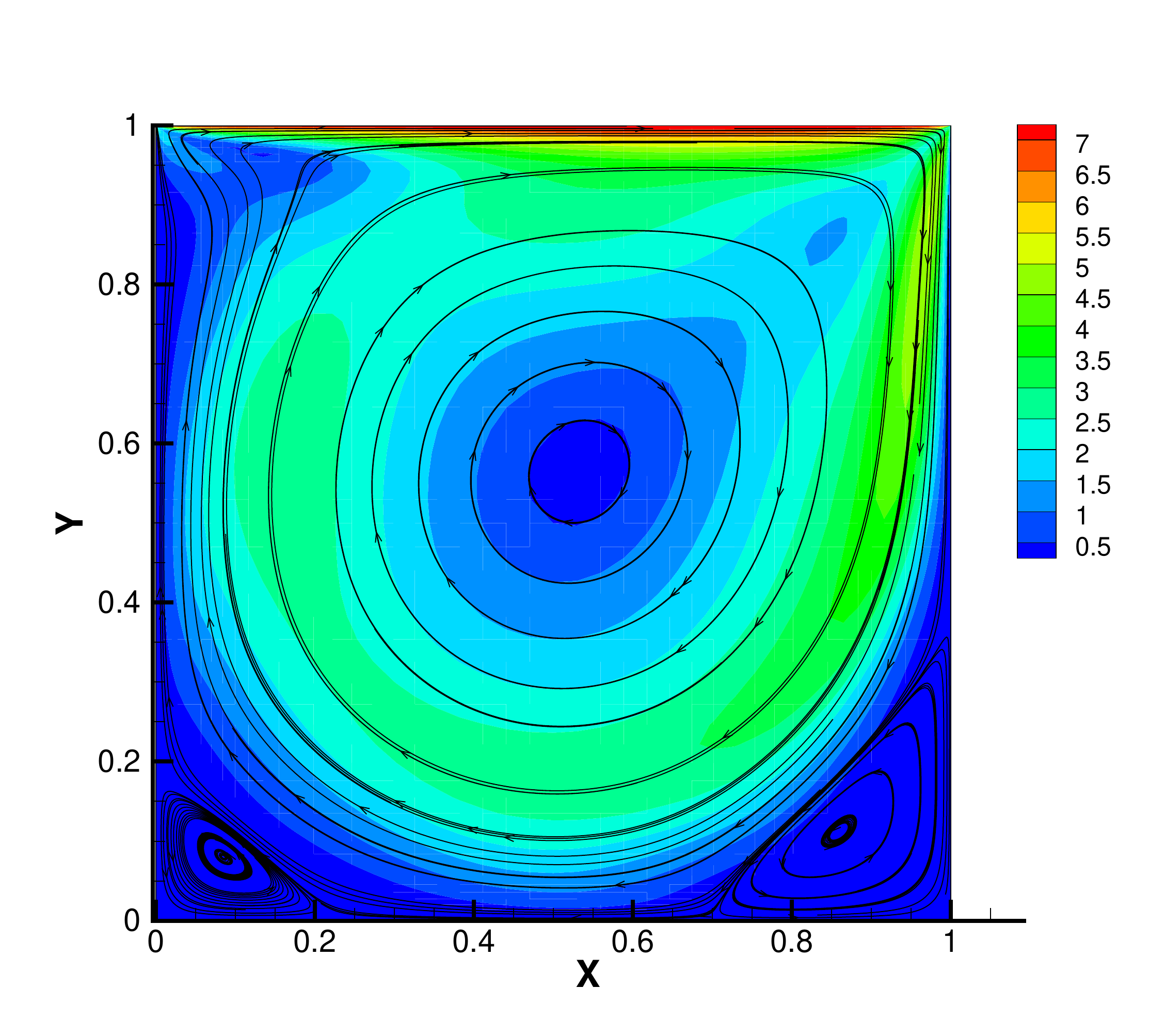}
\end{minipage}
\hfill
\begin{minipage}{0.48\linewidth}
\begin{minipage}{0.48\linewidth}
  \centerline{\includegraphics[width=3.5cm]{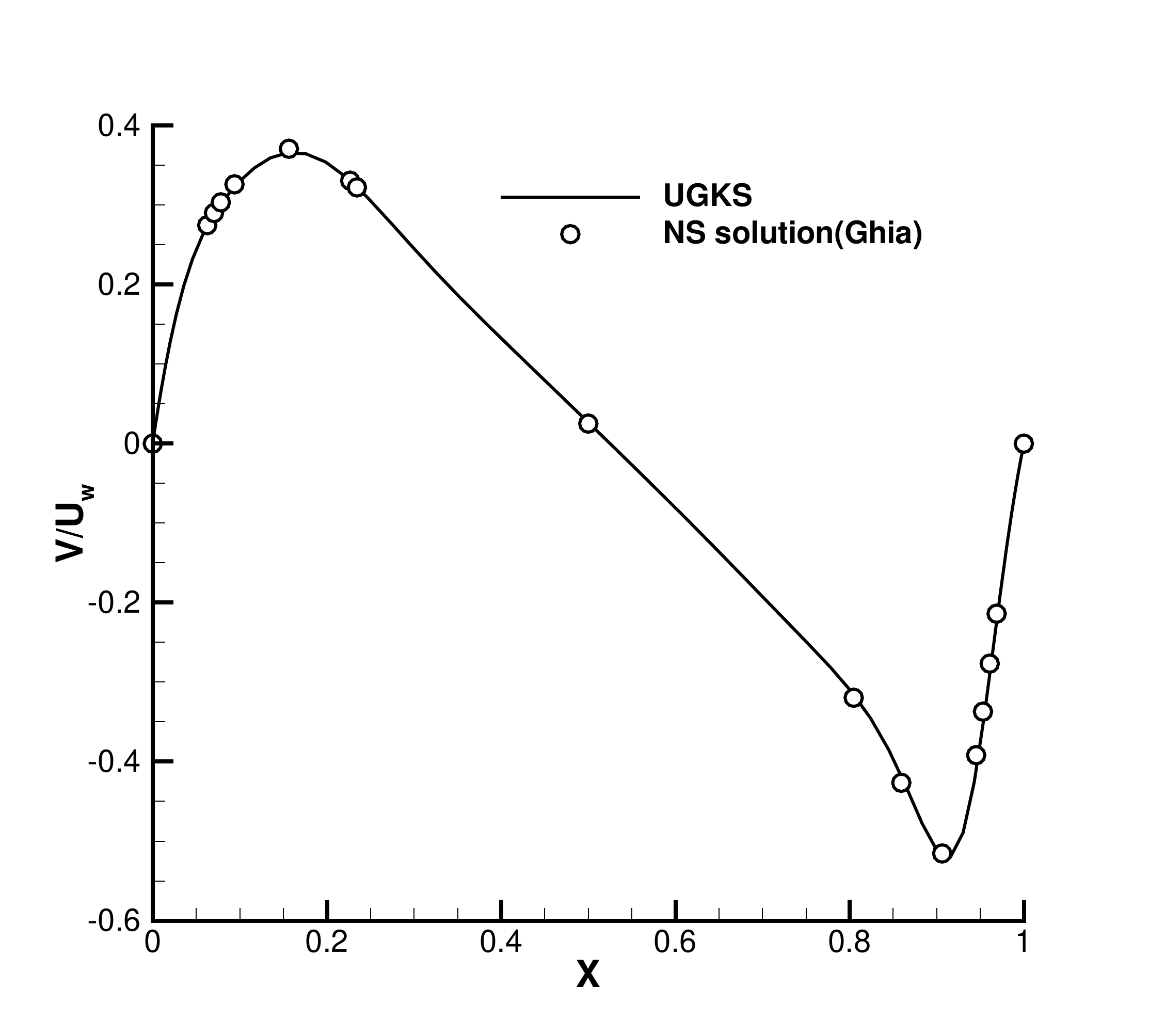}}
\end{minipage}
\hfill
\begin{minipage}{0.48\linewidth}
  \centerline{\includegraphics[width=3.5cm]{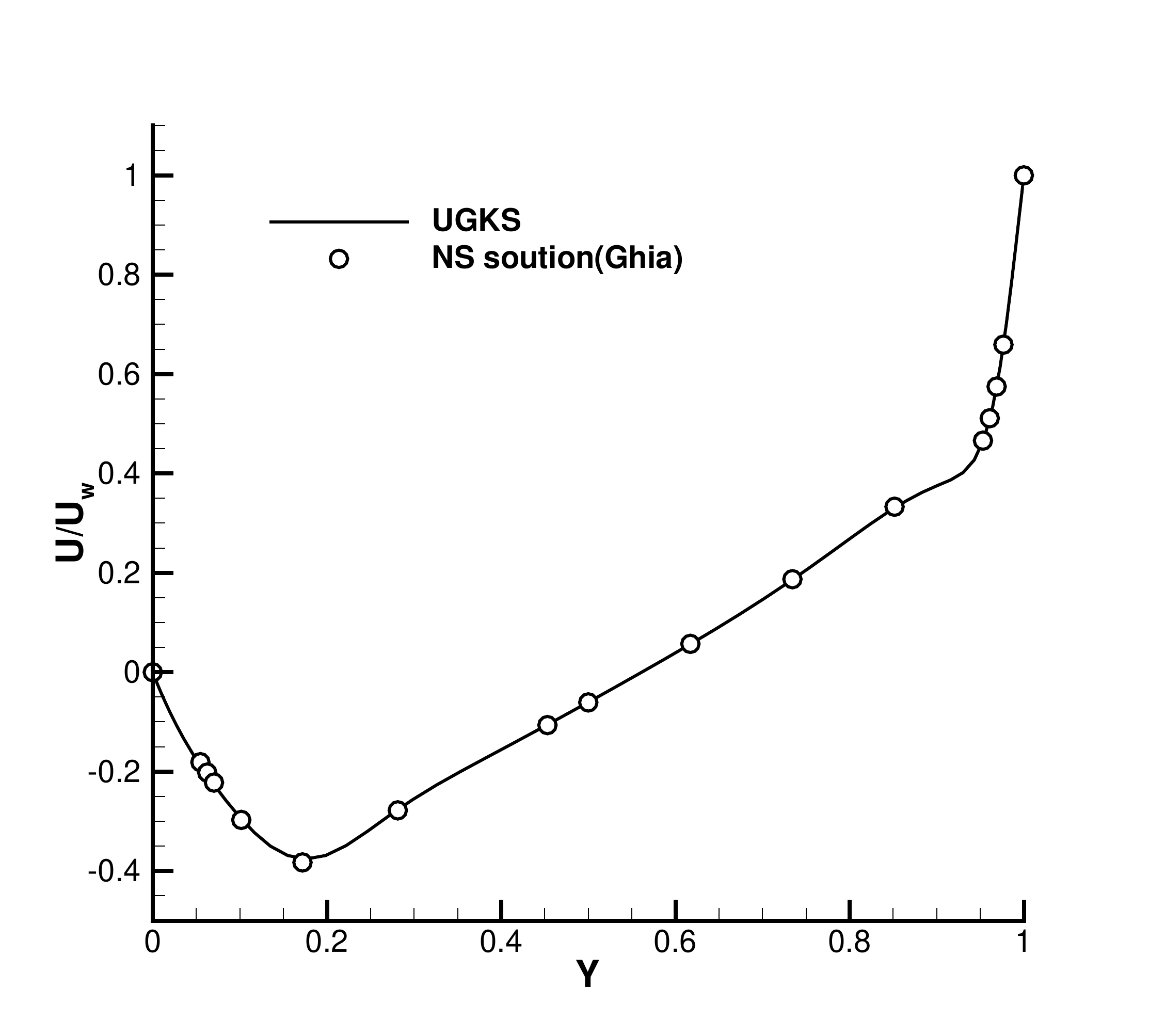}}
\end{minipage}
\vfill
\begin{minipage}{0.48\linewidth}
  \centerline{\includegraphics[width=3.5cm]{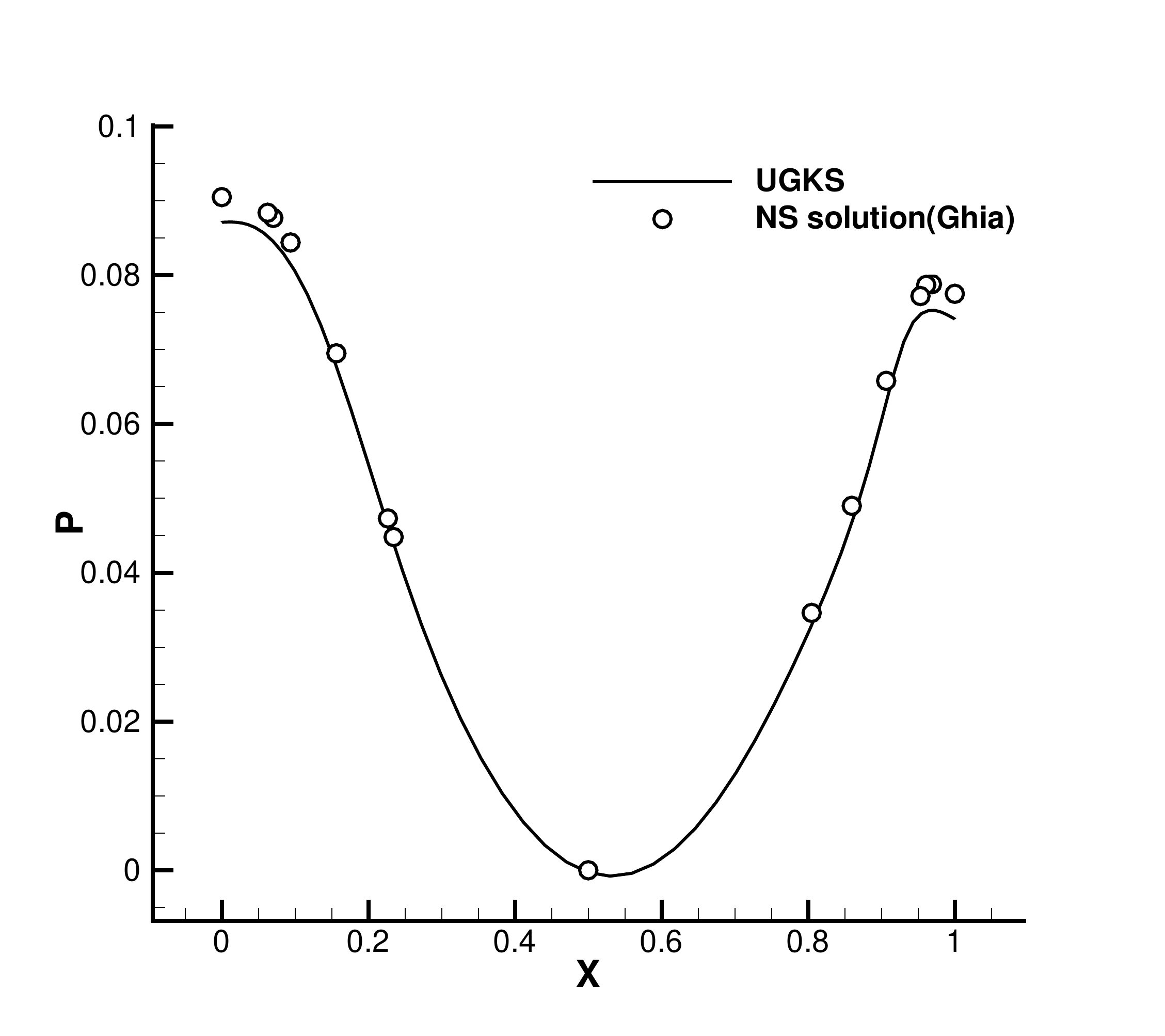}}
\end{minipage}
\hfill
\begin{minipage}{0.48\linewidth}
  \centerline{\includegraphics[width=3.5cm]{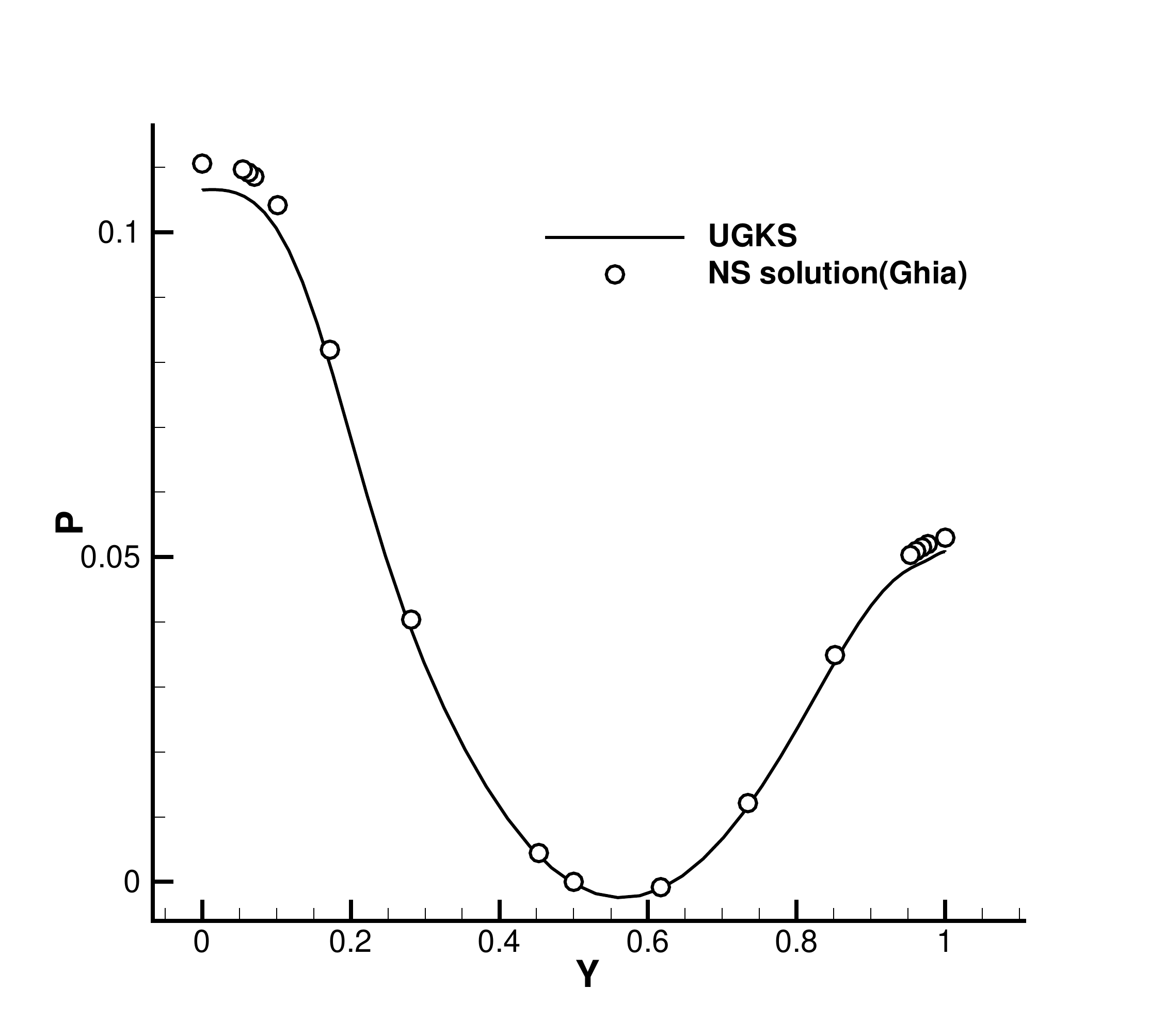}}
\end{minipage}
\end{minipage}\\
\caption{Cavity flow at $\mbox{Kn} = 1.42\times10^{-4}$ and $\mbox{Re} = 1000$.
 (left) velocity stream lines with velocity contour background;
 (right) U-velocity along the central vertical line,
  V-velocity along the central horizontal line,
 pressure along the central vertical line,
 and pressure along the central horizontal line, circles: NS solution, line: UGKS.}
 \label{r1000}
\end{figure}

The next two test cases are numerically to validate the AP property of the current scheme in the continuum flow regime at Knudsen numbers $1.42 \times 10^{-3}$ and $1.42 \times 10^{-4}$ or
$\mbox{Re} =100$ and $1000$.
The computational domain for $\mbox{Re} =100$ and $\mbox{Re} =1000$ is composed of $61\times61$ non-uniform mesh in physical space
and $32\times32$ points in  the velocity space. In both cases, the freedom of molecule is restricted in a 2-D  space in order to get the flow condition
close to the 2-D incompressible flow limit.
Also, the non-slip boundary condition is imposed in these two calculations.
Fig.(\ref{r100}) and (\ref{r1000}) show the UGKS results and reference NS solutions \cite{ghia1982high}.
This clearly demonstrates that the UGKS converges to the NS solutions accurately in the hydrodynamic limit.
 Fig.(\ref{r100}) also shows the switching interface between the full Boltzmann and Shakhov model, a smooth solution is obtained across the interface. At $\mbox{Re}=1000$, the UGKS uses Shakhov model in the whole domain.

\begin{figure}
\center
\includegraphics[width=6cm]{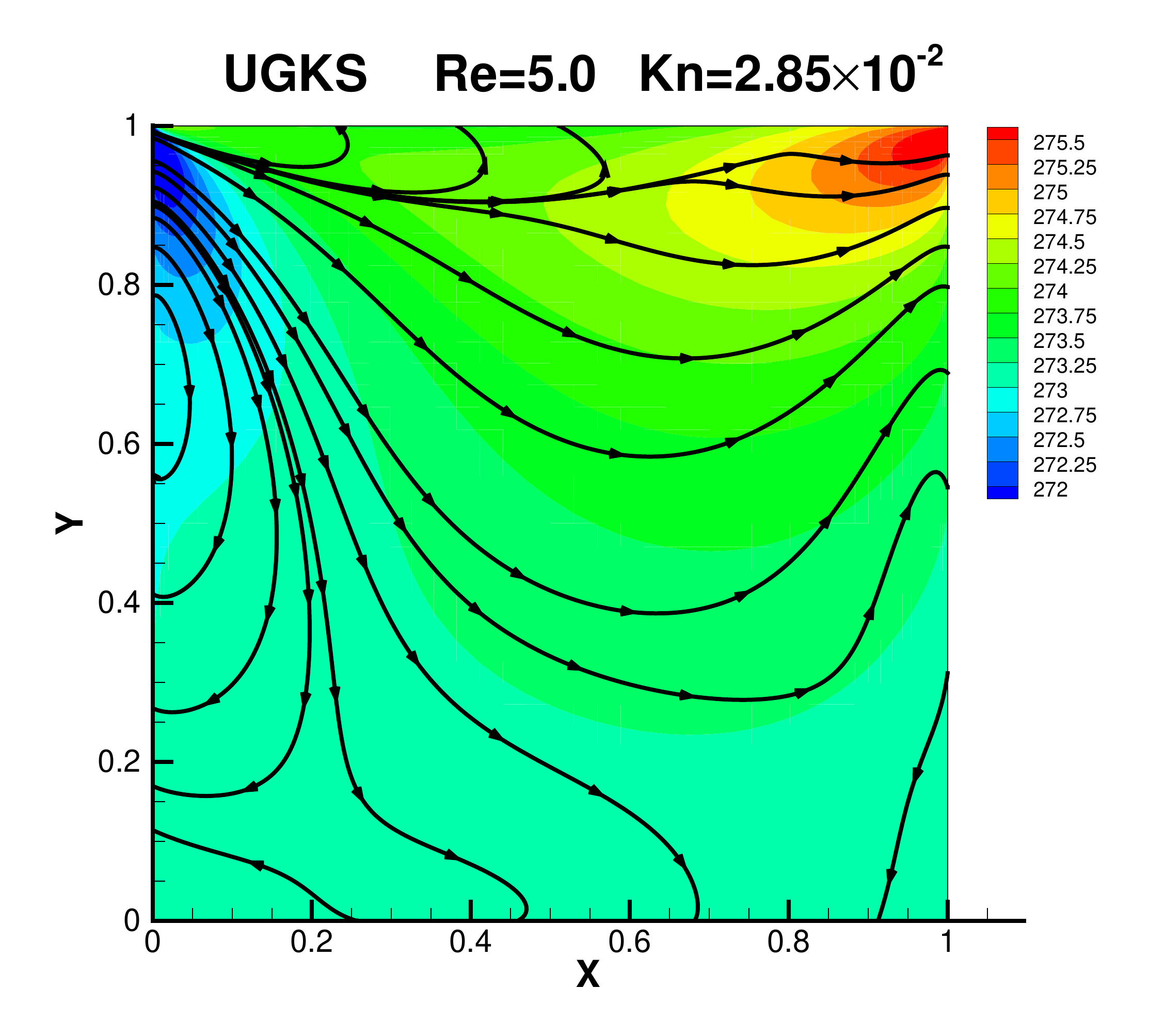}{a}
\includegraphics[width=6cm]{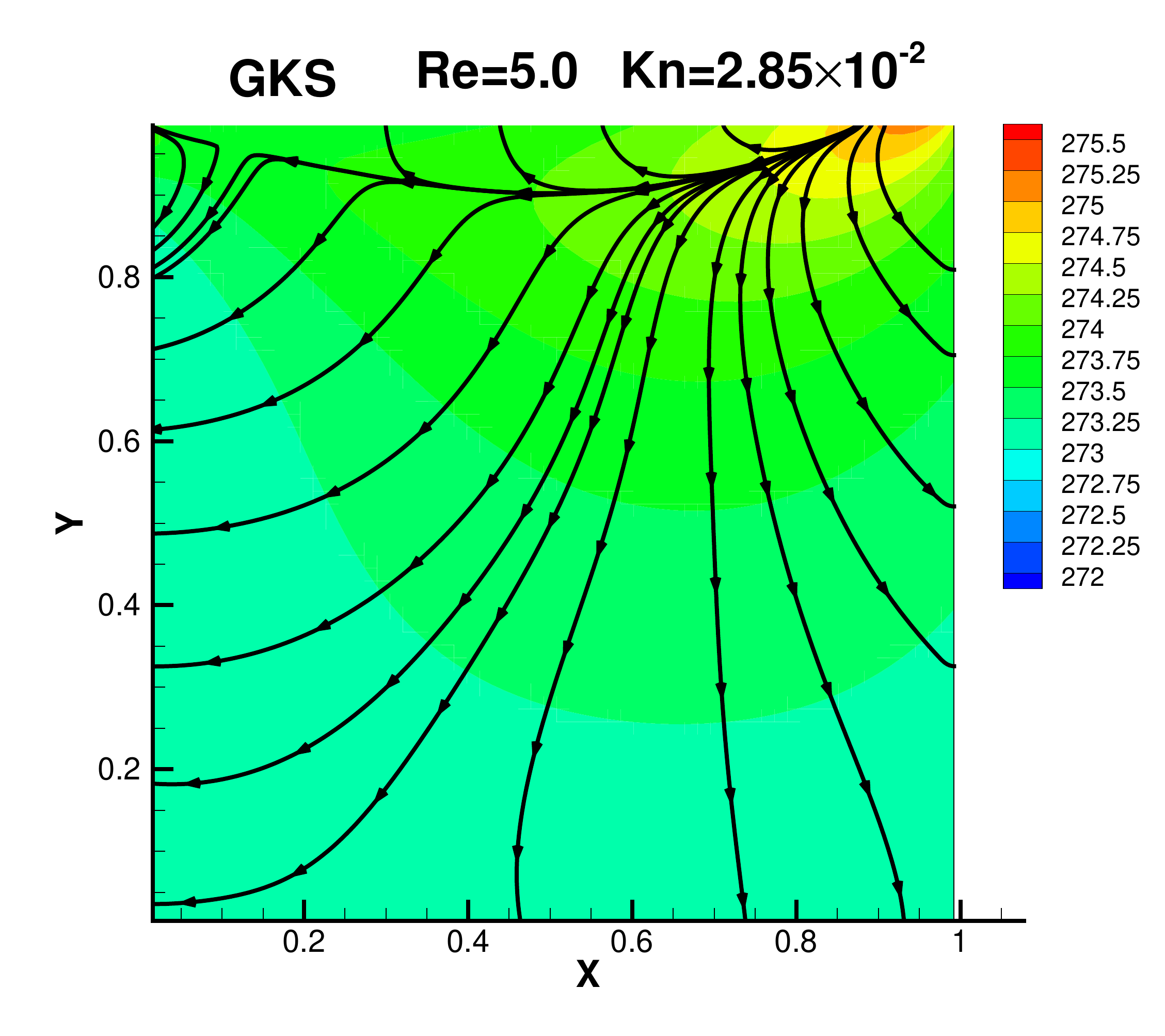}{b}
\includegraphics[width=6cm]{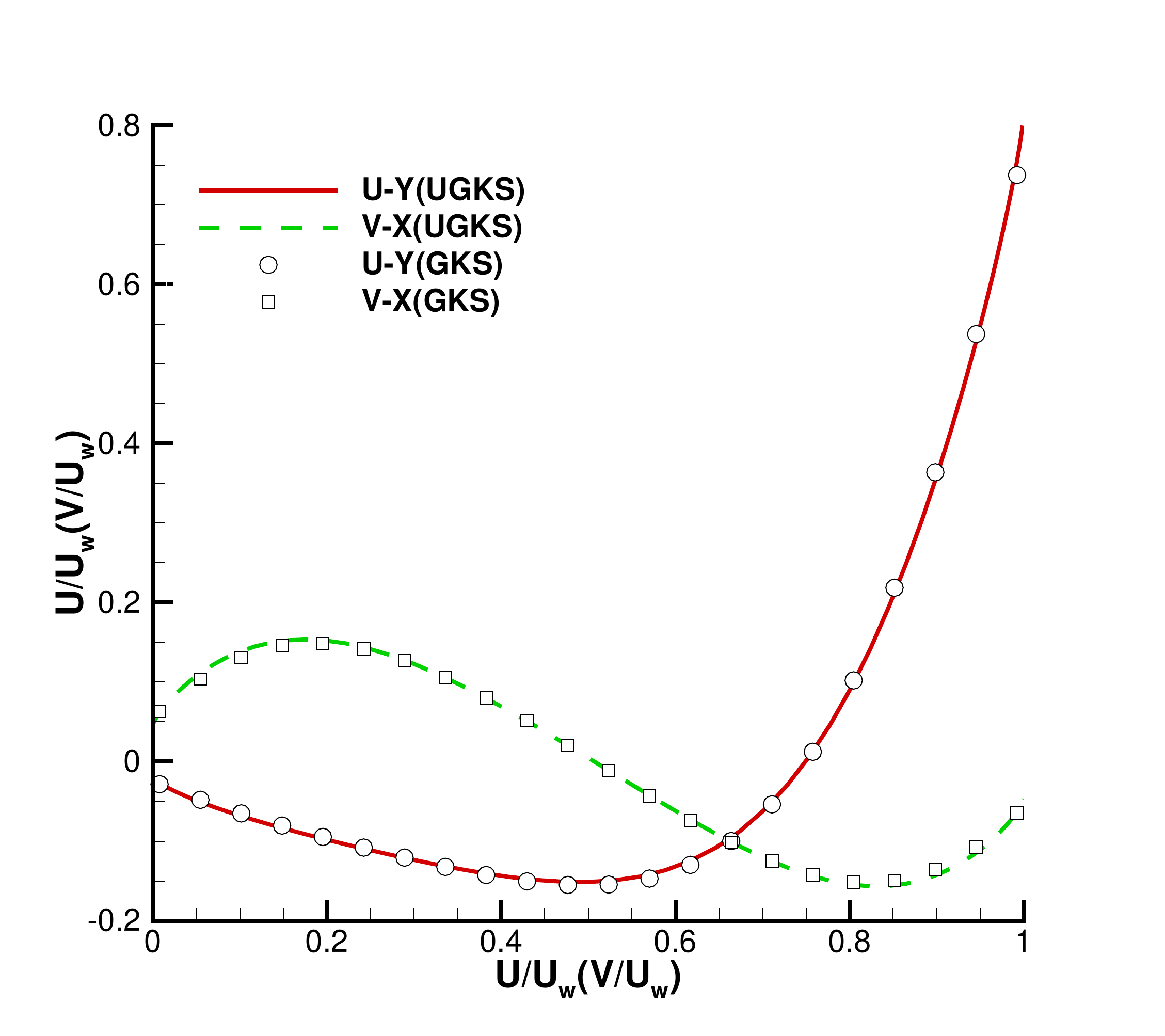}{c}
\caption{Cavity simulation using UGKS and GKS at $\mbox{Kn}=2.85\times10^{-2}$ and $\mbox{Re}=5$.
 (a) temperature contour and heat flux: UGKS;
 (b) temperature contour and heat flux: GKS;
 (c) U-velocity along the central vertical line and V-velocity along the central horizontal line, circles: GKS, line: UGKS.}
 \label{r5}
\end{figure}

\begin{figure}
\center
\includegraphics[width=6cm]{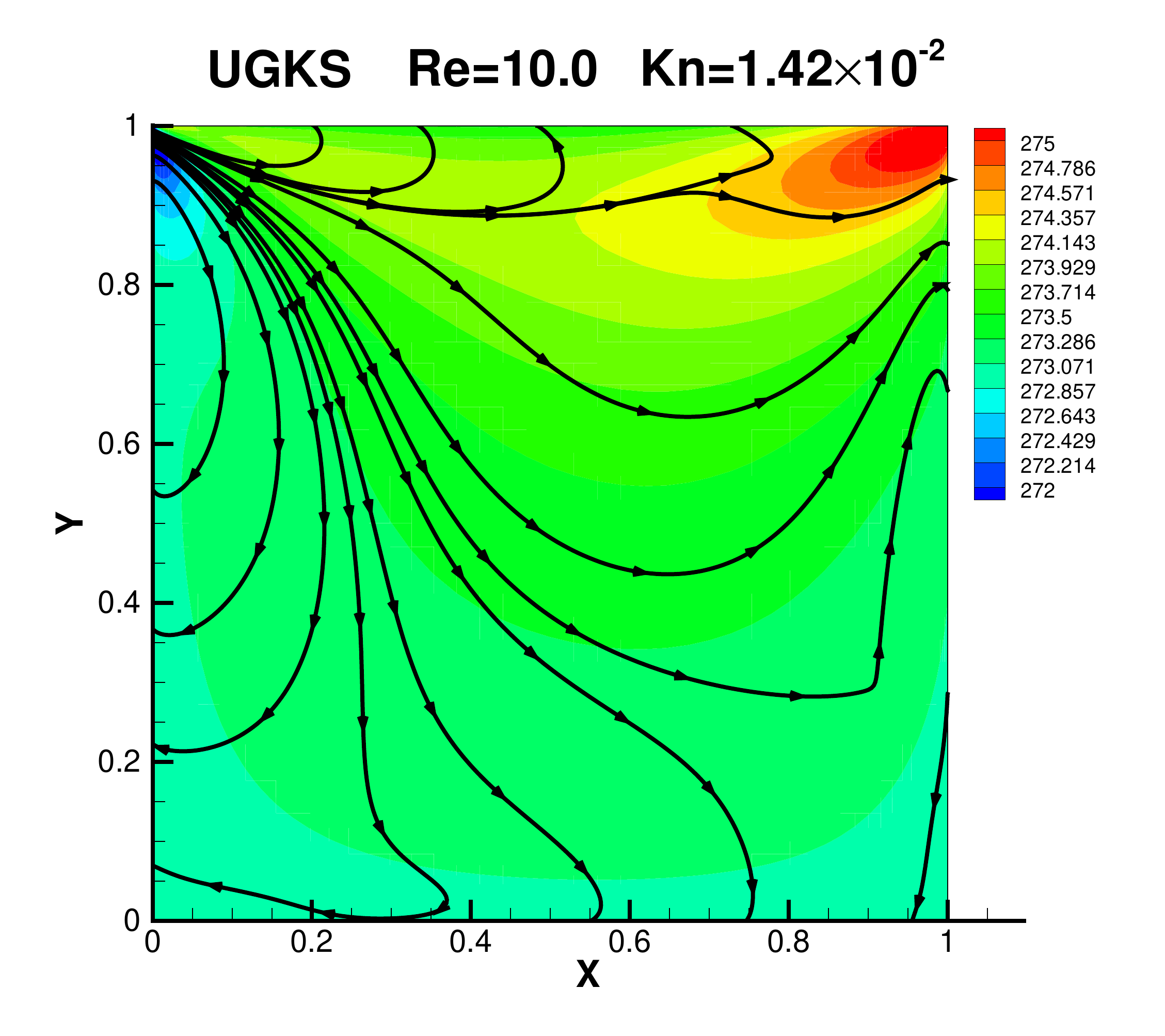}{a}
\includegraphics[width=6cm]{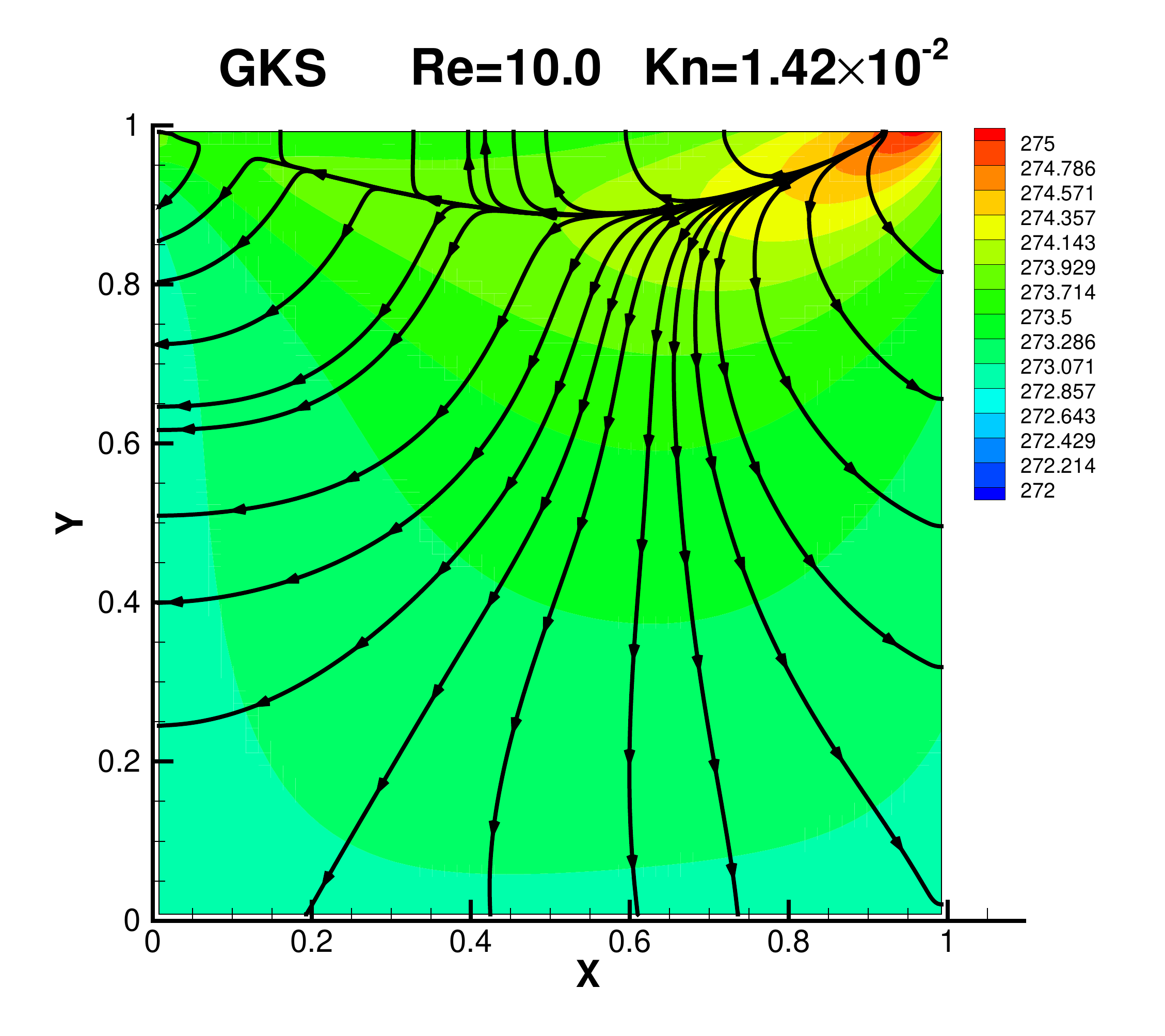}{b}
\includegraphics[width=6cm]{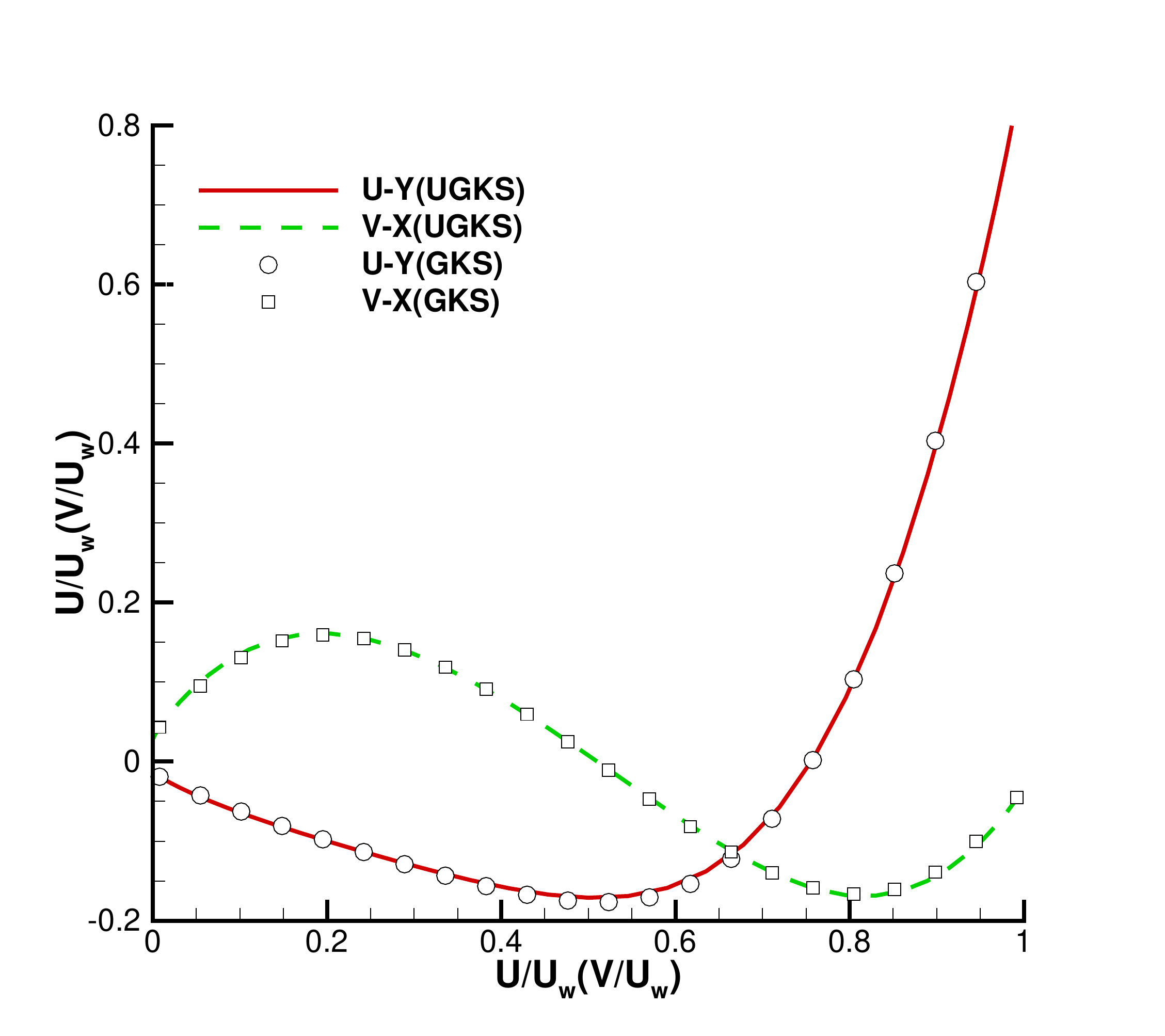}{c}
\caption{Cavity simulation using UGKS and GKS at $\mbox{Kn}=1.42\times10^{-2}$ and $\mbox{Re}=10$.
(a) temperature contour and heat flux: UGKS;
 (b) temperature contour and heat flux: GKS;
 (c) U-velocity along the central vertical line and V-velocity along the central horizontal line, circles: GKS, line: UGKS.}
 \label{r10}
\end{figure}

\begin{figure}
\center
\includegraphics[width=6cm]{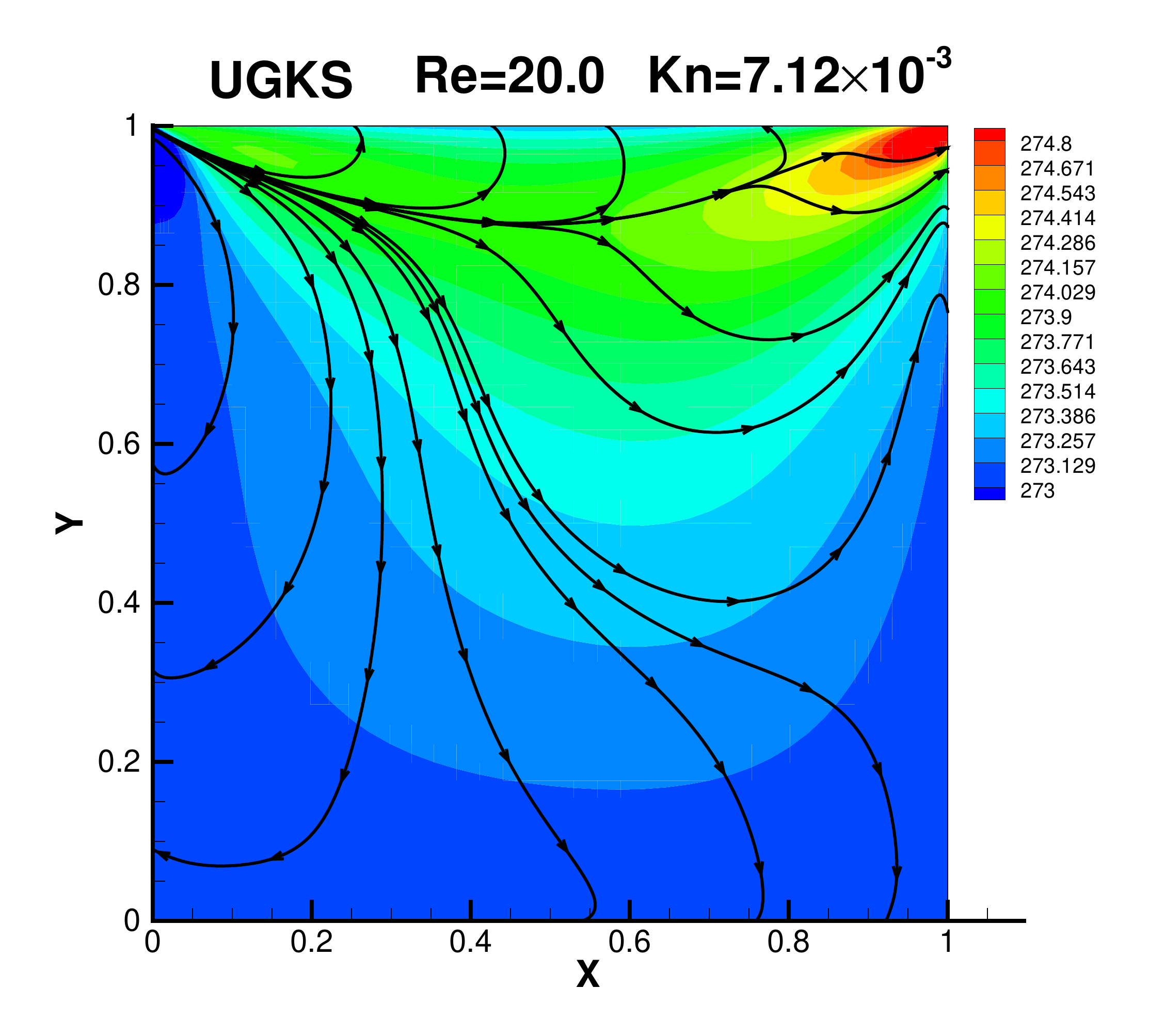}{a}
\includegraphics[width=6cm]{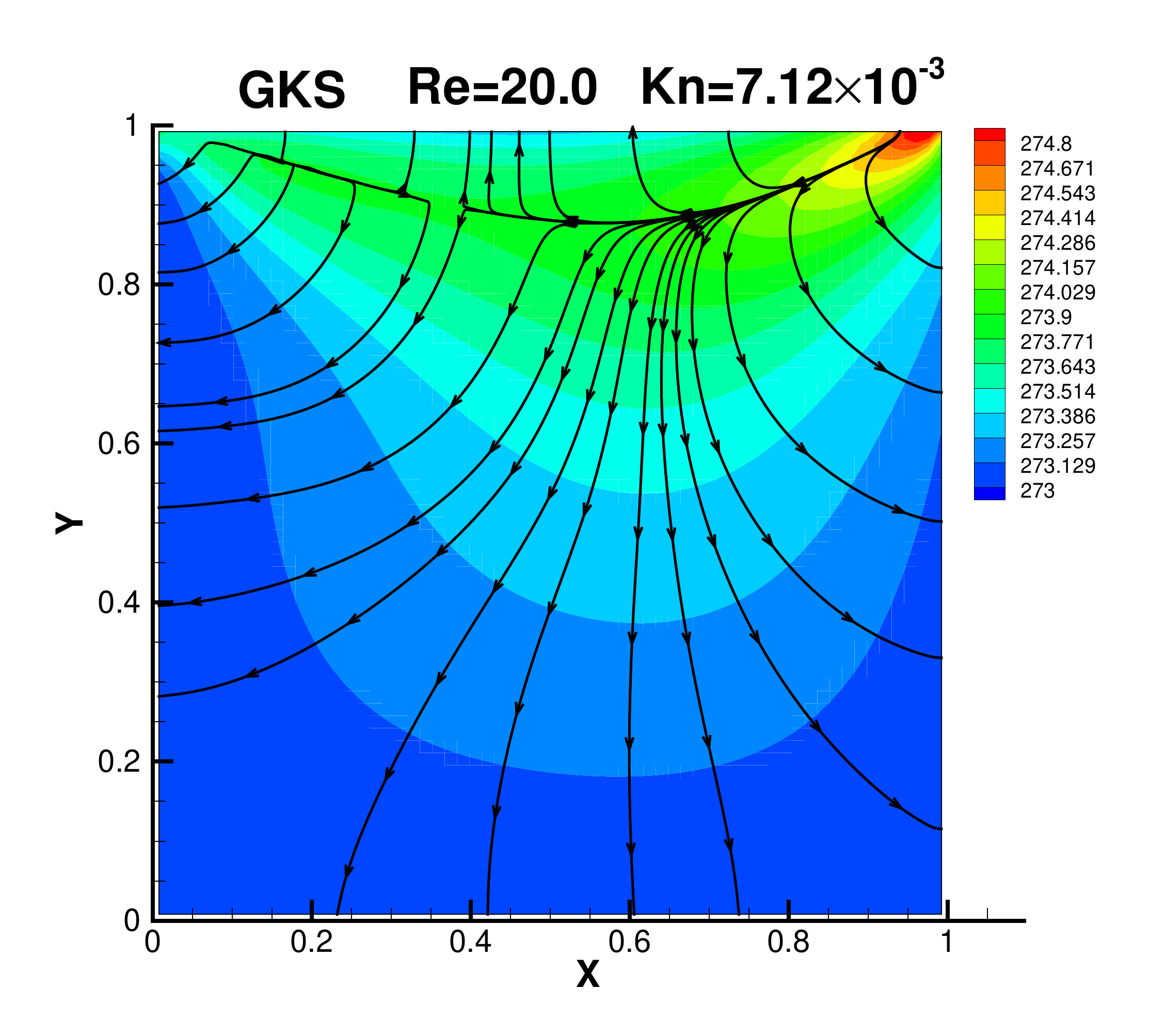}{b}
\includegraphics[width=6cm]{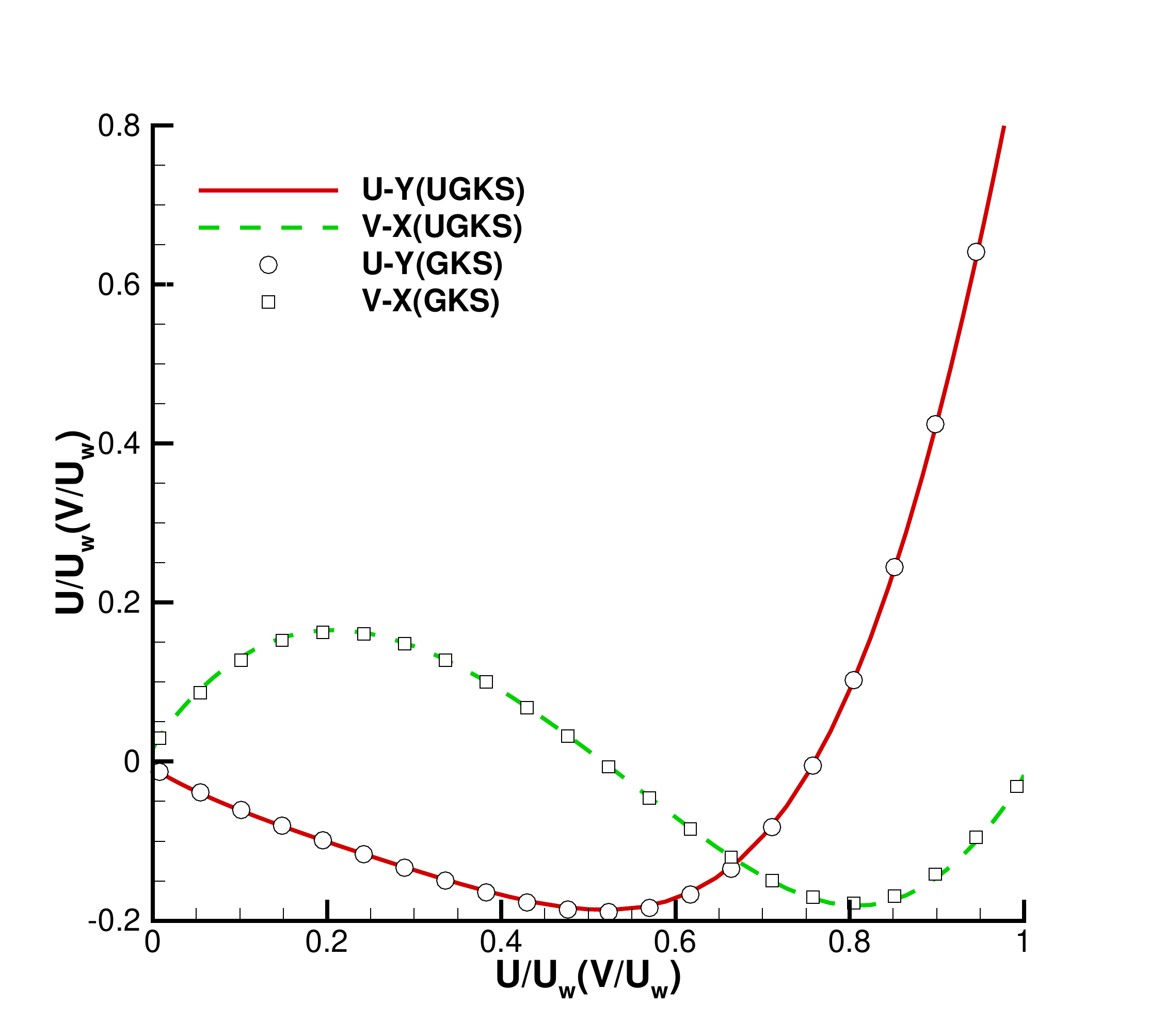}{c}
\caption{Cavity simulation using UGKS and GKS at $\mbox{Kn}=7.12\times10^{-3}$ and $\mbox{Re}=20$.
(a) temperature contour and heat flux: UGKS;
 (b) temperature contour and heat flux: GKS;
 (c) U-velocity along the central vertical line and V-velocity along the central horizontal line, circles: GKS, line: UGKS.}
 \label{r20}
\end{figure}

\begin{figure}
\center
\includegraphics[width=6cm]{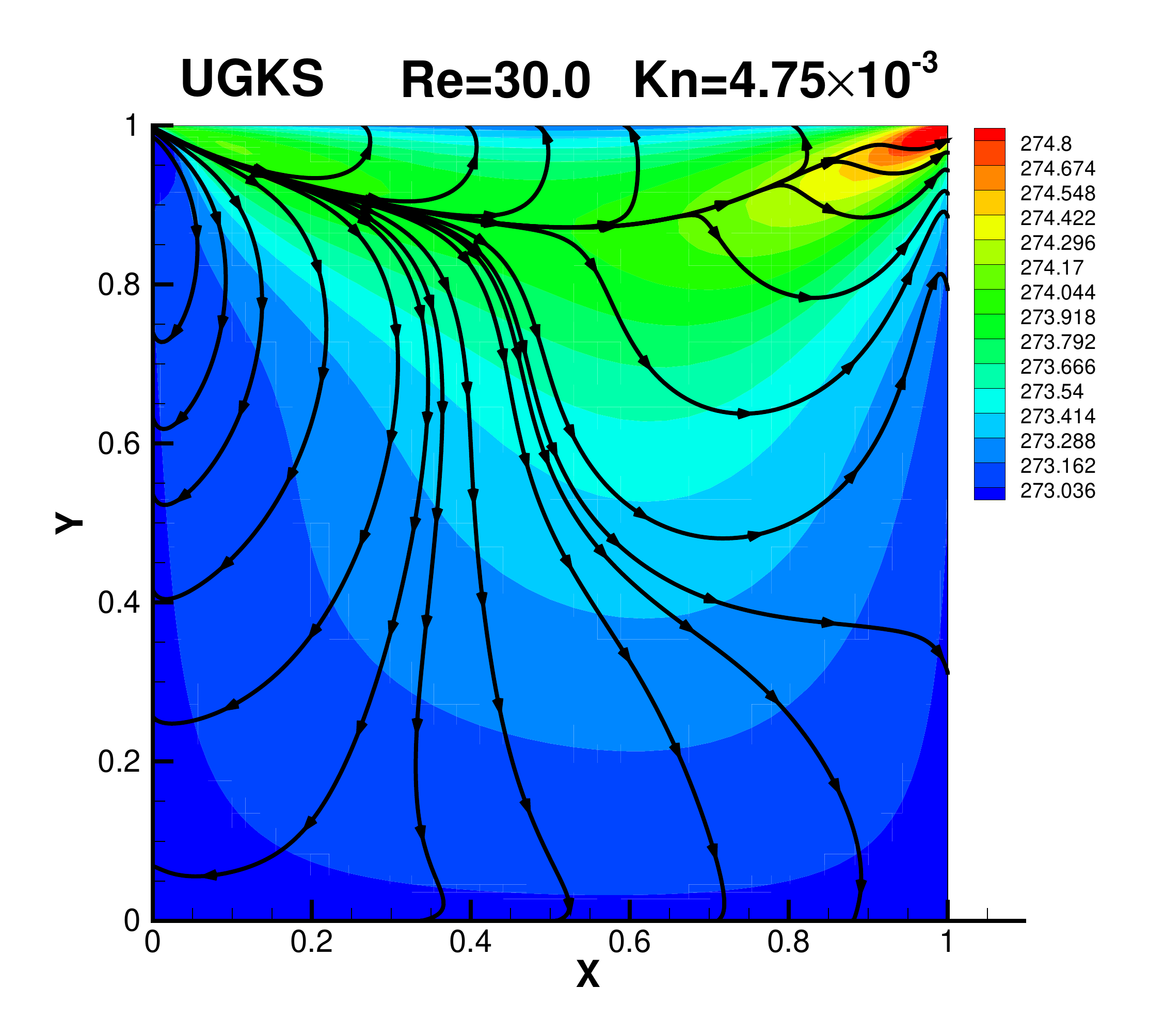}{a}
\includegraphics[width=6cm]{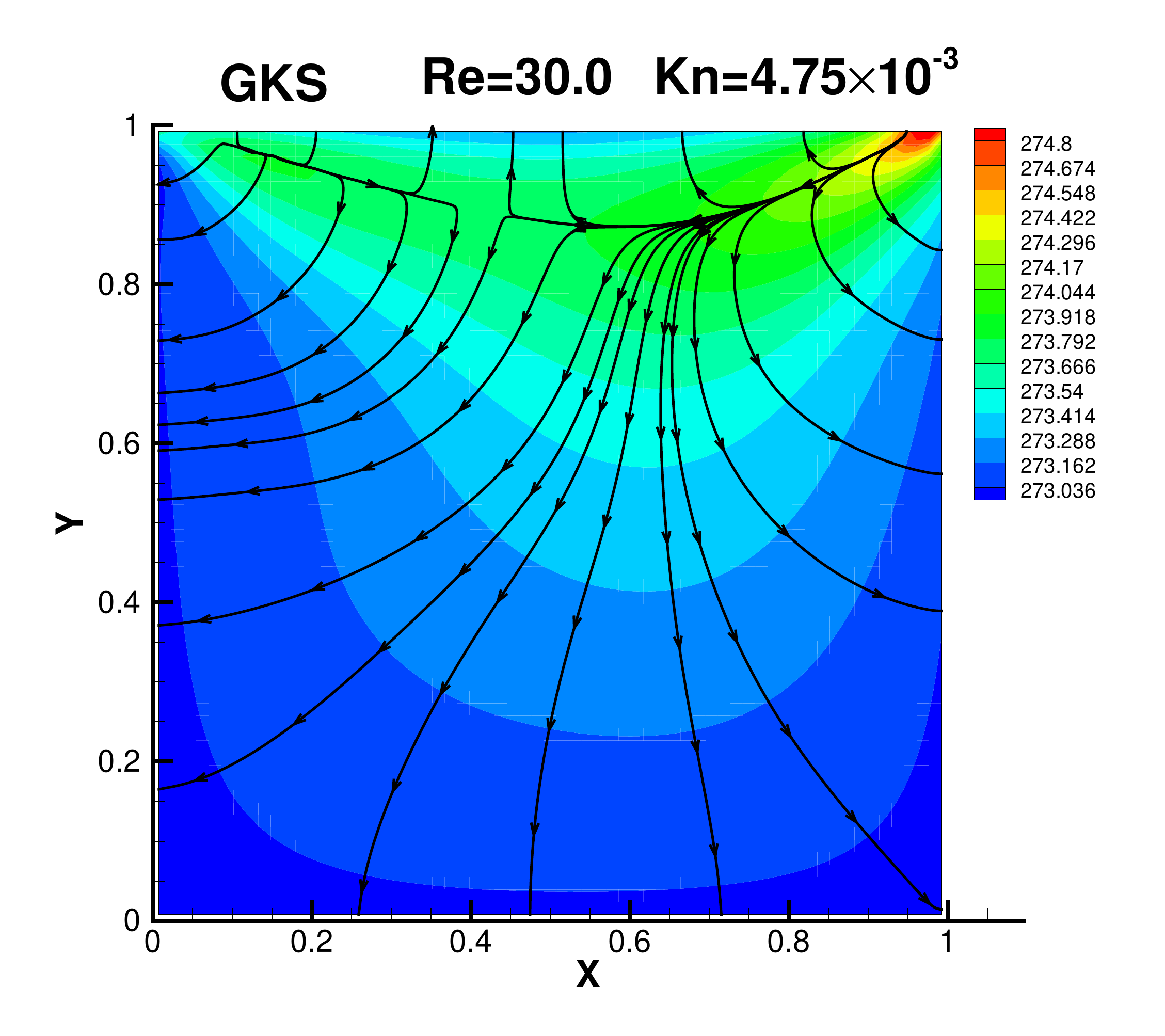}{b}
\includegraphics[width=6cm]{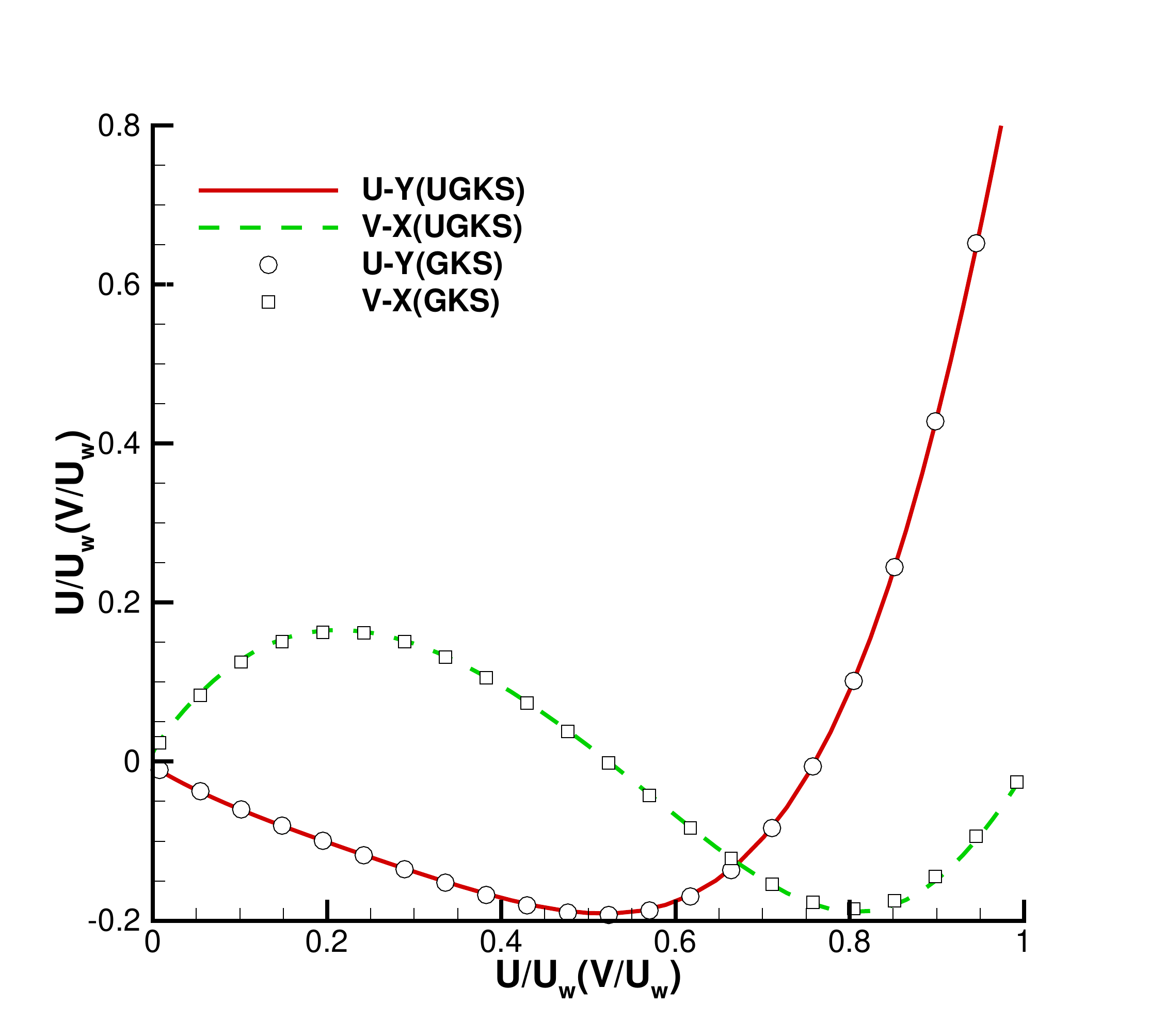}{c}
\caption{Cavity simulation using UGKS and GKS at $\mbox{Kn}=4.75\times10^{-3}$ and $\mbox{Re}=30$.
(a) temperature contour and heat flux: UGKS;
 (b) temperature contour and heat flux: GKS;
 (c) U-velocity along the central vertical line and V-velocity along the central horizontal line, circles: GKS, line: UGKS.}
 \label{r30}
\end{figure}

\begin{figure}
\center
\includegraphics[width=6cm]{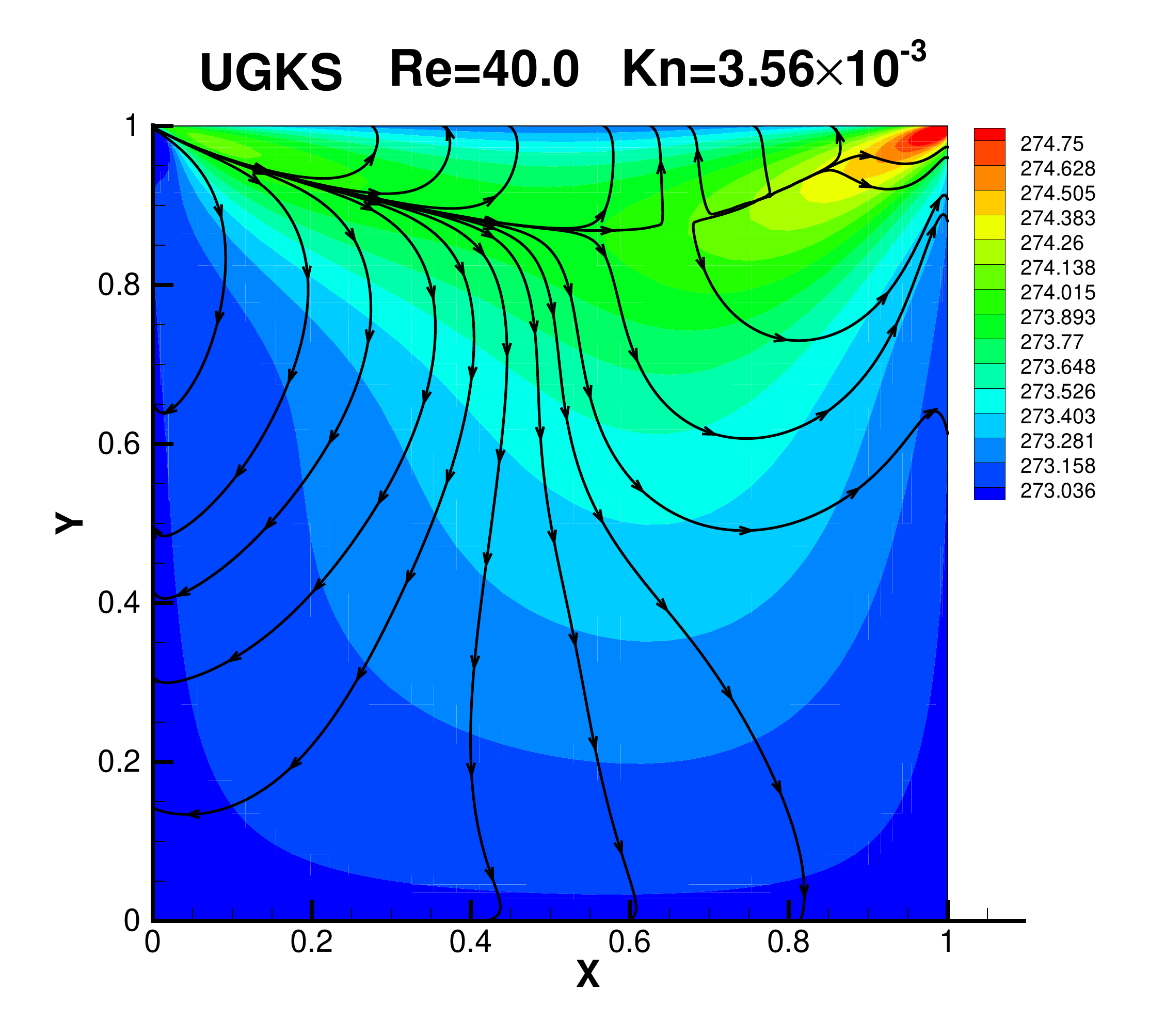}{a}
\includegraphics[width=6cm]{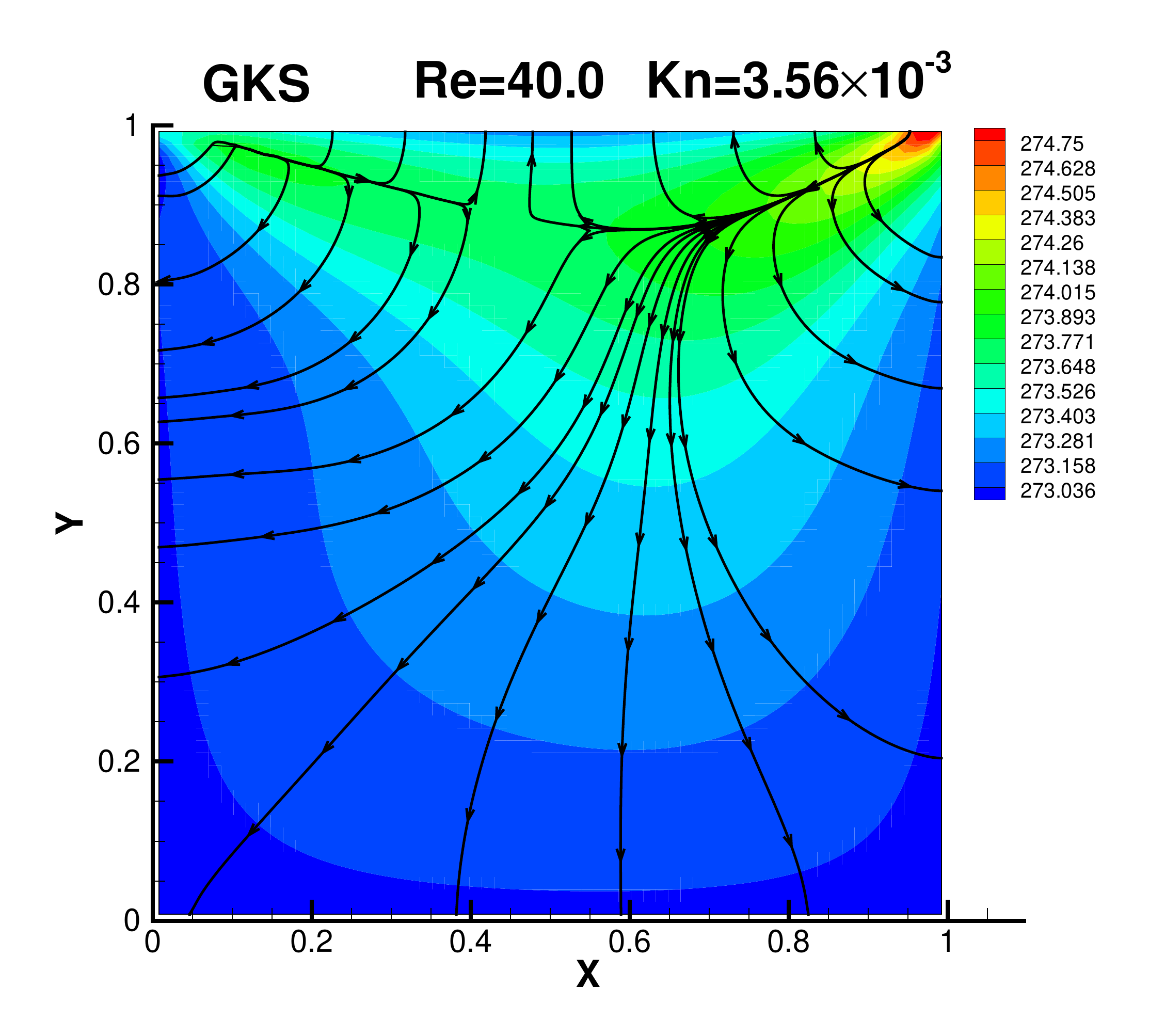}{b}
\includegraphics[width=6cm]{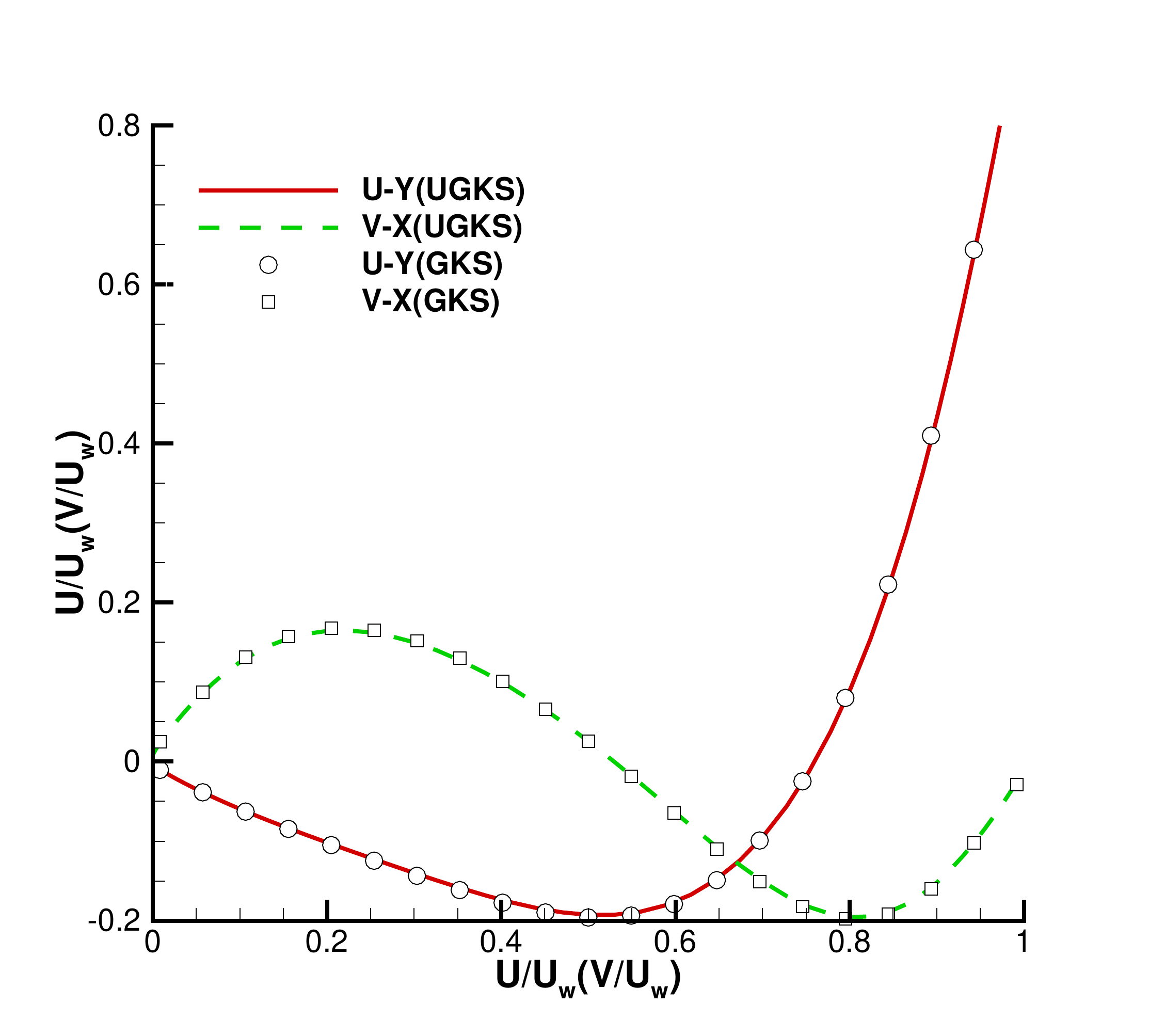}{c}
\caption{Cavity simulation using UGKS and GKS at $\mbox{Kn}=3.56\times10^{-3}$ and $\mbox{Re}=40$.
(a) temperature contour and heat flux: UGKS;
 (b) temperature contour and heat flux: GKS;
 (c) U-velocity along the central vertical line and V-velocity along the central horizontal line, circles: GKS, line: UGKS.}
 \label{r40}
\end{figure}

\begin{figure}
\center
\includegraphics[width=6cm]{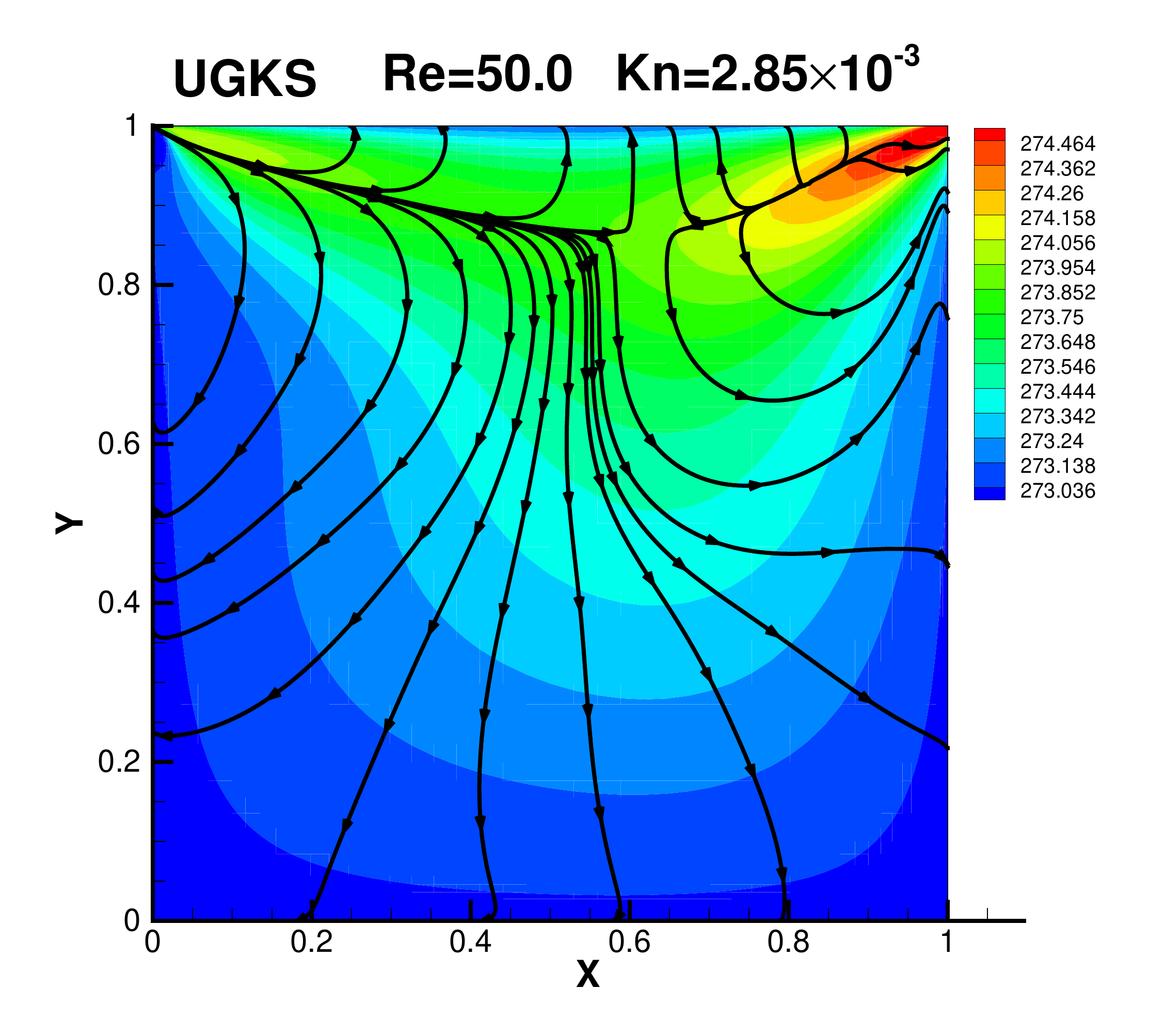}{a}
\includegraphics[width=6cm]{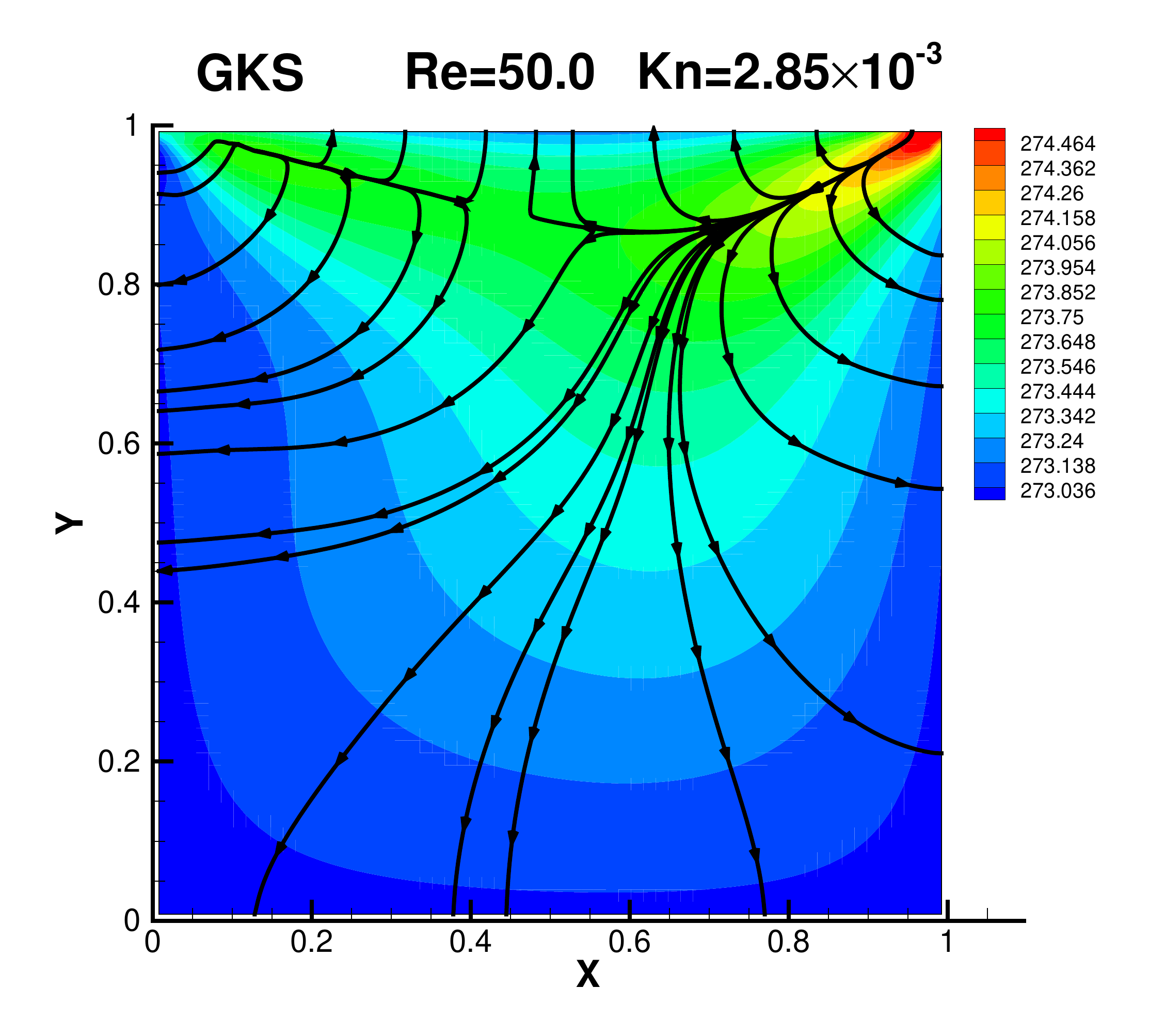}{b}
\includegraphics[width=6cm]{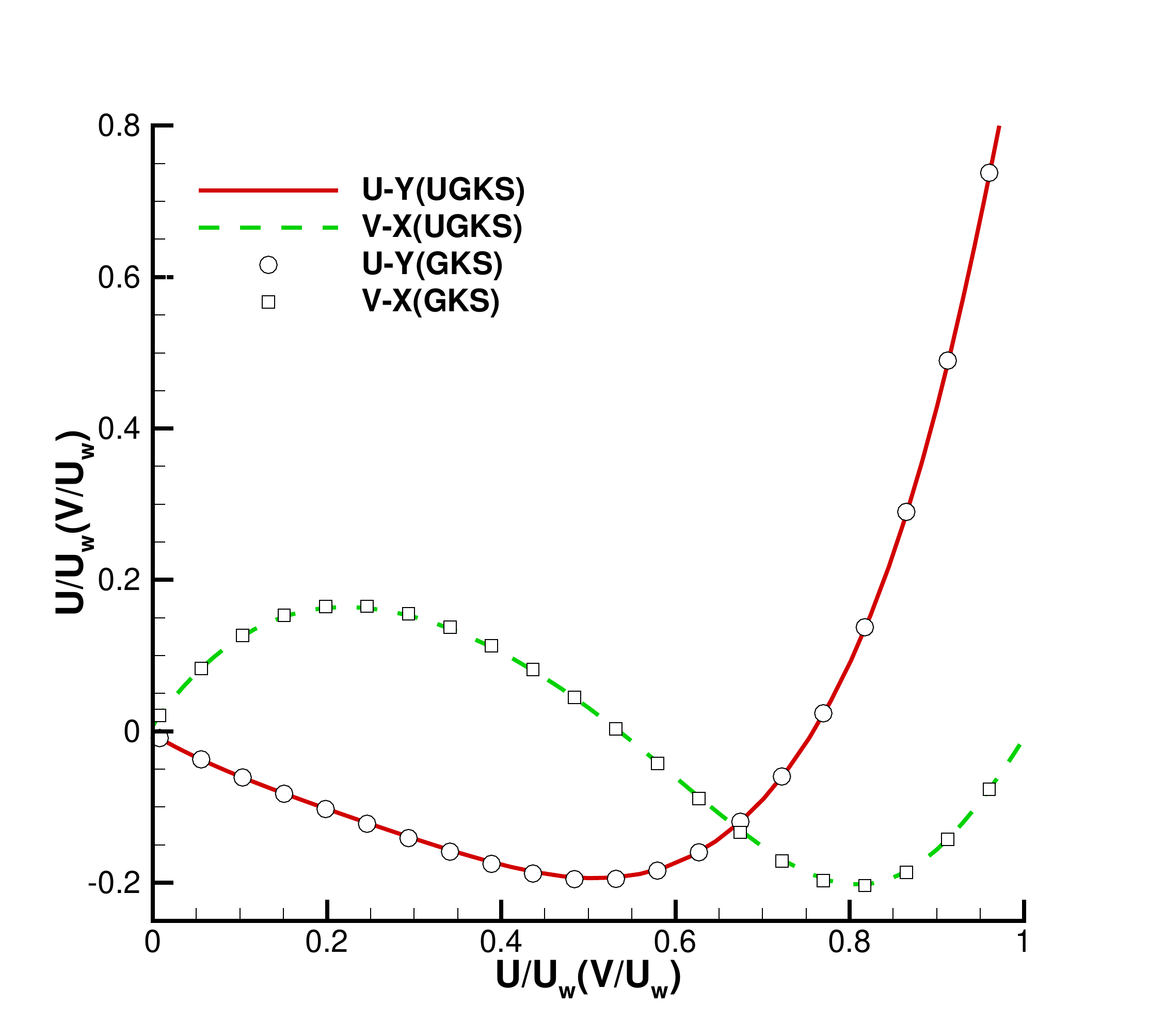}{c}
\caption{Cavity simulation using UGKS and GKS at $\mbox{Kn}=2.85\times10^{-3}$ and $\mbox{Re}=50$.
(a) temperature contour and heat flux: UGKS;
 (b) temperature contour and heat flux: GKS;
 (c) U-velocity along the central vertical line and V-velocity along the central horizontal line, circles: GKS, line: UGKS.}
 \label{r50}
\end{figure}

Based on the above simulations, we get confidence to use the UGKS in the whole flow regime.
In the near continuum regime, it  will be interesting to use UGKS to test the validity of the NS solution.
Before the development of UGKS \cite{xu2010}, an accurate gas-kinetic scheme (GKS) for the NS solutions has been constructed and validated thoroughly \cite{xu2001,xu96}.
The comparison between the solutions from the UGKS and GKS  is basically a comparison of the governing equations of the UGKS and the NS ones.
In the following, we test the cavity case at $\mbox{Re} =5, 10, 20, 30, 40$, and $50$, which are shown
in Fig.(\ref{r5})-(\ref{r50}).
At the above Reynolds numbers, the velocity profiles between UGKS and GKS are basically the same. However, the temperature profiles get close to each other after $\mbox{Re} =20$.
But, the heat flux can keep differences between UGKS and GKS even up to $\mbox{Re} =50$. As shown in these figures,
the heat flux from UGKS is not necessarily perpendicular to the temperature contour level, which is the
basic assumption of the Fourier's law. We believe that the UGKS provides more accurate physical solutions than those from the NS equations.
The UGKS is an indispensable tool in the study of non-equilibrium flow at near continuum flow regime,
which can be used as numerical experiments for the construction of non-equilibrium thermodynamic theory in this regime.

\subsection{Flat-plate boundary layer}

\begin{figure}
\center
\includegraphics[width=10cm]{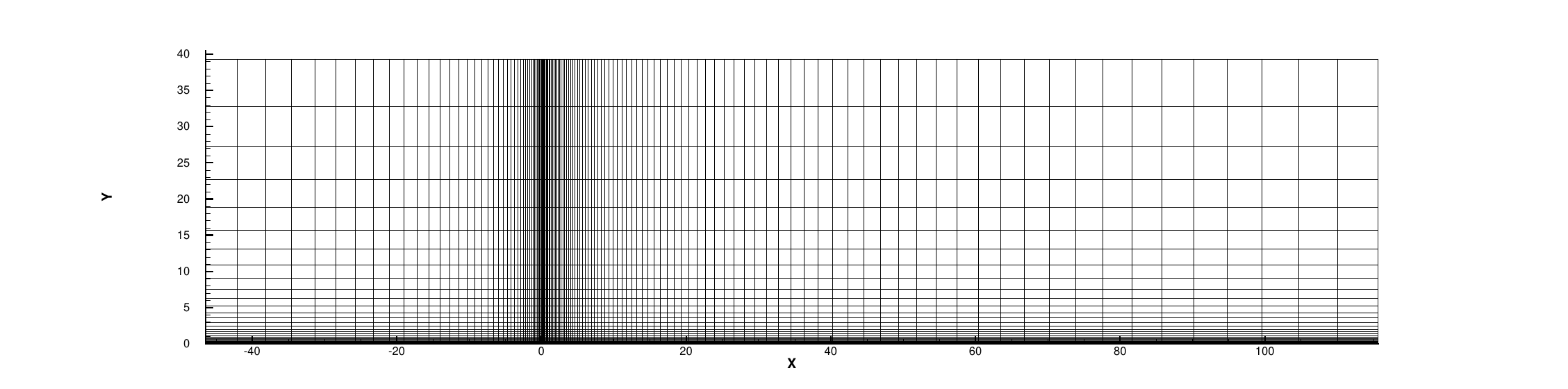}{a}
\includegraphics[width=10cm]{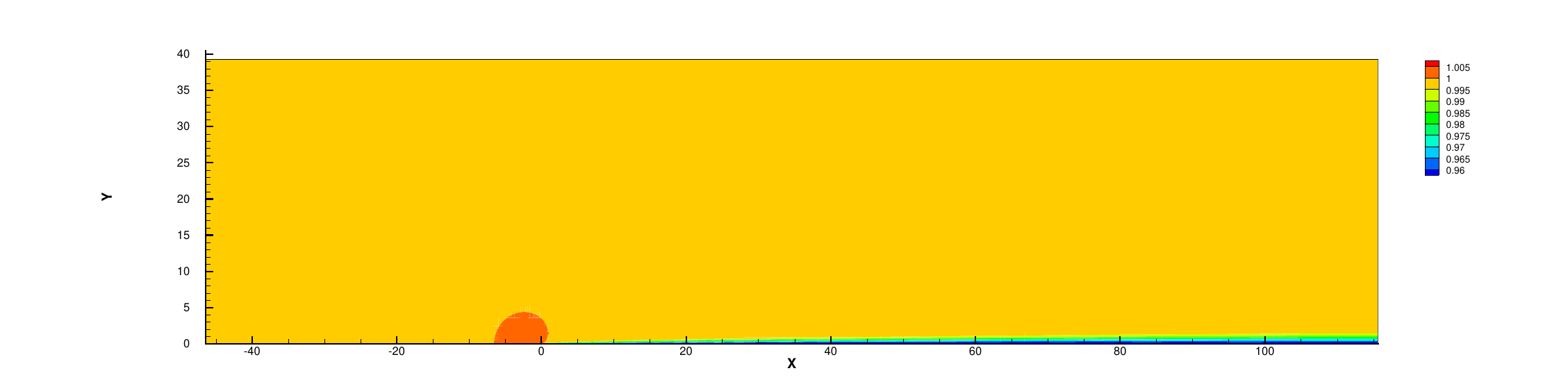}{b}
\includegraphics[width=10cm]{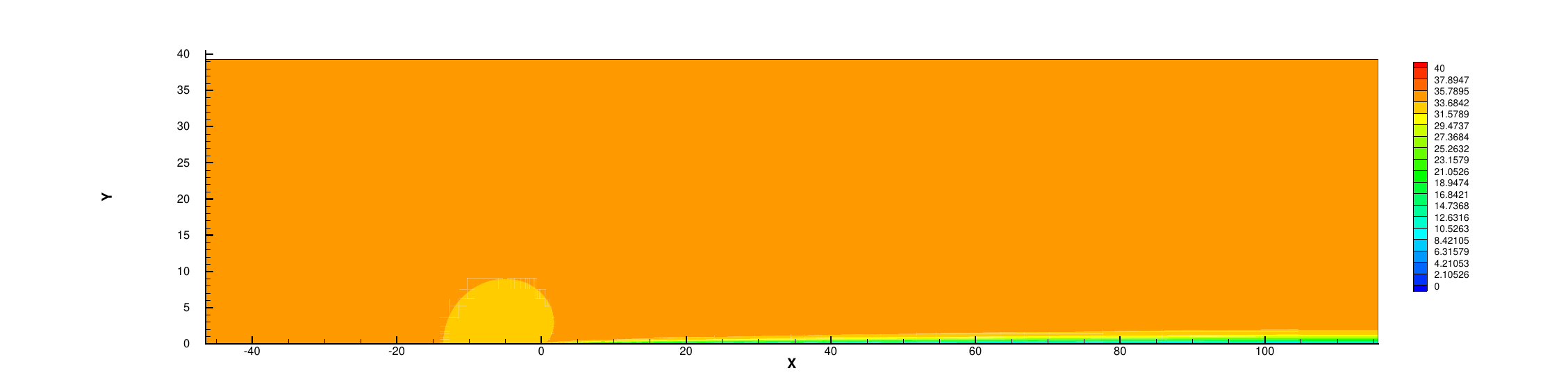}{c}
\includegraphics[width=10cm]{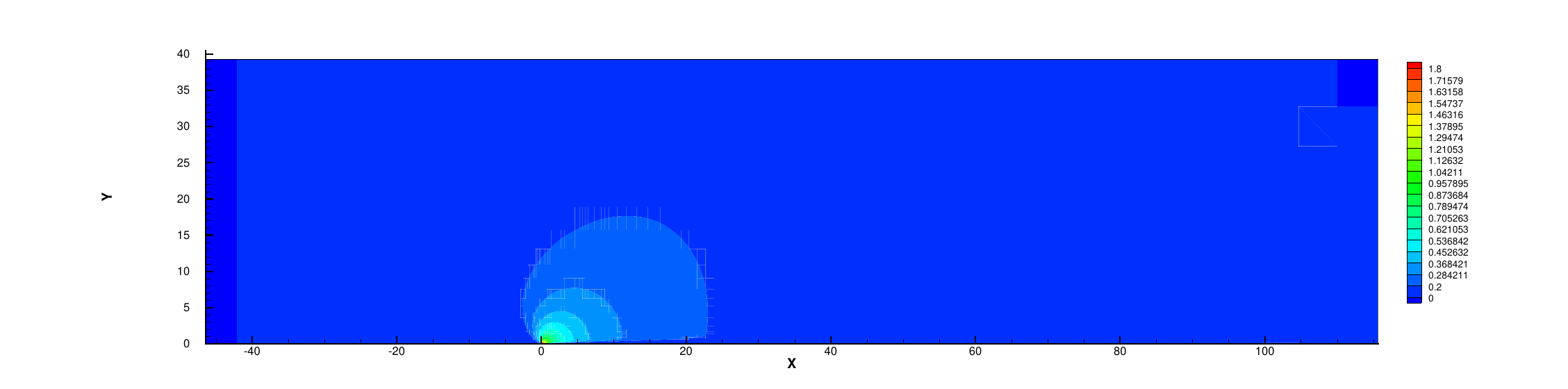}{d}
\caption{Laminar boundary layer computation using UGKS at $\mbox{M}=0.3$ and $\mbox{Re}=10^5$.
(a) mesh distribution; (b) density contours; (c) U velocity contours; (d) V velocity contours.}
\label{layerc}
\end{figure}

\begin{figure}
\center
\includegraphics[width=6cm]{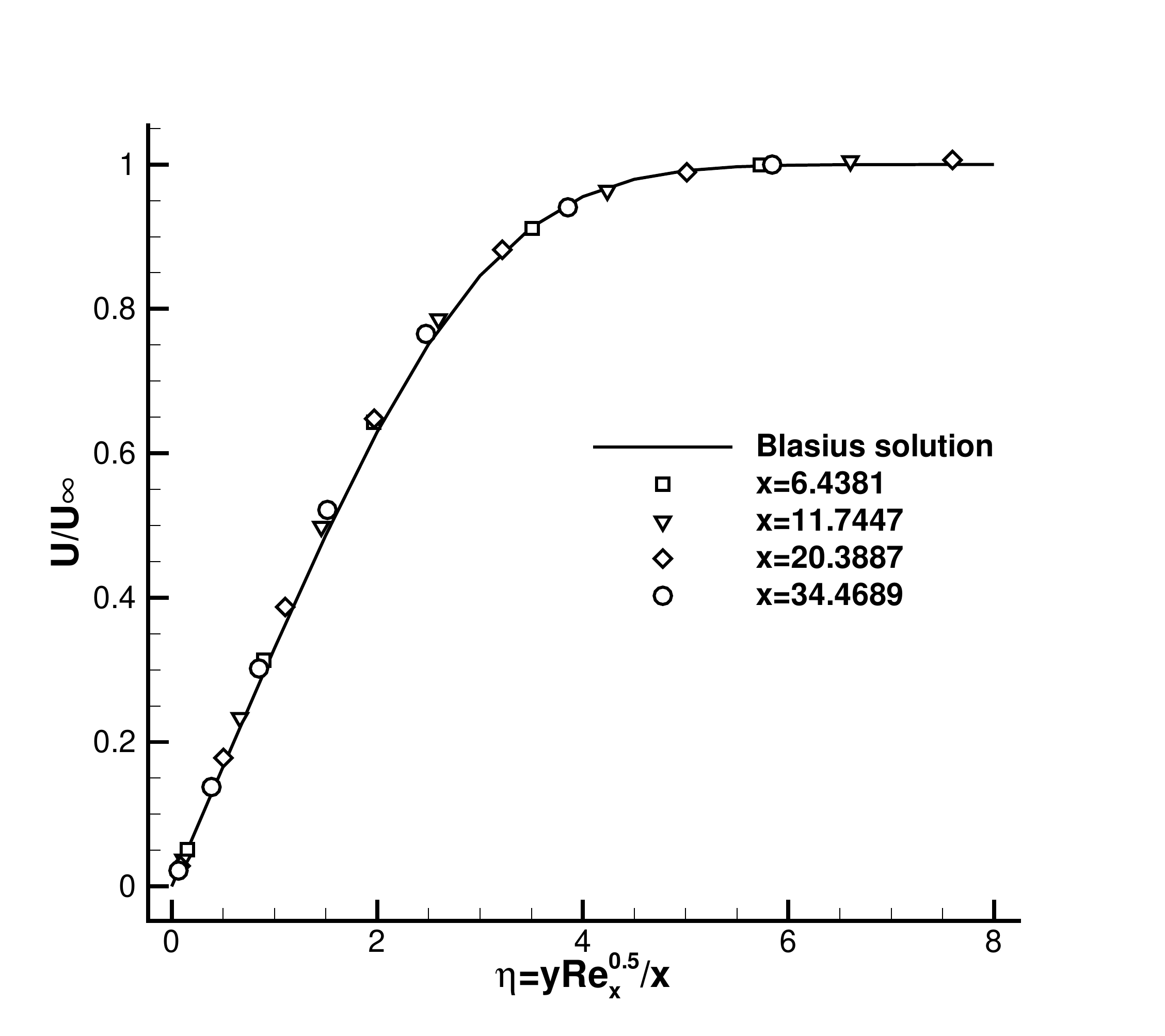}{a}
\includegraphics[width=6cm]{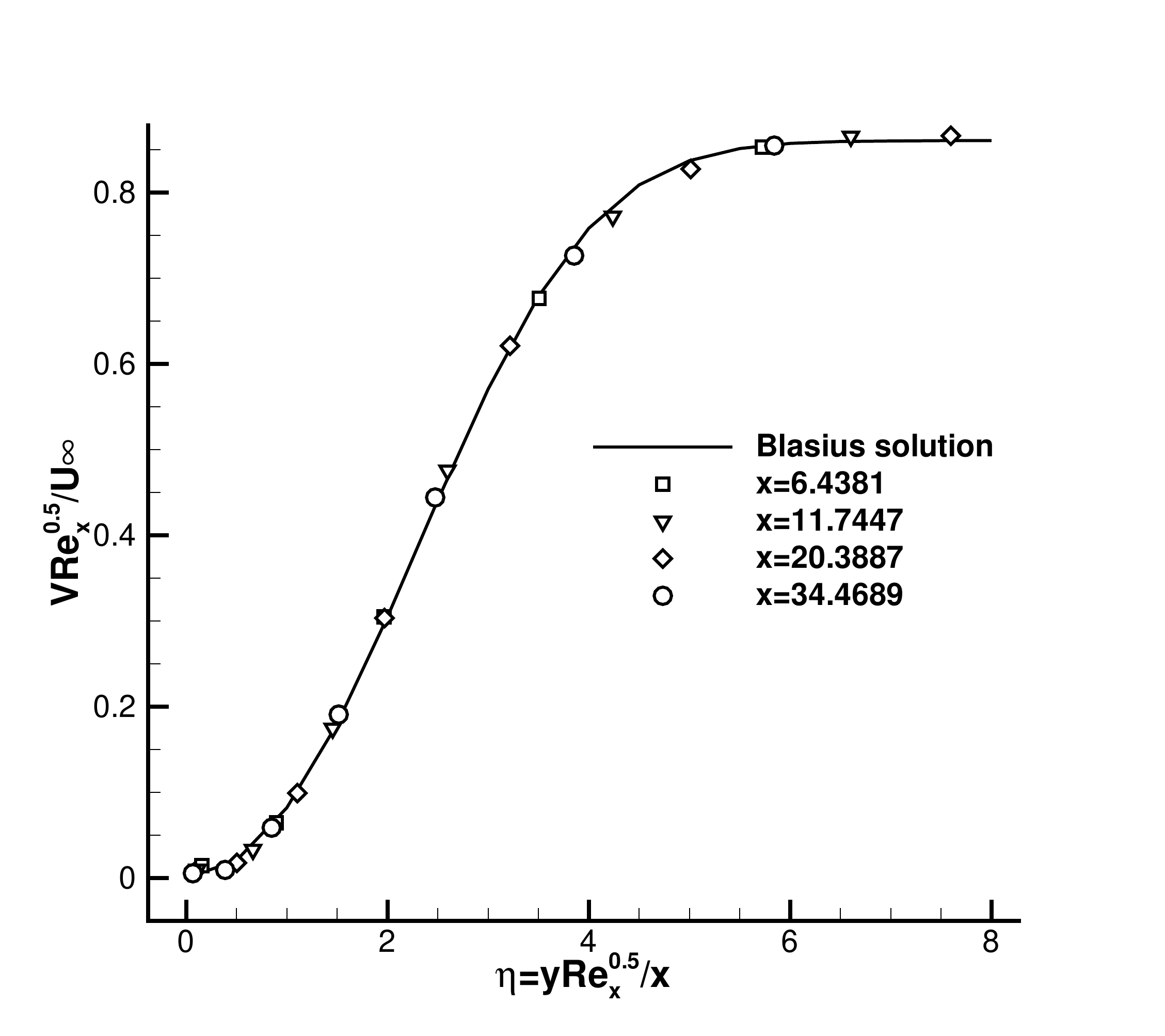}{b}
\caption{UGKS solution. (a) U-velocity distribution at different locations; (b) V-velocity distribution at different locations. Symbols: UGKS, lines: reference solutions.}
\label{layerp}
\end{figure}

The last case is the laminar boundary layer, where the flow is in the fully continuum regime.
The flow at  $\mbox{M}=0.3$ and $\mbox{Re} =10^5$ over a flat plate is simulated. A rectangular mesh with
$120\times30$ non-uniform grid points is used and the mesh distribution is shown in Fig.(\ref{layerc}(a)).
In this case, the local time step is mostly larger than the local particle mean free time and the Shakhov model will be adopted automatically in the current UGKS.
The density, U, and V velocity contours  are shown in Fig.(\ref{layerc}(b))-(\ref{layerc}(d)).
The U and V velocity profiles at different locations are plotted in Fig.(\ref{layerp}),
where the solid lines are the reference Blasius solutions \cite{liao2004explicit}. Even with as less as $5$ mesh points in the boundary layer, both $U$ and $V$ velocity components can be
accurately captured by the UGKS. This is an ideal test case for validating AP properties of kinetic methods.

\section{Conclusion}

In this paper, based on the numerical experiments on the time evolution of a gas distribution function
from the full Boltzmann collision term and the kinetic model equation,
a unified gas-kinetic scheme with the implementation of both collision  models is constructed and tested for multiscale flow problems.
The underlying principle for the development of UGKS is the direct modeling. The modeling scale is the mesh size and time step.
The local flow behavior depends on the ratio of the cell size to the particle mean free path, or the local time step to the particle mean free time.
The principle for the construction of UGKS is different from the traditional CFD methodology,
where fluid dynamic equations are directly discretized.
In UGKS, there is no specific macroscopic or microscopic governing equations to be solved, and the algorithm is based on the modeling of gas evolution in a discretized space directly.
The multiscale nature of UGKS is mainly achieved through the adoption of a time evolution solution for the flux evaluation,
where a dynamic transition from the kinetic scale particle transport to the hydrodynamic scale wave propagation has been incorporated in the flux modeling.
This feature makes UGKS  be able to capture  the corresponding flow physics in different regimes.
With the adoption of a local time step and a switching function between the full Boltzmann collision term and model equation, the UGKS becomes a relatively
efficient method for the study of multiscale flow problems.
In the rarefied flow regime, the UGKS presents the Boltzmann solution, and in hydrodynamic regime it goes to the Navier-Stokes solutions.
In the transition regime, the UGKS itself provides a valuable tool for the study of non-equilibrium flow phenomena.
For example, in the cavity flow simulation the UGKS presents a heat flux which is inconsistent with the Fourier's law even at Reynolds numbers $\mbox{Re} \leq 50$.

Certainly, it is only at the early stage to use the direct modeling concept for the CFD algorithm development \cite{xu-book},
the so-called direct construction of discrete gas evolution equations.
In UGKS, the cell size and time step are not purely related to the truncation errors, but they play dynamic roles. The flow physics described by UGKS depends on the cell resolution.
Therefore, it is hard to use conventional error analysis to evaluate UGKS.
To the current stage, the UGKS has been constructed only for monatomic gas and diatomic gas with rotational degrees of freedom only.
More physical effects, such as
molecular vibration, ionization, even chemical reactions, need to be included into the UGKS.
Mathematical analysis, such as the consistency, stability, and the underlying multiscale governing equations,
needs to be further studied for UGKS.
In conclusion, the UGKS is an extremely useful tool for the study of non-equilibrium thermodynamics.

\section*{Acknowledgement}
The authors would like thank all reviewers for their valuable comments, Prof. S. Jin for his constructive discussion on the stability analysis, and Dr. S.Z. Chen and Dr. L. Wu
for their help on numerical discretization of the full Boltzmann collision term.
The work was supported by Hong Kong research grant council (620813,16211014,16207715) and  NSFC-91330203, and
was partially supported by the open fund of state key laboratory of high-temperature gas dynamics, China (No. 2013KF03).

\clearpage

\clearpage




\end{document}